\pdfoutput=1
\documentclass[12pt,twoside]{article}
\usepackage[utf8]{inputenc}
\usepackage[T1]{fontenc}
\usepackage[english]{babel}
\usepackage{indentfirst}
\usepackage{amsmath,amsthm,amsfonts,amssymb,euscript}
\usepackage{mathtools}
\usepackage{cite}
\usepackage{hyperref}
\hypersetup{
    colorlinks=true,
    linkcolor=blue,
    citecolor=blue,
    urlcolor=blue
}

\newcommand{\Leftclt}{\footnotesize{D.\,L.~Stupin}}

\newcommand{\Rightclt}{\footnotesize{Some properties of extremal functions in the Krzyz problem}}

\theoremstyle{plain}
\newtheorem{Theorem}{Theorem}[section]
\newtheorem{Lemma}{Lemma}[section]
\newtheorem{Corollary}{Corollary}[section]
\newtheorem{Assertion}{Assertion}[section]

\theoremstyle{definition}
\newtheorem{Definition}{Definition}[section]
\newtheorem{Remark}{Remark}[section]
\newtheorem{Example}{Example}[section]
\numberwithin{equation}{section}

\makeatletter
\renewcommand{\@evenhead}{\raisebox{-10pt}[\headheight][0pt]{%
    \vbox{\hbox to\textwidth{\footnotesize\thepage\hfil\Leftclt}\hrule}}}
\renewcommand{\@oddhead}{\raisebox{-10pt}[\headheight][0pt]{%
    \vbox{\hbox to\textwidth{\footnotesize\Rightclt\hfil\thepage}\hrule}}}
\makeatother

\begin{document}
\thispagestyle{empty}
\begin{center}
{\large\textbf{Some properties of extremal functions in the Krzyz problem}}

    \vspace{1ex}

\textbf{Denis L.~STUPIN}\\
    {\footnotesize
Tver State University\\
33 Zhelyabova St., Tver 170100, Russian Federation}
\end{center}%

\vspace{1ex}

\begin{abstract}

    \noindent
\textbf{Abstract.}
Krzyz's conjecture on estimating the moduli of Taylor coefficients in the class of holomorphic, bounded, nonvanishing
functions is considered.
This paper studies both global and locally extremal functions in the introduced
topology, which is stronger than the topology of locally uniform convergence.
The relation between the class under consideration and the Caratheodory class is established, and the general form of
extremal functions in the Krzyz problem is described.
It is shown that every extremal function corresponds to a polynomial $H$ generated by its coefficients and having
positive real part in the unit disk.
These polynomials are studied using a finite-dimensional analogue of the Caratheodory-Toeplitz criterion for
polynomials, following from the Fejer-Riesz theorem.
Necessary extremality conditions are obtained that substantially simplify the study of extremals for fixed $n$.
Conditions for uniqueness of an extremal function are investigated.
In particular, it is proved that uniqueness of a global extremal implies the Krzyz conjecture.
The notion of functions of extremal type, satisfying all the necessary extremality conditions obtained, is introduced.
Every locally extremal function is of extremal type, and the search for such functions at a fixed coefficient index
reduces to solving a finite system of equations and inequalities.
The sets of functions of extremal type are completely described for $n=1,2,3$.
It is proved that the set of functions of extremal type consists of one function for $n=1$, five functions for $n=2$,
and nineteen functions for $n=3$.
Comparing the values of the functional under study on these sets proves the Krzyz conjecture for $n=1,2,3$.

\noindent\textbf{Keywords:}
the Krzyz conjecture, the Krzyz problem, bounded nonvanishing functions,
locally extremal functions, necessary extremality conditions, functions of extremal type,
Caratheodory class, Fejer-Riesz theorem, polynomials with positive real part

\noindent \textbf{Mathematics Subject Classification:} 30C50.

\noindent
This preprint is an English translation of the preprint originally published in Russian on Preprints.ru:
\href{https://doi.org/10.24108/preprints-3113495}{https://doi.org/10.24108/preprints-3113495}.
\end{abstract}
\newpage
\normalsize
\section*{Introduction}

Let us denote the unit disk by $\Delta\coloneqq\{z\in\mathbb{C}\,:\,|z|<1\}$ and the unit circle by $\partial\Delta$.
The Taylor coefficients of a function $f$ will be denoted by $\{f\}_n$, $n\in \{0\}\cup \mathbb N$.

The class $B$ is the set of functions $f$ holomorphic in the unit disk $\Delta$ and satisfying
$0<|f(z)|\leqslant1$, $z\in\Delta$.

In 1968, J.~Krzyz proposed the conjecture~\cite{Krzyz:ref2} that if $f\in B$, then
\begin{equation*}
    |\{f\}_n|\leqslant2/e,\quad n\in\mathbb N,
\end{equation*}
with equality attained only for functions of the form $e^{i\theta}F(e^{i\varphi}z^n,1)$, where
\begin{equation}
    \label{eq:Stupin:BExtrF}
    F(z,t)
    \coloneqq
    e^{-t\frac{1+z}{1-z}},
    \quad\varphi,\theta\in\mathbb R,\quad t\in[0,+\infty).
\end{equation}

Fix $n\in\mathbb N$.
The problem of finding
\begin{equation*}
    m_n \coloneqq \max\limits_{f \in B}{|\{f\}_n|},
\end{equation*}
will be called the Krzyz problem for the coefficient with index $n$ or the Krzyz problem for index $n$.
The problem of finding a sharp estimation of $|\{f\}_n|$ for all $n\in\mathbb N$ on the class $B$ will be called the
Krzyz problem.
The Krzyz problem for fixed $n$ may also be called simply the Krzyz problem if this causes no confusion.
Note that the Krzyz conjecture states that $m_n$ is independent of $n$.

By Montel's compactness principle (see~\cite[p.~201]{Shabat:ref1}), the class $B$ is a compact family.
The existence of extremals in the Krzyz problem for fixed index $n$ is evident, since adjoining the function
$f(z)\equiv0$ to $B$ gives a family $\overline B$ compact in itself in the topology of locally uniform convergence,
while the functional assigning to each function $f$ in $\overline B$ the modulus of its Taylor coefficient with index
$n$ is continuous on $\overline B$ and therefore attains its supremum (see~\cite[p.~203]{Shabat:ref1}).

A function $f\in B$ will be called globally extremal in the Krzyz problem for index $n$ if $|\{f\}_n|=m_n$.
The notion of a locally extremal function will be defined precisely in Section~\ref{sec:Stupin:theTopology}, using the
topology introduced there.
Throughout the paper, the term extremal function (or extremal) means a locally extremal function in this topology unless
global extremality is stated explicitly.
Every globally extremal function is locally extremal.

If $f(z)=F(-z^n,1)$, then $\{f\}_n=2/e$, and hence $m_n\geqslant2/e$.
At present, no functions $f$ with $\{f\}_n>2/e$ are known to us.
If such a function existed, the Krzyz conjecture would be false.

\section*{Purpose and relevance of the paper}

The problem of obtaining an elementary proof of the Krzyz conjecture remains relevant.
One of the main obstacles to such a proof is the large number of locally extremal functions that must be found and
studied~\cite{Prokhorov:ref1},~\cite{Stupin:ref15}.
Unlike many papers in this area, we study both global and local extremals.
Extremal functions are difficult to find directly, so it is important to identify a class of functions that are easier
to construct and satisfy all known necessary extremality conditions.

The purpose of this paper is to study necessary extremality conditions using a polynomial $H$ with positive real part in
the unit disk, generated by the coefficients of a function $f$, and on this basis to develop a method for constructing
functions of extremal type.
Passing from the function $f$ to the polynomial $H$ reduces this problem to a finite-dimensional one.
The results obtained are based mainly on a characterization of the set of such polynomials, which in turn follows from
the classical Fejer-Riesz theorem for trigonometric polynomials.

In this paper, we present and generalize some known necessary extremality conditions in the Krzyz problem and obtain new
results in this area.
We review these conditions and show how they are related.
Using a result on uniqueness of a polynomial with positive real part, we obtain conditions for uniqueness of an extremal
function in the Krzyz problem and show that uniqueness of a global extremal implies the Krzyz conjecture.
On the basis of the conditions considered, we formulate a method for constructing functions satisfying all these
conditions.
We call such functions functions of extremal type; all known necessary extremality conditions carry over to them without
change.
For every fixed $n$, the search for these functions reduces to solving a finite system of equations and inequalities.
The proposed method makes it possible to list all functions of extremal type for a given $n$.
It is used to prove the Krzyz conjecture for $n=1,2,3$.

All results are presented in a complete and self-contained form; references are provided for further reading and
attribution.

\section{Subordinate functions}

By $\Omega$ we denote the class of functions $\omega$ holomorphic in $\Delta$ and satisfying $|\omega(z)|\leqslant1$,
$z\in\Delta$, and by $\Omega_0$ the subclass of $\Omega$ normalized by $\omega(0)=0$.

Let the functions $G$ and $g$ be holomorphic in $\Delta$.
A function $g$ is called subordinate to a function $G$ in $\Delta$ if it can be represented in $\Delta$ as
$g(z)=G(\omega(z))$, where $\omega\in\Omega_0$.
The function $G(z)$ will be called a majorant for $g(z)$ in $\Delta$.
Clearly, $\{g\}_0=\{G\}_0$.

\smallskip
The notion of subordination goes back to E.~Lindelof~\cite{Lindelof:ref1}.
The modern terminology was introduced later by J.~E.~Littlewood~\cite{Littlewood:ref1} and
W.~Rogosinski~\cite{Rogosinski:ref1}, who also developed the method and obtained several results by means of it.
The Littlewood-Rogosinski subordination principle is often used to derive coefficient estimations in the class $B$
(see~\cite{Hummel:ref1, Szapiel:ref1, Stupin:ref13, Peretz:ref1, Stupin:ref5, Stupin:ref6}).

In the Krzyz problem, the difficulty in applying this method lies in the complexity of the coefficients $\{F\}_k(t)$ of
the function $F(z,t)$.

Subordination theory makes it very easy to find sharp estimations of the first and second coefficients on the class of
functions $g(z)$ subordinate to a function $G(z)$.
Let $M_G$ be the class consisting of functions $g(z)=G(\omega(z))$, where $G$ is holomorphic in $\Delta$ and
$\omega\in\Omega_0$.

We state Schwarz's lemma~\cite[p.~29]{Golusin:ref1} in the following differential form:

\begin{Lemma}[Schwarz]
    \label{lem:Stupin:Schwarz}
If $\omega\in\Omega_0$, then $|\{\omega\}_1|\leqslant1$, with equality attained only by rotations of the function
$\omega(z)=z$ in the variable plane $z$ or $w$
    \textnormal($w=\omega(z)$\textnormal).
\end{Lemma}

\begin{Theorem}[Rogosinski]
If $g\in M_G$, then
    \begin{equation}
        \label{eq:Stupin:f_1}
        |\{g\}_1|\leqslant|\{G\}_1|,
    \end{equation}
and estimation~\eqref{eq:Stupin:f_1} is sharp.
If $\{G\}_1\neq0$, then equality in inequality~\eqref{eq:Stupin:f_1} is attained only by rotations of the function
$G(z)$ in the variable plane $z$.
If $\{G\}_1=0$, then $\{g\}_1=0$, and equality is attained by every function $g\in M_G$.
\end{Theorem}

\begin{proof}
Let $g\in M_G$, that is, there exists a function $\omega\in\Omega_0$ such that $g(z)=G(\omega(z))$.
Then $\{g\}_1=\{G\}_1\{\omega\}_1$.
Therefore, $|\{g\}_1|=|\{G\}_1||\{\omega\}_1|$.
Hence, $|\{g\}_1|\leqslant|\{G\}_1|$ by Schwarz's lemma (Lemma~\ref{lem:Stupin:Schwarz}).
If $\{G\}_1\neq0$, the equality cases follow from the same lemma.
If $\{G\}_1=0$, then $\{g\}_1=0$ for every $g\in M_G$.
\end{proof}

This assertion is taken from~\cite[p.~50]{Rogosinski:ref1}, but the degenerate case $\{G\}_1=0$ is omitted
in~\cite{Rogosinski:ref1}.
We now state the form of Schwarz's lemma~\cite[p.~319]{Golusin:ref1} invariant under automorphisms of the unit disk in
the following differential form:

\begin{Lemma}[Schwarz-Pick]
    \label{lem:Stupin:Schwarz-Pick}
If $\omega\in\Omega$, then $|\{\omega\}_1|\leqslant1-|\{\omega\}_0|^2$, with equality attained only by rotations of the
function $\omega(z)=\frac{\{\omega\}_0+z}{1+\overline{\{\omega\}_0}z}$ in the variable plane $z$ or $w$
    \textnormal($w=\omega(z)$\textnormal).
\end{Lemma}

Clearly, for $\{\omega\}_0=0$, the Schwarz-Pick lemma reduces to Schwarz's lemma.
The following assertion follows directly from Lemma~\ref{lem:Stupin:Schwarz-Pick}:

\begin{Corollary}
    \label{cor:Stupin:Schwarz-Pick_2}
If $\omega\in\Omega_0$, then $|\{\omega\}_2|\leqslant1-|\{\omega\}_1|^2$, with equality in this inequality attained only
by rotations of the function $\omega(z)=z\frac{\{\omega\}_1+z}{1+\overline{\{\omega\}_1}z}$ in the variable plane $z$ or
$w$ \textnormal($w=\omega(z)$\textnormal).
\end{Corollary}

\begin{proof}
Set $\varphi(z)\coloneqq\omega(z)/z$ for $z\neq0$ and $\varphi(0)\coloneqq\{\omega\}_1$.
By the maximum principle for holomorphic functions, $\varphi\in\Omega$.
Note that $\varphi(0)=\{\omega\}_1$ and $\varphi'(0)=\{\omega\}_2$.
Applying Lemma~\ref{lem:Stupin:Schwarz-Pick} to the function $\varphi$ gives the required estimation.
The equality conditions in this lemma give rotations of the function
$\varphi(z)=\frac{\{\omega\}_1+z}{1+\overline{\{\omega\}_1}z}$ in the variable plane $z$ or $w$, which implies the
stated form of the function $\omega$.
\end{proof}

The following assertion is taken in part from~\cite[p.~59]{Rogosinski:ref1}:

\begin{Theorem}[Rogosinski]
If $g\in M_G$, then
    \begin{equation}
        \label{eq:Stupin:f_2}
        |\{g\}_2|\leqslant\max(|\{G\}_1|,|\{G\}_2|),
    \end{equation}
and estimation~\eqref{eq:Stupin:f_2} is sharp.
Equality in inequality~\eqref{eq:Stupin:f_2} is attained in the following and only the following cases\textnormal:
    \newline\textnormal1.
If $|\{G\}_1|>|\{G\}_2|$, then $g(z)=G(\eta z^2)$, where $|\eta|=1$\textnormal;
    \newline\textnormal2.
If $|\{G\}_1|<|\{G\}_2|$, then $g(z)=G(\eta z)$, where $|\eta|=1$\textnormal;
    \newline\textnormal3.
If $|\{G\}_1|=|\{G\}_2|=0$, then $g$ is any function in $M_G$\textnormal;
    \newline\textnormal4.
If $|\{G\}_1|=|\{G\}_2|\neq0$, then $g(z)=G(\omega_{r,\eta}(z))$,
        $\omega_{r,\eta}(z)
        \coloneqq
        \eta z\frac{r+\lambda\eta z}{1+r\lambda\eta z}$,
        $\lambda\coloneqq\frac{\{G\}_2}{\{G\}_1}$,
where $0\leqslant r\leqslant1$ and $|\eta|=1$.
In particular, for $r=0$, $g(z)=G(\lambda\eta^2 z^2)$, and for $r=1$, $g(z)=G(\eta z)$.
\end{Theorem}

\begin{proof}
Let $\omega(z)=\{\omega\}_1 z+\{\omega\}_2 z^2+o(z^2)\in\Omega_0$.
Since
    \begin{equation}
        \label{eq:Stupin:g_2}
        \{g\}_2=\{G\}_1\{\omega\}_2+\{G\}_2\{\omega\}_1^2,
    \end{equation}
by Corollary~\ref{cor:Stupin:Schwarz-Pick_2} and Schwarz's lemma~\ref{lem:Stupin:Schwarz}, we obtain
    \[
        |\{g\}_2|
        \leqslant
        |\{G\}_1|(1-|\{\omega\}_1|^2)+|\{G\}_2||\{\omega\}_1|^2
        \leqslant
        \max(|\{G\}_1|,|\{G\}_2|).
    \]

We now consider the equality cases.

    \noindent
    \textbf{1.}
If $|\{G\}_1|>|\{G\}_2|$, equality is possible only if $\{\omega\}_1=0$ and $|\{\omega\}_2|=1$.
The equality condition in Corollary~\ref{cor:Stupin:Schwarz-Pick_2} implies that $\omega(z)=\eta z^2$, where $|\eta|=1$.

    \noindent
    \textbf{2.}
If $|\{G\}_1|<|\{G\}_2|$, equality is possible only if $|\{\omega\}_1|=1$ and $|\{\omega\}_2|=0$.
The equality condition in Schwarz's lemma (Lemma~\ref{lem:Stupin:Schwarz}) implies that $\omega(z)=\eta z$, where
$|\eta|=1$.

    \noindent
    \textbf{3.}
If $|\{G\}_1|=|\{G\}_2|=0$, then $\{g\}_1=\{g\}_2=0$ for every function $g\in M_G$.

    \noindent
    \textbf{4.}
Now let $|\{G\}_1|=|\{G\}_2|\neq0$.
By Schwarz's lemma (Lemma~\ref{lem:Stupin:Schwarz}), $\{\omega\}_1=r\eta$, where $r\coloneqq|\{\omega\}_1|\in[0,1]$,
$|\eta|=1$.
By formula~\eqref{eq:Stupin:g_2}, equality in inequality~\eqref{eq:Stupin:f_2} is possible if and only if
    \[
        |\{\omega\}_2|=1-|\{\omega\}_1|^2
        \wedge
        \big(
        \big(
        r\in(0,1)
        \wedge
        \arg(\{G\}_1\{\omega\}_2)=\arg(\{G\}_2\{\omega\}_1^2)
        \big)
        \vee
        r=0
        \vee
        r=1
        \big).
    \]

Let $r\in(0,1)$.
Rewrite the second equality as $\arg(\{\omega\}_2)=\arg(\lambda\{\omega\}_1^2)$,
    $\lambda\coloneqq\frac{\{G\}_2}{\{G\}_1}$.
It follows from $|\{G\}_1|=|\{G\}_2|\neq0$ that $|\lambda|=1$.
It follows from $\{\omega\}_1=r\eta$ and $\arg(\{\omega\}_2)=\arg(\lambda\{\omega\}_1^2)$ that
    $\{\omega\}_2=(1-r^2)\lambda\eta^2$.
The equality condition in Corollary~\ref{cor:Stupin:Schwarz-Pick_2} implies that
    \[
        \omega(z)
        =
        e^{i\tau}e^{i\theta}z\frac{\{\omega\}_1+e^{i\theta}z}{1+\overline{\{\omega\}_1}e^{i\theta}z},
        \quad
        \tau,\theta\in\mathbb R.
    \]
Comparing the coefficients of $z$ in this formula, we obtain
    $e^{i\tau}e^{i\theta}\{\omega\}_1=\{\omega\}_1$.
Since $\{\omega\}_1\neq0$, we have $e^{i\tau}=e^{-i\theta}$.
Moreover,
    $
        \{\omega\}_2=e^{i\theta}(1-|\{\omega\}_1|^2)=e^{i\theta}(1-r^2)
    $,
but $\{\omega\}_2=(1-r^2)\lambda\eta^2$, so $e^{i\theta}=\lambda\eta^2$.
Hence,
    \[
        \omega(z)
        =
        z\frac{r\eta+\lambda\eta^2 z}{1+r\overline{\eta}\lambda\eta^2 z}
        =
        z\frac{r\eta+\lambda\eta^2 z}{1+r\lambda\eta z}
        =
        \eta z\frac{r+\lambda\eta z}{1+r\lambda\eta z}.
    \]

If $|\{\omega\}_1|=0$, then $|\{\omega\}_2|=1$ and $\omega(z)=\xi z^2$, $|\xi|=1$, by
Corollary~\ref{cor:Stupin:Schwarz-Pick_2}.
Since $|\xi|=|\eta|=|\lambda|=1$, we can always choose $\eta$ so that $\xi=\lambda\eta^2$.

If $|\{\omega\}_1|=1$, then $|\{\omega\}_2|=0$ and $\omega(z)=\eta z$ by Schwarz's lemma~\ref{lem:Stupin:Schwarz}.

Conversely, each such function clearly belongs to $\Omega_0$ and gives equality in inequality~\eqref{eq:Stupin:f_2}.
\end{proof}

Although a general method for describing equality cases is given in~\cite[p.~57]{Rogosinski:ref1}, the equality cases in
inequality~\eqref{eq:Stupin:f_2} are not listed in~\cite{Rogosinski:ref1}.

\section{Subclasses}

The set of functions $h$ holomorphic in the disk $\Delta$, normalized by $h(0)>0$ and having positive real part there,
will be denoted by $C$ and called the Caratheodory class.
Fix $t>0$.
The set of functions $h$ in $C$ normalized by $h(0)=t$ will be denoted by $C_t$ and called the normalized Caratheodory
class, or simply the Caratheodory class if this causes no confusion.
Note that $C=\bigcup\limits_{t>0}C_t$.
The subsets of polynomials of degree $n$ in the classes $C$ and $C_t$ will be denoted by $C^n$ and $C_t^n$,
respectively.

Since the set $B$ is invariant under rotations in the variable plane $w$ ($w=f(z)$), we may restrict ourselves without
loss of generality to functions for which $f(0)>0$.
Since $0<\{f\}_0\leqslant1$, we can set $\{f\}_0=e^{-t}$, where the parameter $t\in[0,+\infty)$.
The corresponding subclasses will be denoted by $B_t$.
As is known from subordination theory~\cite{Rogosinski:ref1}, every function in the class $B_t$ can be represented as
\begin{equation}
    \label{eq:Stupin:B_t_to_C}
    f(z)
    =
    e^{-t\,h(z)},
    \quad
    h\in C_1.
\end{equation}
Note that for every $t>0$ this formula establishes a one-to-one correspondence between the classes $C_1$ and $B_t$.
It is also clear that for every fixed $t>0$ the function $F(z,t)$ is a majorant for every function in the class
$B_t$.

Note also that the class $B_0$ consists of a single function $f\equiv1$, so $B_0$ may be regarded as completely studied.
Below we may state $t\geqslant0$ for completeness, but in fact $t>0$ may be assumed throughout.
This convention allows us, for example, to divide freely by $t$.

The class of functions $h\in C_t$ with real coefficients will be denoted by $C_t^r$, and the class of functions $f\in
B_t$ with real coefficients by $B_t^r$.
For every $t>0$, formula~\eqref{eq:Stupin:B_t_to_C} establishes a one-to-one correspondence between the classes $C_1^r$
and $B_t^r$.
Similarly, define the classes $B^r$, $C^r$, and $\Omega_0^r$ as the subclasses of $B$, $C$, and $\Omega_0$ consisting of
functions with real coefficients.

\section{A brief review of results on the Krzyz conjecture}
\label{sec:Stupin:review}

At present, the Krzyz conjecture has been proved for the first six Taylor coefficients inclusive.
Geometric considerations clearly give $|\{f\}_0|\leqslant1$.
A sharp estimation of $|\{f\}_1|$ can be found in many papers starting in 1934, the first being~\cite{Levin:ref1}.
Estimating $|\{f\}_2|$ has presented no difficulty since 1943 thanks to subordination theory~\cite{Rogosinski:ref1}.
Krzyz formulated the conjecture discussed here on the basis of estimations of the moduli of the first two coefficients.
The proof for $n=3$ was first published in 1977 by J.~Hummel, S.~Scheinberg, and L.~Zalcman~\cite{Hummel:ref1}.
An estimation of the functional $|\{f\}_3|$ on the class $B_t$, sharp for every $t>0$, was obtained
in~\cite{Prokhorov:ref1}.
An analogous estimation on the set of functions in $B_t$ with real coefficients is given in~\cite{Stupin:ref15}.
For $n=4$, we should mention the proof by W.~Szapiel~\cite{Szapiel:ref1}.
The fifth coefficient was first estimated by Szapiel's method in N.~Samaris's paper~\cite{Samaris:ref1}.
Using Szapiel's method, the author of this paper obtained the estimation $|\{f\}_6|\leqslant2/e+0.00116077$
in~\cite{Stupin:ref13}, and the estimation $|\{f\}_6|\leqslant2/e$ by a numerical method in~\cite{Stupin:ref10}.

Research on the Krzyz problem can be divided into several main areas.
One is the estimation of initial coefficients mentioned in the preceding paragraph.
The so-called asymptotic estimations, discussed for example in~\cite{Peretz:ref1, Stupin:ref14}, deserve separate
attention.
There are also estimations uniform in $n$ obtained using Cauchy's integral formula for all positive integers
$n$.
The paper~\cite{HorowitzC_1978} obtained the estimation
$|\{f\}_n|\leqslant1-\frac1{3\pi}+\frac4{\pi}\sin\frac1{12}=0.999877\ldots$,
and the dissertation~\cite[p.~19]{Ermers:ref1} the improved estimation
$|\{f\}_n|\leqslant\frac45+\frac4{\pi}\sin\frac{\pi}{20}=0.999178\ldots$
Another important area concerns the properties of extremal functions and the search for functions in the class $B$
having these properties.
This area includes, in particular,~\cite{Hummel:ref1, Peretz:ref2, Martin:ref1} and the present paper.

Problems in geometric function theory related in one way or another to the Krzyz conjecture are considered
in~\cite{KortramRA_1993, Szapiel:ref1, SzynalJ_2001, SzynalJ_2003, PeretzR_2002, Peretz:ref2, Peretz:ref3, Agler:ref1}.
Some generalizations of the Krzyz conjecture are described in~\cite[p.~187]{Hummel:ref1}.

In~\cite{Krushkal:ref1}, S.~L.~Krushkal proposed a proof of the Hummel-Scheinberg-Zalcman
conjecture~\cite[p.~189]{Hummel:ref1} for functions in the spaces $H^p$, $p=2m$.
This proof uses a new method involving deep properties of Teichmuller spaces.
The Krzyz conjecture follows from the Hummel-Scheinberg-Zalcman conjecture by passing to the limit as $p\to\infty$.

The present paper is devoted to the search for an elementary proof of the Krzyz conjecture.

\section{The Caratheodory-Toeplitz criterion}

For further reference, we state the following classical result, which can be found
in~\cite{Caratheodory:ref2, ToeplitzO_1911, Stupin:ref1}:
\begin{Theorem}[Caratheodory, Toeplitz]
    \label{th:Stupin:nCTEquiv}
Let $n\in\mathbb{N}$, $\{h\}_0>0$, $\{h\}_1,\ldots,\{h\}_{n}\in\mathbb{C}$.
The polynomial
    \begin{equation}
        \label{eq:Stupin:C-T_polynom}
        Q_n(z)
        \coloneqq
        \{h\}_0+\sum\limits_{k=1}^{n}{\{h\}_k z^k}
    \end{equation}
can be extended to a function
    \[
        h(z)
        \coloneqq
        Q_n(z)+o(z^{n})
        ~
        \in
        ~
        C
    \]
if and only if the determinants
    \begin{equation}
        \label{eq:Stupin:C-T_minors}
        M_k
        \coloneqq
        \begin{vmatrix}
            2\{h\}_0               & \{h\}_1                & \cdots & \{h\}_{k-1}        & \{h\}_{k}\\
            \overline{\{h\}_1}     & 2\{h\}_0               & \cdots & \{h\}_{k-2}        & \{h\}_{k-1}\\
            \vdots                 & \vdots                 & \ddots & \vdots             & \vdots\\
            \overline{\{h\}_{k-1}} & \overline{\{h\}_{k-2}} & \cdots & 2\{h\}_0           & \{h\}_1\\
            \overline{\{h\}_{k}}   & \overline{\{h\}_{k-1}} & \cdots & \overline{\{h\}_1} & 2\{h\}_0\\
        \end{vmatrix}
        ,\quad k=0,\ldots,n,
    \end{equation}
are either all positive or positive up to some index $m\leqslant n$, starting from which they all vanish.
In the latter case the extension is unique, and there exist numbers $\alpha_j>0$, $j=1,\ldots,m$, and numbers
$z_j\coloneqq e^{i\varphi_j}$, $j=1,\ldots,m$, $0\leqslant\varphi_1<\ldots<\varphi_m<2\pi$ such that
    \begin{equation}
        \label{eq:Stupin:B_t_extremal_h}
        h(z)
        =
        \sum\limits_{j=1}^m \alpha_j\frac{1+z_j z}{1-z_j z}.
    \end{equation}
\end{Theorem}

Fix $n\in\mathbb N$ and consider the space $\mathbb R\times\mathbb C^n$ whose points are vectors with one real and $n$
complex coordinates.
For every function $h\in C$, set
$
h^{(n+1)}\coloneqq (\{h\}_0,\ldots,\{h\}_{n}).
$
The set of all such points $h^{(n+1)}$ in the space $\mathbb R\times\mathbb C^n$ will be denoted by $C^{(n+1)}$ and
called the $(n+1)$th coefficient body of the class $C$.

The coefficient problem for the class $C$ is formulated as follows: find necessary and sufficient conditions on complex
numbers $\{h\}_0,\{h\}_1,\{h\}_2,\ldots$ for the series $\{h\}_0+\{h\}_1 z+\{h\}_2 z^2+\ldots$ to be the Taylor series
of a function in the class $C$.
The Caratheodory-Toeplitz criterion gives a complete solution to this problem.

\smallskip
Identifying the space $\mathbb R\times\mathbb C^n$ with $\mathbb R^{2n+1}$, we equip it with the usual Euclidean
topology and use the corresponding notions of interior and boundary points.
Using these considerations, one can prove the following theorem (see~\cite[p.~480]{Golusin:ref1}):

\begin{Theorem}[Caratheodory, Fejer]
    \label{th:Stupin:CaratheodoryFejer}
Let $n\in\mathbb N$.
For every function $h$ in the class $C$, there exists a unique fraction of the form
    \begin{equation*}
        h_m(z)
        \coloneqq
        \lambda+\sum\limits_{j=1}^{m}\alpha_j\dfrac{1+z_j z}{1-z_j z},
    \end{equation*}
where
    $m\in\{0,\ldots,n\}$,
    $0\leqslant\lambda\leqslant\{h\}_0$,
    $\alpha_j>0$,
    $\lambda+\sum\limits_{j=1}^{m}\alpha_j = \{h\}_0$,
    $z_j\coloneqq e^{i\varphi_j}$, $j=1,\ldots,m$,
    $0\leqslant\varphi_1<\ldots<\varphi_m<2\pi$,
whose first $n+1$ Taylor coefficients are respectively equal to $\{h\}_0,\ldots,\{h\}_n$.
\end{Theorem}

C.~Caratheodory's paper~\cite{Caratheodory:ref2}, published in 1911, is regarded as the foundational work on the
coefficient problem for the class $C$, although Caratheodory published his first investigations on this subject in 1907
in~\cite{CaratheodoryC_1907}.
Note that the class $C$ is denoted by the first letter of the surname Carath\'eodory and named the Caratheodory class in
his honor; in the literature on univalent functions, the set $C_1$ is often called the Caratheodory class.

Caratheodory established a parametric representation of the boundary of the coefficient body $C^{(n+1)}$ that completely
describes this boundary by a trigonometric polynomial.
In 1911, O.~Toeplitz found~\cite{ToeplitzO_1911} a simple and very beautiful determinant representation of this
boundary, which reduced all Caratheodory's results to algebraic form.

It is also appropriate to mention G.~M.~Golusin~\cite[p.~477]{Golusin:ref1}, who contributed to systematizing and
simplifying the presentation of the solution to the coefficient problem for the class of functions holomorphic and
bounded in $\Delta$.

\section{The general form of extremal functions}

The class $B$ is invariant under rotations in the variable plane $z$ and under rotations in the variable plane $w$
($w=f(z)$).
That is, if $n\in\mathbb{N}$ and the function $f\in B$ is extremal in the Krzyz problem for this index $n$, then for
$|\eta|=|\zeta|=1$ the function $\eta f(\zeta z)\in B$ and is extremal.
Moreover, a rotation in the variable plane $z$ does not affect the coefficient $\{f\}_0$.
Hence, without loss of generality, we may assume that if $f$ is extremal, then $\{f\}_0>0$ and $\{f\}_n>0$.
Below, when speaking of a function $f$ extremal in the Krzyz problem, we assume that $\{f\}_0>0$, $\{f\}_n>0$, and $f\in
B_t$, where $t=-\ln\{f\}_0$.
We will also often call a function extremal in the Krzyz problem for index $n$ simply extremal if the specific $n$ is
clear from the context.

\smallskip
Fix $n\in\mathbb{N}$.

If $f(z)=e^{-h(z)}$ is an extremal function, then, as is known (see~\cite{Stupin:ref1,Stupin:ref7}), the point
$h^{(n+1)}$ belongs to the boundary of $C^{(n+1)}$, which is equivalent to $M_n=0$ ($M_n$ is defined by
formula~\eqref{eq:Stupin:C-T_minors}).
By the Caratheodory-Toeplitz criterion, this means that the extension of the polynomial $Q_n$ defined by
formula~\eqref{eq:Stupin:C-T_polynom} is unique.
Hence every extremal function has the form $f(z)=e^{-h(z)}$, where $h$ is given by
formula~\eqref{eq:Stupin:B_t_extremal_h}.
Thus the following assertion holds:

\begin{Corollary}
    \label{cor:Stupin:B_t_extremal_Th}
Let $n\in\mathbb{N}$, let the function $f$ be extremal in the Krzyz problem for index $n$, and let $t\coloneqq
-\ln\{f\}_0$.
Then there exist numbers $m\leqslant n$, $\alpha_j>0$, $j=1,\ldots,m$, $\sum\limits_{j=1}^m\alpha_j=t$, and numbers
$z_j\coloneqq e^{i\varphi_j}$, $j=1,\ldots,m$, $0\leqslant\varphi_1<\ldots<\varphi_m<2\pi$, such that the function
$f(z)=e^{-h(z)}$, where the function $h$ is given by formula~\eqref{eq:Stupin:B_t_extremal_h}.
\end{Corollary}

This result has been well known since the early twentieth century (see, for example,~\cite[p.~171]{Hummel:ref1}
or~\cite[p.~725]{Martin:ref1}).

\smallskip
Let $n\in\mathbb{N}$, and let the function $f$, $\{f\}_0>0$, $\{f\}_n>0$, be extremal in the Krzyz problem.
Set $t\coloneqq-\ln\{f\}_0>0$, $z_j\coloneqq e^{i\varphi_j}$, $j=1,\ldots,n$.
Assuming $\alpha_j\geqslant0$, $j=1,\ldots,n$, formula~\eqref{eq:Stupin:BExtrF} and
Corollary~\ref{cor:Stupin:B_t_extremal_Th} imply that
\begin{equation*}
    f(z)
    =
    \prod\limits_{j=1}^{n}F(z_j z,\alpha_j),
    \quad
    f(0)=e^{-t}.
\end{equation*}

To find an extremal function using Corollary~\ref{cor:Stupin:B_t_extremal_Th}, one must determine at most $2n$ real
parameters $\alpha_j$, $\varphi_j$, $j=1,\ldots,n$.
If the value $t$ is not fixed in advance, the equality $\sum\limits_{j=1}^n\alpha_j=t$ determines $t$ and does not
reduce the number of independent parameters.
If $t$ is fixed and an extremal function is sought in the class $B_t$, the number of independent parameters does not
exceed $2n-1$.

If $f$ has real coefficients, then the coefficients of $h$ are also real, so for some $j$ such that $\alpha_j>0$ and
$0<\varphi_j<\pi$, formula~\eqref{eq:Stupin:B_t_extremal_h} must contain, besides the term $\alpha_j\frac{1+z_j z}{1-z_j
z}$, the term $\alpha_j\frac{1+\overline{z_j}z}{1-\overline{z_j}z}$.
It follows that if $f$ has real coefficients, finding an explicit representation of $f$ using
Corollary~\ref{cor:Stupin:B_t_extremal_Th} requires determining at most $n$ real parameters.
For fixed $t$ in the class $B_t$, the number of independent parameters does not exceed $n-1$.

\smallskip
Let $n\in\mathbb{N}$.
Consider the families of functions
\begin{multline*}
    \mathcal{B}(n)
    \coloneqq
    \bigg\{
    f
    :
    f=e^{-h},\,
    \{f\}_n>0,\,
    h(z)
    \coloneqq
    \sum\limits_{j=1}^n
    \alpha_j\tfrac{1+z_j z}{1-z_j z},\,
    \\
    \alpha_j>0,\,
    z_j\coloneqq e^{i\varphi_j},\,j=1,\ldots,n,\,
    0\leqslant\varphi_1<\ldots<\varphi_n<2\pi
    \bigg\}
\end{multline*}
and
\begin{multline*}
    \mathfrak{B}(n)
    \coloneqq
    \bigg\{
    f
    :
    f=e^{-h},\,
    \{f\}_n>0,\,
    h(z)
    \coloneqq
    \sum\limits_{j=1}^m
    \alpha_j\tfrac{1+z_j z}{1-z_j z},\,
    \\
    m\in\{1,\ldots,n\},\,
    \alpha_j>0,\,
    z_j\coloneqq e^{i\varphi_j},\,j=1,\ldots,m,\,
    0\leqslant\varphi_1<\ldots<\varphi_m<2\pi
    \bigg\}.
\end{multline*}
We will need the following evident assertion:

\begin{Theorem}
    \label{th:Stupin:denseExtr}
Let $n\in\mathbb{N}$.
The set $\mathcal{B}(n)$ is dense in $\mathfrak{B}(n)$ in the topology of locally uniform convergence in $\Delta$.
\end{Theorem}

\section{Topology}
\label{sec:Stupin:theTopology}

Let $X$ be a set on which we wish to introduce a topology.

\begin{Definition}
Suppose each point $x\in X$ is assigned a nonempty family $\mathcal N(x)$ of subsets of $X$.
A set $G\subset X$ is called open if for every $x\in G$ there exists a set $A\in\mathcal N(x)$ such that $A\subset G$.
\end{Definition}

\begin{Definition}
A set $U\subset X$ is called a neighborhood of a point $x\in X$ if there exists an open set $G$ such that $x\in G\subset
U$.
\end{Definition}

\begin{Definition}
A family $\mathcal N(x)$ of neighborhoods of a point $x\in X$ is called a neighborhood base at $x$ if for every
neighborhood $U$ of $x$ there exists a set $A\in\mathcal N(x)$ such that
    $A\subset U$.
\end{Definition}

\begin{Definition}
A family $\{\mathcal N(x)\}_{x\in X}$ is called a neighborhood system for a set $X$ if $\mathcal N(x)$ is a neighborhood
base at $x$ for every $x\in X$.
\end{Definition}

\smallskip
To introduce a topology using a neighborhood system, we use the following assertion
(see~\cite[p.~25]{Kantorovich:ref1}):

\begin{Assertion}
    \label{st:Stupin:neighborhoodSystemTopology}
Let $X$ be a set, and suppose each point $x\in X$ is assigned a nonempty family $\mathcal N(x)$ of subsets of $X$
satisfying the following conditions\textnormal:
    \newline\textnormal1.
$x\in A$ for every $A\in\mathcal N(x)$\textnormal;
    \newline\textnormal2.
if $A_1,A_2\in\mathcal N(x)$, there exists a set $A_3\in\mathcal N(x)$ such that $A_3\subset A_1\cap A_2$\textnormal;
    \newline\textnormal3.
if $A_x\in\mathcal N(x)$, there exists a set $A_x'\in\mathcal N(x)$ such that $A_x'\subset A_x$, and for every point
$y\in A_x'$ there exists a set $A_y\in\mathcal N(y)$ such that $A_y\subset A_x$.
    \newline
Then the collection of open subsets of $X$ is a topology on $X$ in which $\{\mathcal N(x)\}_{x\in X}$ is a neighborhood
system.
\end{Assertion}

Thus every open set is a union of sets from the neighborhood system.

\begin{Definition}
The topology introduced in this way will be called the topology generated by the neighborhood system $\{\mathcal
N(x)\}_{x\in X}$, and the set $X$ a topological space.
\end{Definition}

For a point $h\in C_1$ and $0<\varepsilon\leqslant1$, define the set
\[
    U_\varepsilon(h)
    \coloneqq
    \big\{
    (1-\alpha)\cdot h+\alpha\cdot g
    \,:\,
    g\in C_1,\ 0\leqslant\alpha<\varepsilon
    \big\}.
\]
Set $\mathcal N_{C_1}(h)\coloneqq\{U_\varepsilon(h):0<\varepsilon\leqslant1\}$.

\begin{Assertion}
The set $\{\mathcal N_{C_1}(h)\}_{h\in C_1}$ generates a topology on $C_1$ and is a neighborhood system in this
topology.
\end{Assertion}

\begin{proof}
We check that $\{\mathcal N_{C_1}(h)\}_{h\in C_1}$ satisfies all the conditions of
Assertion~\ref{st:Stupin:neighborhoodSystemTopology}.

    \noindent
    \textbf{1.}
The first condition holds because $h\in U_\varepsilon(h)$.

    \noindent
    \textbf{2.}
The second condition holds because
    $U_{\min\{\varepsilon_1,\varepsilon_2\}}(h)\subset U_{\varepsilon_1}(h)\cap U_{\varepsilon_2}(h)$.

    \noindent
    \textbf{3.}
To check the third condition, it suffices to prove that for every point $q\in U_\varepsilon(h)$ there exists a number
$\mu\in(0,1]$ such that $U_\mu(q)\subset U_\varepsilon(h)$.
Since $q\in U_\varepsilon(h)$, we have
    $q
    =
    (1-\alpha)h+\alpha g$,
    $g\in C_1$,
    $0\leqslant\alpha<\varepsilon$.
We show that we can take $\mu\coloneqq\tfrac{\varepsilon-\alpha}{1-\alpha}$.
Since $0\leqslant\alpha<\varepsilon\leqslant1$, we have $0<\mu\leqslant\varepsilon\leqslant1$.

By definition, the set $U_\mu(q)$ consists of functions
    \[
        u
        =
        (1-\beta)q+\beta p
        =
        (1-\beta)(1-\alpha)h+(1-\beta)\alpha g+\beta p,
        \quad
        p\in C_1,
        \quad
        0\leqslant\beta<\mu\leqslant1.
    \]
The coefficients of $h$, $g$, and $p$ are nonnegative and sum to $1$.
Hence $u$ is a convex combination of three functions $h$, $g$, and $p$.

We show that $u=(1-\gamma)h+\gamma w$, where $w\in C_1$.
Let $\gamma\coloneqq(1-\beta)\alpha+\beta$ be the sum of the coefficients of $g$ and $p$.
Clearly, $0\leqslant\gamma<1$.
For $\gamma=0$ we have $\alpha=\beta=0$, $u=h$.
Thus we may set $w\equiv1$.
For $\gamma>0$, set
        $w
        \coloneqq
        \frac{(1-\beta)\alpha}{\gamma}g
        +
        \frac{\beta}{\gamma}p$.
The coefficients of $g$ and $p$ are nonnegative and sum to $1$, so $w\in C_1$ by convexity of the class $C_1$.

By definition, $u\in U_\varepsilon(h)$ if $0\leqslant\gamma<\varepsilon$.
If $\gamma=0$, this inequality clearly holds, and for $\gamma>0$ the inequality holds for
    $\beta<\mu$.
If $u\in U_\mu(q)$, then $\beta<\mu$, so $\gamma<\varepsilon$ and $u\in U_\varepsilon(h)$.
Hence $U_\mu(q)\subset U_\varepsilon(h)$.

In Assertion~\ref{st:Stupin:neighborhoodSystemTopology} we can take $A_h=A_h'=U_\varepsilon(h)$, so its third condition
holds.

Thus all conditions of Assertion~\ref{st:Stupin:neighborhoodSystemTopology} hold, implying the required conclusion.
\end{proof}

\begin{Definition}
The topology generated by the neighborhood system $\{\mathcal N_{C_1}(h)\}_{h\in C_1}$ will be called the convex
topology on $C_1$ and denoted by $\tau_{C_1}$.
\end{Definition}

Note that the set $U_\varepsilon(h)$ is a neighborhood of $h$ in the topology $\tau_{C_1}$.
We call it an $\varepsilon$-neighborhood of $h$.

Let functions $h,g\in C_1$ and $\varepsilon\in(0,1)$ be given.
The functions $(1-\alpha)h+\alpha g$, $0\leqslant\alpha<\varepsilon$, form the initial part of the segment joining $h$
and $g$.
Thus geometrically $U_\varepsilon(h)$ is the union of the initial parts of all these segments.

\smallskip
Let $B_+\coloneqq\bigcup\limits_{t>0}B_t$.
Every function $f\in B_+$ has a unique representation $f=e^{-th}$, where $t>0$ and $h\in C_1$.
Set $\mathbb R^+\coloneqq\{t\in\mathbb R:t>0\}$.
Equip $\mathbb R^+$ with the usual topology induced from $\mathbb{R}$, $C_1$ with the convex topology, and $\mathbb
R^+\times C_1$ with the product topology.
Then the mapping $\Phi(t,h)\coloneqq e^{-th}$ is a bijection between the sets $\mathbb R^+\times C_1$ and $B_+$.

\begin{Definition}
A set $G\subset B_+$ is called open if the set $\Phi^{-1}(G)$ is open in $\mathbb R^+\times C_1$.
The resulting topology will be denoted by $\tau_{B_+}$.
\end{Definition}

For a function $f=e^{-th}\in B_+$, set

\[
    V_{\delta,\varepsilon}(f)
    \coloneqq
    \big\{
    e^{-\tau g}:
    \tau>0,\ |\tau-t|<\delta,\ g\in U_\varepsilon(h)
    \big\},
    \quad
    \mathcal N_{B_+}(f)
    \coloneqq
    \{V_{\delta,\varepsilon}(f):\delta>0,\ 0<\varepsilon\leqslant1\}.
\]
We have $V_{\delta,\varepsilon}(f)=\Phi\bigl(((t-\delta,t+\delta)\cap\mathbb R^+)\times U_\varepsilon(h)\bigr)$.
Since the families $\{(t-\delta,t+\delta)\cap\mathbb R^+:\delta>0\}$ and $\mathcal N_{C_1}(h)$ are neighborhood bases at
$t$ and $h$, respectively, the family $\mathcal N_{B_+}(f)$ is a neighborhood base at $f$.

\begin{Definition}
A function $f\in B_+$ will be called locally extremal in the Krzyz problem for index $n$ (with respect to $\tau_{B_+}$)
if there exist numbers $\delta>0$ and $0<\varepsilon\leqslant1$ such that
        $|\{g\}_n|
        \leqslant
        |\{f\}_n|$,
        $g\in V_{\delta,\varepsilon}(f)$.
An arbitrary nonconstant function $f\in B$ will be called locally extremal in the Krzyz problem for index $n$ if the
function $e^{-i\arg f(0)}f$ is locally extremal in $B_+$ for this index.
\end{Definition}

\begin{Lemma}
    \label{lem:Stupin:compactEstimateC1}
For every nonempty compact set $K\subset\Delta$ there exists a number $M_K>0$, depending only on $K$, such that
$\max\limits_{z\in K}|q(z)|\leqslant M_K$ and $\max\limits_{z\in K}|q(z)-h(z)|<2\varepsilon M_K$, where $h,g\in C_1$,
$q=(1-\alpha)h+\alpha g$, $0\leqslant\alpha<\varepsilon$, $0<\varepsilon\leqslant1$.
\end{Lemma}

\begin{proof}
Set $r\coloneqq\max\{|z|:z\in K\}<1$, $M_K\coloneqq\tfrac{1+r}{1-r}$.
By Caratheodory's estimation~\cite[p.~480]{GrigoriciucES_2021}, for every function $u\in C_1$ we have
$|u(z)|\leqslant\frac{1+|z|}{1-|z|}\leqslant M_K$, $z\in K$.

The first of the required estimations follows immediately.
Furthermore, $q-h=\alpha(g-h)$, so
    $\max\limits_{z\in K}|q(z)-h(z)|
    =
    \alpha\max\limits_{z\in K}|g(z)-h(z)|
    \leqslant
    \alpha\big(
    \max\limits_{z\in K}|g(z)|
    +
    \max\limits_{z\in K}|h(z)|
    \big)
    \leqslant
    2\alpha M_K
    <
    2\varepsilon M_K$.
\end{proof}

\begin{Lemma}
    \label{lem:Stupin:expLipschitz}
The mapping $w\mapsto e^{-w}$ is Lipschitz with constant $1$ in the right half-plane.
That is, if the numbers $a,b\in\mathbb C$ have positive real parts, then
    $|e^{-a}-e^{-b}|\leqslant|a-b|$.
\end{Lemma}

\begin{proof}
Let $\epsilon(\lambda)\coloneqq e^{-((1-\lambda)b+\lambda a)}$ and $0\leqslant\lambda\leqslant1$.
Then
    \[
        e^{-a}-e^{-b}
        =
        \epsilon(1)-\epsilon(0)
        =
        \int\limits_0^1 d\epsilon(\lambda)
        =
        (b-a)\int\limits_0^1 e^{-((1-\lambda)b+\lambda a)}\,d\lambda.
    \]
The real part of $(1-\lambda)b+\lambda a$ is positive, so the modulus of the integrand does not exceed one, implying the
required inequality.
\end{proof}

\begin{Assertion}
    \label{st:Stupin:segmentTopologyImpliesLocalUniform}
The topology $\tau_{B_+}$ is no weaker than the topology of locally uniform convergence in $\Delta$ induced on the set
$B_+$.
\end{Assertion}

\begin{proof}
Fix a function $f=e^{-th}\in B_+$, a nonempty compact set $K\subset\Delta$, and a number $\xi>0$.
It suffices to find numbers $\delta>0$ and $0<\varepsilon\leqslant1$ such that every function $\widetilde f\in
V_{\delta,\varepsilon}(f)$ satisfies the inequality
    $\max\limits_{z\in K}|\widetilde f(z)-f(z)|<\xi$.

Let $\widetilde f=e^{-\tau q}\in V_{\delta,\varepsilon}(f)$.
Then $|\tau-t|<\delta$ and $q\in U_\varepsilon(h)$.
By Lemma~\ref{lem:Stupin:compactEstimateC1}, there exists a number $M_K>0$ for which
        $\max\limits_{z\in K}|q(z)|
        \leqslant
        M_K$
and
        $\max\limits_{z\in K}|q(z)-h(z)|
        <
        2\varepsilon M_K$.

For every point $z\in K$, the numbers $\tau q(z)$ and $th(z)$ have positive real parts, since $q,h\in C_1$ and
$\tau,t>0$.
By Lemma~\ref{lem:Stupin:expLipschitz}, for $a=\tau q(z)$ and $b=th(z)$ we have
    \[
        |\widetilde f(z)-f(z)|
        \leqslant
        |\tau q(z)-th(z)|.
    \]
Hence $\max\limits_{z\in K}|\widetilde f(z)-f(z)|\leqslant\max\limits_{z\in K}|\tau q(z)-th(z)|$.
Adding and subtracting $tq(z)$, we obtain
            $\tau q(z)-th(z)
            =
            \tau q(z)-tq(z)+tq(z)-th(z)
            =
            (\tau-t)q(z)+t(q(z)-h(z))$.
Therefore, by the triangle inequality,
        $|\tau q(z)-th(z)|
        \leqslant
        |\tau-t|\,|q(z)|
        +
        t|q(z)-h(z)|$.
Taking the maximum over $z\in K$, we obtain
            $\max\limits_{z\in K}|\tau q(z)-th(z)|
            \leqslant
            |\tau-t|\max\limits_{z\in K}|q(z)|
            +
            t\max\limits_{z\in K}|q(z)-h(z)|
            <
            \delta M_K+2t\varepsilon M_K.$
Thus $\max\limits_{z\in K}|\widetilde f(z)-f(z)|<\delta M_K+2t\varepsilon M_K$.
Choose $\delta>0$ and $0<\varepsilon\leqslant1$ so that $\delta M_K+2t\varepsilon M_K<\xi$.
Then
    $V_{\delta,\varepsilon}(f)\subset\big\{\widetilde f\in B_+:\max\limits_{z\in K}|\widetilde f(z)-f(z)|<\xi\big\}$.
Hence every neighborhood of the function $f$ in the topology of locally uniform convergence contains a neighborhood of
it in the topology $\tau_{B_+}$, which proves the assertion.
\end{proof}

Thus convergence in the topology $\tau_{B_+}$ implies locally uniform convergence.
Also, local extremality in the topology of locally uniform convergence implies local extremality in the topology
$\tau_{B_+}$.
In general, the converse does not follow, since neighborhoods may be smaller in a stronger topology.

\smallskip
For every $t>0$, the mapping $\Phi_t(h)\coloneqq e^{-th}$ is a bijection between the sets $C_1$ and $B_t$.

\begin{Definition}
A set $G\subset B_t$ is called open if the set $\Phi_t^{-1}(G)$ is open in $C_1$.
The resulting topology will be denoted by $\tau_{B_t}$.
\end{Definition}

The topology $\tau_{B_t}$ coincides with the subspace topology of $B_+$.
This follows from the equality $B_t=\Phi(\{t\}\times C_1)$ and the definition of the product topology.
This topology and the corresponding notion of local extremality were considered in~\cite{Stupin:ref8}.

\smallskip
By Corollary~\ref{cor:Stupin:B_t_extremal_Th}, if $f$ is a function extremal in the Krzyz problem for index $n$,
$\{f\}_0>0$, and $\{f\}_n>0$, then $f\in\mathfrak{B}(n)$.
Let $\mathfrak E(n)$ denote the subset of $\mathfrak{B}(n)$ consisting of functions locally extremal in the topology
$\tau_{B_+}$.

\begin{Assertion}
For every index $n$ the set $\mathfrak E(n)$ is nonempty.
\end{Assertion}

\begin{proof}
As stated in the introduction, a global extremal exists.
Rotations in the variable planes $z$ and $w$ can make its coefficients with indices $0$ and $n$ positive.
Every global extremal is a local extremal, so a global extremal normalized by rotations in the variable planes $z$ and
$w$ belongs to $\mathfrak E(n)$.
\end{proof}

Another way to establish that $\mathfrak E(n)$ is nonempty is given below in the proof of
Theorem~\ref{th:Stupin:uniqunessOfExtremalFunctionInB}.

\section{Some properties of extremal functions}

The results of this section largely resemble those in~\cite[p.~726]{Martin:ref1}.
The difference is that~\cite{Martin:ref1} considers only properties of global extremals, whereas here local extremals
are also studied.

\smallskip
Note that if $n\in\mathbb{N}$ and $f$ and $g$ are functions holomorphic in $\Delta$, then
\begin{equation}
    \label{eq:Stupin:prodCoeff}
    \{f \cdot g\}_n
    =
    \{f\}_n\{g\}_0+\{f\}_{n-1}\{g\}_1+\ldots+\{f\}_0\{g\}_n.
\end{equation}
\begin{Lemma}
    \label{lem:Stupin:varCoeff}
Let $n\in\mathbb{N}$, let $f$ and $g$ be functions holomorphic in $\Delta$, let $\varepsilon\in\mathbb{R}$, and let
    \begin{equation}
        \label{eq:Stupin:var}
        f_{\varepsilon}(z)\coloneqq f(z)e^{-\varepsilon g(z)},
    \end{equation}
then
    \begin{equation}
        \label{eq:Stupin:varCoeff}
        \{f_{\varepsilon}\}_n
        =
        \{f\}_n
        -
        \varepsilon
        \{f \cdot g\}_n
        +
        o(\varepsilon),
        \quad
        \varepsilon\to0.
    \end{equation}
\end{Lemma}

\begin{proof}
For the variation $f_{\varepsilon}$ of the function $f$ by the function $e^{-\varepsilon g}$,
$\varepsilon\in\mathbb{R}$, for every fixed $z\in\Delta$ we have
    \[
        f_{\varepsilon}(z)
        \coloneqq
        f(z)e^{-\varepsilon g(z)}
        =
        f(z)(1-\varepsilon g(z)+o(\varepsilon g(z)))
        =
        f(z) - \varepsilon f(z) g(z)+o(\varepsilon g(z)),
        \quad
        \varepsilon\to0.
    \]
Moreover, $\varepsilon g(z)\to0$ locally uniformly in $\Delta$ as $\varepsilon\to0$.
Fixing a number $r\in(0,1)$ and calculating
    \[
        \{f_{\varepsilon}\}_n
        =
        \frac1{2\pi i}\int\limits_{|z|=r}\frac{f_{\varepsilon}(z)}{z^{n+1}}dz
        =
        \frac1{2\pi i}\int\limits_{|z|=r}\frac{f(z)}{z^{n+1}}dz
        -
        \frac{\varepsilon}{2\pi i}\int\limits_{|z|=r}\frac{f(z)g(z)}{z^{n+1}}dz
        +
        o(\varepsilon),
    \]
we obtain formula~\eqref{eq:Stupin:varCoeff}.
\end{proof}

\begin{Theorem}
    \label{th:Stupin:scalar_prod}
Let $n\in\mathbb{N}$, $f\in\mathfrak E(n)$; then for every function $g\in C$
    \begin{equation}
        \label{eq:Stupin:scalar_prod}
        \operatorname{Re}\{f\cdot g\}_n
        \geqslant
        0,
    \end{equation}
with
    \begin{equation}
        \label{eq:Stupin:orthogonalVects}
        \operatorname{Re}\{f\cdot h\}_n
        =
        0
\quad\text{for}\quad
        h=-\ln f.
    \end{equation}
\end{Theorem}

\begin{proof}
    \textbf{1.}
We prove inequality~\eqref{eq:Stupin:scalar_prod}.
Let $g\in C$.
The variation $f_{\varepsilon}$ of the function $f$ defined by formula~\eqref{eq:Stupin:var} is internal for
$\varepsilon>0$, that is, $f_{\varepsilon} \in B$.

Set $t_h\coloneqq\{h\}_0$ and $t_g\coloneqq\{g\}_0$.
Then $t_\varepsilon\coloneqq t_h+\varepsilon t_g$ and
    \[
        h_{\varepsilon}
        \coloneqq
        \frac{h+\varepsilon g}{t_\varepsilon}
        =
        \frac{t_h}{t_\varepsilon}\cdot\frac{h}{t_h}
        +
        \frac{\varepsilon t_g}{t_\varepsilon}\cdot\frac{g}{t_g},
        \quad
        h_{\varepsilon},\frac{h}{t_h},\frac{g}{t_g}\in C_1,
        \quad
        \frac{t_h}{t_\varepsilon}+\frac{\varepsilon t_g}{t_\varepsilon}=1.
    \]
Let $\delta>0$ and $0<\eta\leqslant1$.
For sufficiently small $\varepsilon>0$, we have $|t_\varepsilon-t_h|=\varepsilon t_g<\delta$ and
    $\tfrac{\varepsilon t_g}{t_\varepsilon}<\eta$.
Hence $h_{\varepsilon}\in U_\eta(\tfrac{h}{t_h})$ and
    $f_{\varepsilon}=e^{-t_\varepsilon h_{\varepsilon}}\in V_{\delta,\eta}(f)$.
Thus every neighborhood of the function $f$ in the topology $\tau_{B_+}$ contains all functions $f_{\varepsilon}$ for
sufficiently small $\varepsilon>0$.

Lemma~\ref{lem:Stupin:varCoeff} holds for $\varepsilon\to0$, and since by assumption $\{f\}_n>0$ and $f$ is extremal, we
have $\operatorname{Re}{\{f_{\varepsilon}\}_n}\leqslant\{f\}_n$.
This implies (see~\eqref{eq:Stupin:varCoeff}) inequality~\eqref{eq:Stupin:scalar_prod}.

    \textbf{2.}
We prove equality~\eqref{eq:Stupin:orthogonalVects}.
Now take $g=h$.
In this case, the variation $f_{\varepsilon}$ of the function $f$ defined by formula~\eqref{eq:Stupin:var} is internal
for $-1<\varepsilon<0$, that is, $f_{\varepsilon}=f\cdot f^{\varepsilon}\in B$.
Set $t_{\varepsilon}\coloneqq(1+\varepsilon)t_h$ and
    $h_{\varepsilon}\coloneqq\tfrac{h}{t_h}$.
Then $f_{\varepsilon}=e^{-t_{\varepsilon}h_{\varepsilon}}$.
Let $\delta>0$ and $0<\eta\leqslant1$.
For sufficiently small $\varepsilon<0$, we have $|t_{\varepsilon}-t_h|=-\varepsilon t_h<\delta$ and
    $h_{\varepsilon}=\tfrac{h}{t_h}\in U_\eta(\tfrac{h}{t_h})$.
Hence $f_{\varepsilon}\in V_{\delta,\eta}(f)$.
Thus every neighborhood of the function $f$ in the topology $\tau_{B_+}$ contains all functions $f_{\varepsilon}$ for
sufficiently small $\varepsilon<0$.
Since by assumption $\{f\}_n>0$ and $f$ is extremal, formula~\eqref{eq:Stupin:varCoeff} gives $\operatorname{Re}\{f\cdot
h\}_n\leqslant0$ for $\varepsilon<0$, while for $\varepsilon>0$ inequality~\eqref{eq:Stupin:scalar_prod} holds, that is,
$\operatorname{Re}\{f\cdot h\}_n\geqslant0$, implying equality~\eqref{eq:Stupin:orthogonalVects}.
\end{proof}

\begin{Theorem}
    \label{th:Stupin:f_n_ge_2f_0}
Let $n\in\mathbb{N}$, $f\in\mathfrak E(n)$; then
    \begin{equation}
        \label{eq:Stupin:f_n_ge_2f_0}
        \{f\}_n\geqslant 2\{f\}_0.
    \end{equation}
\end{Theorem}

\begin{proof}
Substituting the coefficients of the function
    \[
        g(z)=\frac{1-z^n}{1+z^n}=1-2z^n+\ldots
    \]
into inequality~\eqref{eq:Stupin:scalar_prod} gives inequality~\eqref{eq:Stupin:f_n_ge_2f_0}.
\end{proof}

\begin{Theorem}
    \label{th:Stupin:extrPolynomial}
Let $n\in\mathbb{N}$, $f\in\mathfrak E(n)$, and
    \begin{equation}
        \label{eq:Stupin:extrPolynomialH_n}
        H(z) \coloneqq \{f\}_n+2\{f\}_{n-1}z+\ldots+2\{f\}_0 z^n,
    \end{equation}
then
    \begin{equation}
        \label{eq:Stupin:extrPolynomial}
        H\in C^n.
    \end{equation}
Furthermore, let $m\leqslant n$, and let the numbers $z_j\coloneqq e^{i\varphi_j}$, $j=1,\ldots,m$,
$0\leqslant\varphi_1<\ldots<\varphi_m<2\pi$ partially determine the function $h\coloneqq -\ln f$ \textnormal{(}see
formula~\eqref{eq:Stupin:B_t_extremal_h}\textnormal{)}.
Then
    \begin{equation}
        \label{eq:Stupin:extrPolynomial_phi}
        \operatorname{Re}{H(z_j)}=0,
        \quad
        j=1,\ldots,m.
    \end{equation}
\end{Theorem}

\begin{proof}
    \textbf{1.}
We prove statement~\eqref{eq:Stupin:extrPolynomial}.
Fix $\zeta\in\overline\Delta$ and set
    \[
        g(z)
        \coloneqq
        \frac{1 + \zeta z}{1 - \zeta z}=1+2 \zeta z+2 \zeta^2 z^2+\ldots
    \]
Clearly, $g\in C_1$.
Substituting the coefficients of the function $g$ into inequality~\eqref{eq:Stupin:scalar_prod} gives
    \[
        \operatorname{Re}\{f\cdot g\}_n
        =
        \operatorname{Re} {(\{f\}_n+2\{f\}_{n-1}\zeta+\ldots+2\{f\}_0\zeta^n)}
        \geqslant
        0,
        \quad
        |\zeta|
        \leqslant
        1,
    \]
which is equivalent to $\operatorname{Re}{H(z)}\geqslant0$, $z\in\overline\Delta$.
Since $H(0)>0$, the open mapping principle gives $\operatorname{Re}H(z)>0$, $z\in\Delta$, which is equivalent to
statement~\eqref{eq:Stupin:extrPolynomial}.

    \textbf{2.}
We prove equalities~\eqref{eq:Stupin:extrPolynomial_phi}.
Fix $j\in\{1,\ldots,m\}$ and set $z_j\coloneqq e^{i\varphi_j}$, where the number $\varphi_j$ partially determines the
function $h$ (see formula~\eqref{eq:Stupin:B_t_extremal_h}), and consider the function
    \[
        g(z)
        \coloneqq
        \frac{1+z_j z}{1-z_j z}
        =
        1+2z_j z+2z_j^2 z^2+\ldots
    \]
Then $f_{\varepsilon}=f\cdot e^{-\varepsilon g}\in B$ both for $\varepsilon\geqslant0$ and for sufficiently small
    $\varepsilon<0$.
That is, for $-\alpha_j<\varepsilon<0$, where $\alpha_j$ partially determines the function $h$ (see
formula~\eqref{eq:Stupin:B_t_extremal_h}).
Set $t_h\coloneqq\{h\}_0$.
Then $\{g\}_0=1$.
Let $\delta>0$ and $0<\eta\leqslant1$.

First consider the case $\varepsilon>0$.
Then $t_\varepsilon\coloneqq t_h+\varepsilon$ and
    \[
        h_\varepsilon
        \coloneqq
        \frac{h+\varepsilon g}{t_\varepsilon}
        =
        \frac{t_h}{t_\varepsilon}\cdot\frac{h}{t_h}
        +
        \frac{\varepsilon}{t_\varepsilon}g,
        \quad
        h_\varepsilon,\frac{h}{t_h},g\in C_1,
        \quad
        \frac{t_h}{t_\varepsilon}+\frac{\varepsilon}{t_\varepsilon}=1.
    \]
For sufficiently small $\varepsilon>0$, we have $|t_\varepsilon-t_h|=\varepsilon<\delta$ and
    $\tfrac{\varepsilon}{t_\varepsilon}<\eta$.
Hence $h_\varepsilon\in U_\eta(\tfrac{h}{t_h})$ and
    $f_\varepsilon=e^{-t_\varepsilon h_\varepsilon}\in V_{\delta,\eta}(f)$.

Now consider the case $-\alpha_j<\varepsilon<0$.
Then $t_\varepsilon\coloneqq t_h+\varepsilon$, $h_\varepsilon\coloneqq\frac{h+\varepsilon g}{t_\varepsilon}$, and
    $f_\varepsilon=e^{-t_\varepsilon h_\varepsilon}$.

    \noindent\textbf{a.}
If $t_h-\alpha_j=0$ (which is possible only for $m=1$), then $h=\alpha_j g$ and $h_\varepsilon=g=\tfrac{h}{t_h}$ for all
$-\alpha_j<\varepsilon<0$.
Thus for sufficiently small $\varepsilon<0$ we have $|t_\varepsilon-t_h|=-\varepsilon<\delta$, $h_\varepsilon\in
U_\eta(\tfrac{h}{t_h})$, and $f_\varepsilon\in V_{\delta,\eta}(f)$.

    \noindent\textbf{b.}
If $t_h-\alpha_j>0$, then the function $p\coloneqq\frac{h-\alpha_j g}{t_h-\alpha_j}\in C_1$ and
    \begin{multline*}
        h_\varepsilon
        =
        \frac{h+\varepsilon g}{t_\varepsilon}
        =
        \frac{h+\frac{\varepsilon}{\alpha_j}h-\frac{\varepsilon}{\alpha_j}h+\varepsilon g}
        {t_\varepsilon}
        =
        \frac{(1+\frac{\varepsilon}{\alpha_j})h-\frac{\varepsilon}{\alpha_j}(h-\alpha_j g)}
        {t_\varepsilon}
        =\\=
        \frac{\alpha_j+\varepsilon}{\alpha_j t_\varepsilon}t_h\frac{h}{t_h}
        -
        \varepsilon\frac{t_h-\alpha_j}{\alpha_j t_\varepsilon}p
        =
        (1-\gamma_\varepsilon)\frac{h}{t_h}+\gamma_\varepsilon p,
    \end{multline*}
where $\gamma_\varepsilon\coloneqq -\tfrac{\varepsilon(t_h-\alpha_j)}{\alpha_j t_\varepsilon}$.
For $\varepsilon\to0-$, we have $\gamma_\varepsilon\to0$.
Thus for sufficiently small $\varepsilon<0$ the inequalities $|t_\varepsilon-t_h|=-\varepsilon<\delta$ and
$\gamma_\varepsilon<\eta$ hold.
Hence $h_\varepsilon\in U_\eta(\tfrac{h}{t_h})$ and
    $f_\varepsilon\in V_{\delta,\eta}(f)$.

Thus every neighborhood of the function $f$ in the topology $\tau_{B_+}$ contains all functions $f_\varepsilon$ for
sufficiently small $|\varepsilon|$.
Lemma~\ref{lem:Stupin:varCoeff}, the condition $\{f\}_n>0$, and extremality of $f$ imply that for $\varepsilon<0$
    \[
        \operatorname{Re}\{f\cdot g\}_n
        =
        \operatorname{Re}{(\{f\}_n+\{f\}_{n-1}2z_j+\ldots+\{f\}_0 2z_j^n)}\leqslant0,
    \]
and for $\varepsilon>0$
    \[
        \operatorname{Re}\{f\cdot g\}_n
        =
        \operatorname{Re}{(\{f\}_n+\{f\}_{n-1}2z_j+\ldots+\{f\}_0 2z_j^n)}\geqslant0,
    \]
which is equivalent to the equality $\operatorname{Re}{H(z_j)}=0$.
\end{proof}

\smallskip
Let $\mathfrak E^r(n)$ denote the subset of functions $f\in\mathfrak E(n)$ with real coefficients, and reformulate
Theorem~\ref{th:Stupin:scalar_prod} for $f\in B^r$.

\begin{Corollary}
Let $n\in\mathbb{N}$, $f\in\mathfrak E^r(n)$; then $\{f\cdot g\}_n\geqslant0$ for every function $g\in C^r$.
In particular, equality in this inequality is attained if $g=-\ln f$.
\end{Corollary}

\begin{proof}
The required inequality follows from inequality~\eqref{eq:Stupin:scalar_prod}.
If $f$ has real coefficients, so does $h$; hence the required equality follows immediately from
equality~\eqref{eq:Stupin:orthogonalVects}.
\end{proof}

\section{On variations on $B_t$}

Let $n\in\mathbb N$, $f\in\mathfrak E(n)$, $h\coloneqq-\ln f\in C_{t_h}$, $t_h>0$, and let $g\in C_{t_g}$, $t_g>0$.
Consider the following variation of the function $f$, internal in the class $B$ for $\varepsilon\in[0,1]$:
\[
    f_{\varepsilon}(z)
    \coloneqq
    f(z)^{1-\varepsilon} e^{-\varepsilon g(z)}
    =
    e^{-(1-\varepsilon)h(z)}e^{-\varepsilon g(z)}
    =
    e^{-h(z)}e^{-\varepsilon(g(z)-h(z))}
    =
    f(z)e^{-\varepsilon(g(z)-h(z))}.
\]
Set $t_\varepsilon\coloneqq(1-\varepsilon)t_h+\varepsilon t_g$.
The limit of the mapping $\varepsilon\mapsto t_\varepsilon$ as $\varepsilon\to0+$ is $t_h$, and
\[
    h_{\varepsilon}
    \coloneqq
    \frac{(1-\varepsilon)h+\varepsilon g}{t_\varepsilon}
    =
    \frac{(1-\varepsilon)t_h}{t_\varepsilon}\frac{h}{t_h}
    +
    \frac{\varepsilon t_g}{t_\varepsilon}\frac{g}{t_g},
    \quad
    h_{\varepsilon},\frac{h}{t_h},\frac{g}{t_g}\in C_1.
\]
Hence the limit of the mapping $\varepsilon\mapsto h_\varepsilon$ in the convex topology as $\varepsilon\to0+$ is
$h/t_h$.
Therefore the limit of the mapping $\varepsilon\mapsto f_\varepsilon$ in the topology $\tau_{B_+}$ as $\varepsilon\to0+$
is $f$.
By Lemma~\ref{lem:Stupin:varCoeff}, we have
    $\{f_{\varepsilon}\}_n
    =
    \{f\}_n
    -
    \varepsilon\{f\cdot(g-h)\}_n
    +
    o(\varepsilon)$,
    $\varepsilon\to0$.

Extremality of the function $f$ and the condition $\{f\}_n>0$ give
$
    \operatorname{Re}\{f_{\varepsilon}\}_n\leqslant\{f\}_n
$,
whence
$
    \operatorname{Re}\{f\cdot(g-h)\}_n\geqslant0
$
or
$
    \operatorname{Re}\{f\cdot g\}_n
    \geqslant
    \operatorname{Re}\{f\cdot h\}_n
$.
But
$
    \operatorname{Re}\{f\cdot h\}_n=0
$
by equality~\eqref{eq:Stupin:orthogonalVects}.
Thus we again obtain $\operatorname{Re}\{f\cdot g\}_n\geqslant0$, $g\in C_{t_g}$, that is,
inequality~\eqref{eq:Stupin:scalar_prod}, since $t_g$ is an arbitrary positive number.

Similarly, if $t\coloneqq t_h=t_g$ and $f\in B_t$ is a function extremal on $B_t$, then $h,g\in C_t$ and
$h_\varepsilon\in C_1$, $0\leqslant\varepsilon\leqslant1$, since $C_t$ has a convex structure.
Hence the function $f_{\varepsilon}=e^{-t_h h_\varepsilon}\in B_t$, that is, the variation is internal in the class
$B_t$.
The limit of the mapping $\varepsilon\mapsto h_\varepsilon$ in the convex topology as $\varepsilon\to0+$ is $h/t_h$, so
the limit of the mapping $\varepsilon\mapsto f_\varepsilon$ in the topology $\tau_{B_t}$ as $\varepsilon\to0+$ is $f$.
By Lemma~\ref{lem:Stupin:varCoeff}, we have
\[
    \{f_{\varepsilon}\}_n
    =
    \{f\}_n
    -
    \varepsilon\{f\cdot(g-h)\}_n
    +
    o(\varepsilon),
    \quad
    \varepsilon\to0.
\]
Extremality of the function $f$ under the normalization $\{f\}_n>0$ gives
$
    \operatorname{Re}\{f_{\varepsilon}\}_n\leqslant\{f\}_n
$,
whence for every function $g\in C_t$ we have
$
    \operatorname{Re}\{f\cdot g\}_n
    \geqslant
    \operatorname{Re}\{f\cdot h\}_n
$.

Note that equality~\eqref{eq:Stupin:orthogonalVects} does not hold in general on $B_t$.
For example, for $n=1$, the function $f=e^{-h}$, where $h(z)=t(1-z)/(1+z)$, is extremal on $B_t$
(see~\cite[p.~88]{Stupin:ref14}); for this function $\{f\}_0=e^{-t}$, $\{f\}_1=2te^{-t}$, $\{h\}_0=t$, and
$\{h\}_1=-2t$.
Hence,
$\{f\}_1\{h\}_0+\{f\}_0\{h\}_1=2te^{-t}(t-1)$,
which is nonzero for all positive $t\neq1$.

Thus if $f\in\mathfrak E(n)$, then on the class $C$ the function $h=-\ln f$ is a point of global minimum of the real
linear functional
$L_f(g)\coloneqq\operatorname{Re}\{f\cdot g\}_n$.
If $f$ is extremal on $B_t$, then for fixed $t>0$ the same holds on the class $C_t$.

The inequality $L_f(g)\geqslant L_f(h)$, where $g\in C$ (or $g\in C_t$), means that the function $h=-\ln f$ is a
supporting point from below of the functional $L_f$ on the class $C$ (or $C_t$) and belongs to the corresponding
supporting face defined by $L_f(g)=L_f(h)$, where $g\in C$ (or $g\in C_t$).

\section{Functions of extremal type}

In this section we introduce the notion of a function of extremal type, important for all that follows, and show in
particular that statement~\eqref{eq:Stupin:extrPolynomial} implies inequality~\eqref{eq:Stupin:scalar_prod}, while the
representation of the function $h$ by formula~\eqref{eq:Stupin:B_t_extremal_h} and
equalities~\eqref{eq:Stupin:extrPolynomial_phi} imply equality~\eqref{eq:Stupin:orthogonalVects}.

Let $n\in\mathbb{N}$.
Let $\mathbb E(n)$ denote the subset of $\mathfrak{B}(n)$ consisting of functions $f$ for which $H\in C^n$ and
$\operatorname{Re}{H(z_j)}=0$, $j=1,\ldots,m$, where $m\in\{1,\ldots,n\}$.

Recall that $H(z)\coloneqq\{f\}_n+2\{f\}_{n-1}z+\ldots+2\{f\}_0 z^n$, $z_j\coloneqq e^{i\varphi_j}$,
$j=1,\ldots,m$, $0\leqslant\varphi_1<\ldots<\varphi_m<2\pi$.
Note that the parameters $z_j$ partially determine the function $f$ (see the definition of $\mathfrak B(n)$).

Functions $f\in\mathbb E(n)$ will be called functions of extremal type in the Krzyz problem for the coefficient with
index $n$, or simply functions of extremal type if this does not obscure the meaning.

\begin{Theorem}
    \label{th:Stupin:localExtremalIsExtremalType}
Let $n\in\mathbb N$.
Every function in $\mathfrak E(n)$ is of extremal type, that is,
    $\mathfrak E(n)\subset\mathbb E(n)$.
\end{Theorem}

\begin{proof}
By definition, $f\in\mathfrak E(n)$ if $f\in\mathfrak B(n)$ and $f$ is locally extremal in the topology $\tau_{B_+}$.
If $f\in\mathfrak E(n)$, then $H\in C^n$ and $\operatorname{Re}{H(z_j)}=0$, $j=1,\ldots,m$, by
Theorem~\ref{th:Stupin:extrPolynomial}.
The definition of $\mathbb E(n)$ immediately implies that if $f\in\mathfrak E(n)$, then $f\in\mathbb E(n)$, as required.
\end{proof}

\begin{Lemma}
    \label{lem:Stupin:fg}
Let $n\in\mathbb{N}$, let $f$ be a function holomorphic in $\Delta$, and let $g\in C$.
Then there exist numbers $\lambda\geqslant0$ and $m'\in\{0,\ldots,n\}$, and for $m'>0$ also numbers $\beta_k>0$,
$k=1,\ldots,m'$, and $0\leqslant\theta_1<\ldots<\theta_{m'}<2\pi$, such that
    \begin{equation*}
        \{f\cdot g\}_n
        =
        \lambda\{f\}_n
        +
        \sum\limits_{j=1}^{m'}\beta_j H(e^{i\theta_j}),
    \end{equation*}
where $H(z)\coloneqq\{f\}_n+2\{f\}_{n-1}z+\ldots+2\{f\}_0z^n$.
\end{Lemma}

\begin{proof}
By Theorem~\ref{th:Stupin:CaratheodoryFejer}, there exists a unique function of the form
    \[
        g_{m'}(z)
        \coloneqq
        \lambda+
        \sum\limits_{j=1}^{m'}\beta_j\dfrac{1+e^{i\theta_j}z}{1-e^{i\theta_j}z}
        =
        \lambda+
        \sum\limits_{j=1}^{m'}\beta_j(1+2e^{i\theta_j}z+\ldots+2e^{in\theta_j}z^n+\ldots),
    \]
where $m'\in\{0,\ldots,n\}$, $\lambda\geqslant0$, $\beta_j>0$, $j=1,\ldots,m'$,
$0\leqslant\theta_1<\ldots<\theta_{m'}<2\pi$, whose first $n+1$ Taylor coefficients are respectively equal to
$\{g\}_0,\ldots,\{g\}_n$.
Note that the numbers $\beta_j$ and $\theta_j$ are meaningful only for $m'>0$.
Thus, applying formula~\eqref{eq:Stupin:prodCoeff}, we obtain
    \begin{equation*}
        \{f\cdot g\}_n
        =
        \lambda\{f\}_n
        +
        \sum\limits_{k=1}^{m'}\beta_k(\{f\}_n+\{f\}_{n-1}2e^{i\theta_k}+\ldots+\{f\}_0 2e^{in\theta_k})
        =
        \lambda\{f\}_n
        +
        \sum\limits_{k=1}^{m'}\beta_k H(e^{i\theta_k}).
    \end{equation*}
This proves the assertion.
\end{proof}

\begin{Theorem}
    \label{th:Stupin:scalar_prod_criteria}
Let $n\in\mathbb{N}$, $f\in\mathbb E(n)$, $g\in C$.
Then $\operatorname{Re}\{f\cdot g\}_n\geqslant0$, with equality attained if and only if
    $g(z)\coloneqq\sum\limits_{j=1}^{m'}\beta_j\dfrac{1+e^{i\theta_j}z}{1-e^{i\theta_j}z}$,
    $m'\in\{1,\ldots,n\}$, $\beta_j>0$, $0\leqslant\theta_1<\ldots<\theta_{m'}<2\pi$,
    $\operatorname{Re}H(e^{i\theta_j})=0$, $j=1,\ldots,m'$.
\end{Theorem}

\begin{proof}
We must show that $\operatorname{Re}\{f\cdot g\}_n\geqslant0$ and establish all equality cases in this inequality.

Since $H\in C^n$ and $H$ is continuous in $\overline\Delta$, we have $\{f\}_n=H(0)>0$,
$\operatorname{Re}{H(z)}\geqslant0$ for $z\in\overline\Delta$.
Since $\{f\}_n>0$ and $\beta_j>0$ for all $j=1,\ldots,m'$, Lemma~\ref{lem:Stupin:fg} gives
    $\operatorname{Re}\{f\cdot g\}_n
    =
    \lambda\{f\}_n+
    \sum\limits_{j=1}^{m'}\beta_j\operatorname{Re}H(e^{i\theta_j})
\geqslant0$, with equality attained if and only if $\lambda=0$ and $\operatorname{Re}H(e^{i\theta_j})=0$ for all $j$.

We describe all equality cases.
Let $\mathfrak g\in C$ and $\operatorname{Re}\{f\cdot\mathfrak g\}_n=0$; then $\lambda=0$.
By Theorem~\ref{th:Stupin:nCTEquiv}, the polynomial $Q$ formed from the first $n+1$ coefficients of the function
$\mathfrak g$ has a unique extension to a function in the class $C$.
Since the functions $\mathfrak g$ and $g$ belong to the class $C$ and their first $n+1$ coefficients coincide,
$\mathfrak g=g$.

Conversely, if $\operatorname{Re}H(e^{i\theta_k})=0$, $k=1,\ldots,m'$, and the function $g$ has the form stated in this
theorem, then $\operatorname{Re}\{f\cdot g\}_n=0$ by Lemma~\ref{lem:Stupin:fg}.
\end{proof}

\smallskip
Thus, for a function $f\in\mathfrak{B}(n)$, statement~\eqref{eq:Stupin:extrPolynomial} is equivalent to
inequality~\eqref{eq:Stupin:scalar_prod} holding for all functions $g\in C$.
Indeed, statement~\eqref{eq:Stupin:extrPolynomial} and Lemma~\ref{lem:Stupin:fg} imply
inequality~\eqref{eq:Stupin:scalar_prod}.
Conversely, substituting the functions $g(z)\coloneqq\frac{1+\zeta z}{1-\zeta z}$, $|\zeta|\leqslant1$ into
inequality~\eqref{eq:Stupin:scalar_prod} gives $\operatorname{Re}H(\zeta)\geqslant0$.
Since $H(0)=\{f\}_n>0$, we have $H\in C^n$.

\begin{Assertion}
Let $n\in\mathbb N$, $f\in\mathbb E(n)$.
Then for every function $g\in C$, $\operatorname{Re}\{f\cdot g\}_n\geqslant0$, with $\operatorname{Re}\{f\cdot h\}_n=0$
for $h=-\ln f$.
\end{Assertion}

\begin{proof}
Let $n\in\mathbb{N}$, $f\in\mathbb E(n)$, $g\in C$.
The inequality $\operatorname{Re}\{f\cdot g\}_n\geqslant0$ follows from Theorem~\ref{th:Stupin:scalar_prod_criteria}.

By the definition of the set $\mathbb E(n)$, the function $h=-\ln f=\sum\limits_{j=1}^m\alpha_j\frac{1+z_j z}{1-z_j z}$,
    $m\in\{1,\ldots,n\}$, $\alpha_j>0$, $z_j\coloneqq e^{i\varphi_j}$,
    $0\leqslant\varphi_1<\ldots<\varphi_m<2\pi$,
    $\operatorname{Re}H(e^{i\varphi_j})=0$, $j=1,\ldots,m$.
Hence the equality $\operatorname{Re}\{f\cdot h\}_n=0$ holds by Theorem~\ref{th:Stupin:scalar_prod_criteria} for $g=h$.
\end{proof}

\smallskip
Statement~\eqref{eq:Stupin:extrPolynomial} together with equalities~\eqref{eq:Stupin:extrPolynomial_phi} is strictly
stronger than inequality~\eqref{eq:Stupin:scalar_prod}, as the following example shows:

\begin{Example}
Let $n\in\mathbb N$, $t>1$, and $f(z)\coloneqq F(-z^n,t)$.
Then $f\in\mathcal B(n)$, and the polynomial $H$ generated by the coefficients of the function $f$ according to
formula~\eqref{eq:Stupin:extrPolynomialH_n} belongs to the class $C^n$, but
equalities~\eqref{eq:Stupin:extrPolynomial_phi} do not hold.
\end{Example}

\begin{proof}
Let $z_j\coloneqq e^{i\varphi_j}$, $\varphi_j\coloneqq\frac{(2j-1)\pi}{n}$,
    $\alpha_j\coloneqq\frac{t}n$, $j=1,\ldots,n$.
Then
    \[
        h(z)
        \coloneqq
        \sum\limits_{j=1}^n\alpha_j\frac{1+z_j z}{1-z_j z}
        =
        t\frac{1-z^n}{1+z^n}
\quad\text{and}\quad
        f=e^{-h}=F(-z^n,t).
    \]
Hence $f(z)=e^{-t}+2te^{-t}z^n+o(z^n)$, so $f\in\mathcal B(n)$ and $H(z)=2e^{-t}(t+z^n)$.
For $z\in\Delta$ we have $\operatorname{Re}H(z)\geqslant2e^{-t}(t-|z|^n)>0$, that is, $H\in C^n$.
Since $z_j^n=-1$, we have
        $\operatorname{Re}H(z_j)
        =
        2e^{-t}(t-1)
        >
        0$,
        $j=1,\ldots,n$.
This means that equalities~\eqref{eq:Stupin:extrPolynomial_phi} do not hold, whereas
inequality~\eqref{eq:Stupin:scalar_prod} holds since $H\in C^n$.
\end{proof}

\section{Properties of coefficient derivatives}

\begin{Lemma}
    \label{lem:Stupin:f_n'}
Let $c$ be a smooth complex-valued function of $t\in\mathbb{R}$, and let $c(t^*)>0$ at some point $t=t^*$.
Then $\operatorname{Re}c'(t^*)=0$ if and only if $|c|'(t^*)=0$.
\end{Lemma}

\begin{proof}
Let $u\coloneqq\operatorname{Re}c$, $v\coloneqq\operatorname{Im}c$.
Then $c=u+iv$, $|c|=\sqrt{u^2+v^2}$, and
    \begin{equation}
        \label{eq:Stupin:|c|'}
        |c|'=\frac{uu'+vv'}{|c|},
        \quad
        c\neq0.
    \end{equation}

Since by assumption $c(t^*)>0$, we have $v(t^*)=0$ and $u(t^*)=|c(t^*)|>0$.
Therefore~\eqref{eq:Stupin:|c|'} gives $|c|'(t^*)=u'(t^*)=\operatorname{Re}c'(t^*)$.
Hence the equalities $\operatorname{Re}c'(t^*)=0$ and $|c|'(t^*)=0$ are equivalent.
\end{proof}

\begin{Lemma}
    \label{lem:Stupin:scalarProd_f_n'}
Let $h$ be holomorphic in $\Delta$, let $t\in\mathbb R$ and $f(z,t)\coloneqq e^{-th(z)}$; then
    \begin{equation*}
        \{f\}_n'
        =
        -(\{f\}_n\{h\}_0+\ldots+\{f\}_0\{h\}_n),
        \quad
        n\in\mathbb{N}\cup\{0\}.
    \end{equation*}
\end{Lemma}

\begin{proof}
Fix $n\in\mathbb{N}\cup\{0\}$ and $0<r<1$.
By Cauchy's formula,
    \[
        \{f\}_n'(t)
        =
        \Bigg(
        \frac1{2\pi i}\int\limits_{|z|=r}
        \frac{f(z,t)}{z^{n+1}}\,dz
        \Bigg)'
        =
        \frac{1}{2\pi i}\int\limits_{|z|=r}
        \frac{f_t'(z,t)}{z^{n+1}}\,dz
        =
        \{f_t'\}_n.
    \]
The second equality clearly holds by the theorem on differentiation under the integral sign with respect to a parameter.
Furthermore, $f_t'(z,t)=-h(z)f(z,t)$, and formula~\eqref{eq:Stupin:prodCoeff} gives the required conclusion.
\end{proof}

\smallskip
Lemmas~\ref{lem:Stupin:scalarProd_f_n'} and~\ref{lem:Stupin:f_n'} immediately imply the following criterion:

\begin{Theorem}
    \label{th:Stupin:max_f_n'}
Let $n\in\mathbb{N}$, $h\in C_1$, $f(z,t)=e^{-th(z)}$, and let $\{f\}_n(t^*)>0$ at some point $t=t^*$.
Then $t=t^*$ is a critical point of the function $|\{f\}_n(t)|$ if and only if
    $\operatorname{Re}(\{f\}_n(t^*)\{h\}_0+\ldots+\{f\}_0(t^*)\{h\}_n)=0$.
\end{Theorem}

In particular, for an extremal function, Theorem~\ref{th:Stupin:max_f_n'} gives
equality~\eqref{eq:Stupin:orthogonalVects}.
Theorems~\ref{th:Stupin:scalar_prod_criteria} and~\ref{th:Stupin:max_f_n'} show that
equality~\eqref{eq:Stupin:orthogonalVects} is not a characteristic property of extremal functions.
For example, let $n=1$, $h(z)\coloneqq1-z\in C_1$, and $f(z)\coloneqq e^{-h(z)}=e^{-1+z}\in B_1$.
Then $\{f\}_0=\{f\}_1=e^{-1}$, $\{h\}_0=1$, $\{h\}_1=-1$, and
$\{f\}_1\{h\}_0+\{f\}_0\{h\}_1=0$.
However, the function $f$ is not locally extremal.
Indeed, take $g(z)\coloneqq\frac{1-z}{1+z}\in C_1$ and consider the variation $f_\varepsilon\coloneqq
e^{-h_\varepsilon}$ of the function $f$, where
$h_\varepsilon\coloneqq(1-\varepsilon)h+\varepsilon g$, $0<\varepsilon<1$.
Then $h_\varepsilon\in C_1$, $f_\varepsilon\in B_1$, and
$\{f_\varepsilon\}_1=e^{-1}(1+\varepsilon)>e^{-1}=\{f\}_1$.
Moreover, $t_\varepsilon=t_h=1$, so for any $\delta>0$ and $0<\eta\leqslant1$, for $0<\varepsilon<\eta$ we have
$|t_\varepsilon-t_h|=0<\delta$, $h_\varepsilon\in U_\eta(h)$, and $f_\varepsilon\in V_{\delta,\eta}(f)$.
Hence every neighborhood of the function $f$ in the topology $\tau_{B_+}$ contains functions $f_\varepsilon$ such that
$\{f_\varepsilon\}_1>\{f\}_1$, and therefore $f$ is not locally extremal.

Apparently, the most promising properties are statement~\eqref{eq:Stupin:extrPolynomial} and
equalities~\eqref{eq:Stupin:extrPolynomial_phi}.

\section{The class $C^n$}

In this section we assume that $n\in\mathbb{N}$, $h_0>0$, and the polynomial $H\in C^n$ has the form
\begin{equation}
    \label{eq:Stupin:polynomFromC}
    H(z)
    \coloneqq
    h_0 + 2\sum\limits_{k=1}^n h_k z^k.
\end{equation}
The polynomial $H$ given by formula~\eqref{eq:Stupin:extrPolynomialH_n} is a special case of the polynomial $H$ given by
formula~\eqref{eq:Stupin:polynomFromC}, so the same notation is used without causing confusion.
Every function $f\in\mathfrak E(n)$ satisfies $\{f\}_0>0$ and $\{f\}_n>0$, so we are primarily interested here in
polynomials $H$ with $h_n>0$ and $h_0>0$.

Let $n\in\mathbb{N}$, $a_0,a_k,b_k\in\mathbb{R}$, $k=1,\ldots,n$, $a_n^2+b_n^2>0$.
The function
\begin{equation*}
    T(\varphi)
    \coloneqq
    a_0 + \sum\limits_{k=1}^n (a_k\cos{k\varphi} - b_k\sin{k\varphi})
\end{equation*}
of the real argument $\varphi$ will be called a trigonometric polynomial of degree $n$.

Note that $T$ and $H$ are related.
Indeed, if $H$ is a polynomial of degree $n$, the restriction of $\operatorname{Re} H$ to the unit circle
$\partial\Delta$ is a trigonometric polynomial of degree $n$,

with

    $h_0\coloneqq a_0\in\mathbb{R}$,

    $h_k\coloneqq\frac{a_k+ib_k}2$,

    $a_k,b_k\in\mathbb{R}$,

    $k=1,\ldots,n$.

The converse also holds (see~\cite{Stupin:ref16}).

\smallskip
A Laurent polynomial is an expression of the form $Q(z)=\sum\limits_{k=-n_1}^{n_2}q_k z^k$, where
$n_1,n_2\in\mathbb{N}\cup\{0\}$.
We state the assertion on the complex form of a trigonometric polynomial from~\cite{Stupin:ref16}.

\begin{Assertion}
    \label{st:Stupin:T_complex}
An arbitrary trigonometric polynomial of degree at most $n$ admits the representation
        $T(\varphi)
        \coloneqq
        a_0 + \sum\limits_{k=1}^n (a_k\cos{k\varphi} - b_k\sin{k\varphi})$
if and only if it can be written in the form
        $T(\varphi)=\sum\limits_{k=-n}^n h_k e^{ik\varphi}$
and the following conditions hold:
    \begin{equation}
        \label{eq:Stupin:T_complex_case}
        h_{-k}=\overline{h_k},
        \quad
        k=0,\ldots,n,
    \end{equation}
where
        $h_0\coloneqq a_0\in\mathbb{R}$,
        $h_k\coloneqq\frac{a_k+ib_k}2$,
        $a_k,b_k\in\mathbb{R}$,
        $k=1,\ldots,n$.%
\end{Assertion}

The following assertion was proved in~\cite{Stupin:ref16} in the proof of the complex version of the Fejer-Riesz
theorem.

\begin{Assertion}
    \label{st:Stupin:LaurentBoundaryZerosEven}
Let the Laurent polynomial $Q\not\equiv0$ and let $Q(z)\geqslant0$ for $|z|=1$.
Then each of its roots on the unit circle has even multiplicity.
\end{Assertion}

\smallskip
We need the following assertion (see~\cite{Martin:ref1, Stupin:ref16}):

\begin{Assertion}
    \label{st:Stupin:trigPolynomRootNumber}
A trigonometric polynomial $T$ of degree $n$ has at most $2n$ roots on $[0, 2\pi)$, counted with multiplicity.
\end{Assertion}

The article~\cite{Rovnyak:ref1} notes that in the early twentieth century L.~Fejer~\cite{Fejer:ref1} first drew
attention to the importance of the class of trigonometric polynomials taking nonnegative values on the entire real line.
His conjecture on the structure of these polynomials was proved by F.~Riesz~\cite{Riesz:ref1} and is now known as the
Fejer-Riesz theorem.
A proof can be found in~\cite{Hussen:ref1}, in~\cite[p.~154]{Tsuji:ref1}, or in~\cite{Stupin:ref16}.

\begin{Theorem}[Riesz, Fejer]
    \label{th:Stupin:FR_complex_case}
If $n\in\mathbb{N}$, $T$ is a trigonometric polynomial of degree $n$, and $T(\varphi)\geqslant0$ for all
$\varphi\in\mathbb{R}$, then there exists a polynomial $P(z)\coloneqq\sum\limits_{j=0}^n p_j z^j$, $p_n\neq0$, such that
$T(\varphi)=|P(e^{i\varphi})|^2$, $\varphi\in\mathbb{R}$.
Moreover, $P$ can be chosen so that all its zeros lie outside the unit disk.
In this case $P$ is uniquely determined up to a unimodular multiplicative constant.
\end{Theorem}

The criterion below is an analogue of the Caratheodory-Toeplitz criterion~\ref{th:Stupin:nCTEquiv} for polynomials and
follows from the Fejer-Riesz theorem~\ref{th:Stupin:FR_complex_case} (see~\cite{Stupin:ref16}):

\begin{Theorem}
    \label{th:Stupin:characterization_of_polynomials_from_C}
Let $n\in\mathbb{N}$.
The polynomial $H(z) \coloneqq h_0 + 2\sum\limits_{k=1}^n h_k z^k$, where $h_0>0$, $h_k\in\mathbb{C}$, has positive real
part in $\Delta$ if and only if there exist numbers $p_j\in\mathbb{C}$, not all zero, such that
    \begin{equation}
        \label{eq:Stupin:characterization_of_polynomials_from_C}
        h_k
        =
        \sum\limits_{j=0}^{n-k} p_{j+k}\overline{p_j},
        \quad
        k=0,\ldots,n.
    \end{equation}
Moreover, if $h_k\in\mathbb{R}$, $k=0,\ldots,n$, there exist $p_j\in\mathbb{R}$, $j=0,\ldots,n$ such that
equalities~\eqref{eq:Stupin:characterization_of_polynomials_from_C} hold\textnormal{;}
if $p_j\in\mathbb{R}$, $j=0,\ldots,n$, then $h_k\in\mathbb{R}$, $k=0,\ldots,n$.
\end{Theorem}

The following identity was obtained in~\cite{Stupin:ref16} in the proof of this criterion.

\begin{Assertion}
    \label{st:Stupin:polynomialModulusIdentity}
Let $n\in\mathbb N$, $P(z)\coloneqq\sum\limits_{j=0}^n p_j z^j$, and
    $H(z)\coloneqq h_0+2\sum\limits_{k=1}^n h_kz^k$,
where $p_j\in\mathbb C$, $j=0,\ldots,n$, and the numbers $h_k$, $k=0,\ldots,n$ are defined by
formulas~\eqref{eq:Stupin:characterization_of_polynomials_from_C}.
Then
    \[
        \operatorname{Re}H(e^{i\varphi})=|P(e^{i\varphi})|^2,
        \quad \varphi\in\mathbb R.
    \]
\end{Assertion}

\begin{proof}
We have
    \[
        |P(e^{i\varphi})|^2
        =\sum_{j=0}^n\sum_{\ell=0}^n p_j\overline{p_\ell}e^{i(j-\ell)\varphi}
        =
        \sum_{j=0}^n|p_j|^2
        +2\operatorname{Re}\sum_{k=1}^n
        \sum_{j=0}^{n-k}p_{j+k}\overline{p_j}e^{ik\varphi}
        =\operatorname{Re}H(e^{i\varphi}).
    \]
This proves the assertion.
\end{proof}

Note that for a fixed index $n$ the Caratheodory-Toeplitz criterion~\ref{th:Stupin:nCTEquiv} does not determine whether
a holomorphic function $h$ belongs to the class $C$.
For each index $k$, this criterion determines only whether the point $h^{(k+1)}$ lies on the boundary of the coefficient
body $C^{(k+1)}$, inside it, or outside it.
If $h$ is holomorphic in $\Delta$, then membership of $h$ in the class $C$ is equivalent to $h^{(k+1)}\in C^{(k+1)}$ for
all $k\in\mathbb{N}$.
If $h^{(k+1)}\in\operatorname{int}C^{(k+1)}$ for all $k\in\mathbb{N}$, then $h$ may be called an interior function of
the class $C$.
Theorem~\ref{th:Stupin:characterization_of_polynomials_from_C} determines whether a polynomial $H$ belongs to the
Caratheodory class.
By the Caratheodory-Toeplitz criterion, if $H\in C^n$, then $H^{(k+1)}\in\operatorname{int}C^{(k+1)}$, $k\in\mathbb{N}$,
that is, if $H\in C^n$, then $H$ is an interior function of the class $C$.
Nevertheless, an interior polynomial can be defined as follows: if $\operatorname{Re}H(z)>0$, $z\in\partial\Delta$, then
$H$ is an interior polynomial.

\smallskip
All assertions stated in this section without proof are taken from~\cite{Stupin:ref16} and follow from
Theorem~\ref{th:Stupin:characterization_of_polynomials_from_C}.
These assertions are also adapted to the case $h_n>0$ under study.

Recall that $C^n$ is the set of polynomials $H$ of degree $n$ with $\operatorname{Re}{H(z)}>0$, $z\in\Delta$,

and the class $C_t^n$ is the subset of $C^n$ normalized by $H(0)=t$, where $t=\operatorname{const}>0$.

\begin{Theorem}
    \label{th:Stupin:uniqunessInC}
There exists a unique polynomial $H$ of the form~\eqref{eq:Stupin:polynomFromC} satisfying the conditions\textnormal:
    $H\in C_{h_0}^n$, $h_0=2h_n>0$.
Moreover, $H(z)=h_0(1+z^n)$.
\end{Theorem}

\begin{Theorem}
    \label{th:Stupin:unic}
There exists a unique polynomial $H$ of the form~\eqref{eq:Stupin:polynomFromC} satisfying the conditions\textnormal:
    $H\in C_{h_0}^n$,
    $h_n>0$,
all roots $\zeta_1,\ldots,\zeta_n$ of the polynomial $H$ lie on the unit circle $\partial\Delta$.
Moreover, $H(z)=h_0(1+z^n)$.
\end{Theorem}

Note that the roots $\zeta_k$ of the polynomial $H(z)=h_0(1+z^n)$ are precisely all the $n$th roots of $-1$, that is,
$\zeta_k=e^{i\frac{\pi+2\pi k}n}$, $k=0,\ldots,n-1$.

\begin{Assertion}
    \label{st:Stupin:zero}
Let $m\in\mathbb{N}$, $m\leqslant n$, $z_1,\ldots,z_m\in\partial\Delta$, and let $H$ and $P$ be polynomials of degree
$n$ such that $\operatorname{Re}H(e^{i\varphi})=|P(e^{i\varphi})|^2$, $\varphi\in\mathbb{R}$.
Then the condition $\operatorname{Re}H(z_j)=0$ is equivalent to $P(z_j)=0$ for all $j=1,\ldots,m$.
\end{Assertion}

As mentioned above, the function $T(\varphi)\coloneqq\operatorname{Re} H(e^{i\varphi})$ is a trigonometric polynomial,
and by the Fejer-Riesz theorem (Theorem~\ref{th:Stupin:FR_complex_case}), if $H\in C^n$, there exists a polynomial
$P(z)\coloneqq\sum\limits_{j=0}^n p_j z^j$ of degree $n$ with roots outside $\Delta$ such that $\operatorname{Re}
H(z)=|P(z)|^2$, $|z|=1$.

A polynomial is called normalized (or monic) if its leading coefficient equals $1$.
Along with the polynomials $H$ and $P$, it is convenient to introduce their normalized versions: the polynomial
$\mathcal H(z)\coloneqq \eta_0+\eta_1 z+\ldots+\eta_n z^n$ defined by

    $H=2h_n\mathcal H$

and the polynomial $\Psi(z)\coloneqq\psi_0+\psi_1 z+\ldots+\psi_n z^n$ defined by $P=p_n\Psi$.
Clearly, $\eta_n=\psi_n=1$.
If $\zeta_1,\ldots,\zeta_n$ are the roots of $\mathcal H$, then $\eta_0\coloneqq(-1)^n \zeta_1\cdot\ldots\cdot\zeta_n$.
Similarly, if $z_1,\ldots,z_n$ are the roots of $\Psi$, then $\psi_0\coloneqq(-1)^n z_1\cdot\ldots\cdot z_n$.

A polynomial $P$ will be called external for a polynomial $H$ if $\operatorname{Re} H(z)=|P(z)|^2$ for $|z|=1$ and all
roots of $P$ lie outside the unit disk.
A polynomial $\Psi$ will be called normalized external for a polynomial $H$ if $\operatorname{Re}
H(z)=|p_n|^2\cdot|\Psi(z)|^2$ for $|z|=1$, $\psi_n=1$, and all roots of $\Psi$ lie outside the unit disk.
The term ``external'' for the polynomial $P$ or $\Psi$ is chosen solely because their roots lie outside the unit disk.

\begin{Theorem}
    \label{th:Stupin:h_0_eq_h_n_norm_Psi2}
If $n\in\mathbb{N}$, the polynomial $H$ is defined by formula~\eqref{eq:Stupin:polynomFromC}, $h_n>0$, and $H\in C$,
then\textnormal:
there exists a normalized external polynomial $\Psi(z)\coloneqq\psi_0+\psi_1 z+\ldots+\psi_n z^n$ of degree
$n$\textnormal;
the following equalities hold:
    \begin{equation}
        \label{eq:Stupin:eta_0}
        \eta_0
        \coloneqq
        \frac{h_0}{2h_n}
        =
        \frac{|\psi_0|^2+\ldots+|\psi_n|^2}{2\psi_0};
    \end{equation}
the following sharp inequalities hold\textnormal:
    \begin{equation}
        \label{eq:Stupin:rootsOfH_prod_arg}
        h_0
        \geqslant
        2h_n,
        \quad
        \eta_0
        \geqslant
        1,
        \quad
        \psi_0
        \geqslant
        1;
    \end{equation}
the equality $h_0=2h_n$ is equivalent to $\eta_0=1$, which in turn is equivalent to the equalities
    $\psi_0=1\land\psi_1=0\land\ldots\land\psi_{n-1}=0$\textnormal;
    \newline
all these equalities hold if and only if $H(z)=2h_n(1+z^n)$\textnormal;
    \newline
the equality $\eta_0=1$ implies $\eta_1=0\land\ldots\land\eta_{n-1}=0$\textnormal;
    \newline
Here $\mathcal H(z)\coloneqq H(z)/2h_n=\eta_0+\eta_1 z+\ldots+\eta_n z^n$.
\end{Theorem}

\smallskip
Theorem~\ref{th:Stupin:h_0_eq_h_n_norm_Psi2} implies

\begin{Corollary}
    \label{col:Stupin:max_h_n}
Among all polynomials $H\in C_{h_0}^n$ of the form~\eqref{eq:Stupin:polynomFromC} with $h_n>0$, the maximum value $h_n$
is attained only by the polynomial $H(z)=2h_n(1+z^n)$.
\end{Corollary}

Note that inequality~\eqref{eq:Stupin:f_n_ge_2f_0} also follows from Corollary~\ref{col:Stupin:max_h_n}.

\begin{Theorem}
    \label{th:Stupin:nullsOfP_and_H}
If $n\in\mathbb{N}$, the polynomial $H$ is defined by formula~\eqref{eq:Stupin:polynomFromC} with $h_n>0$ and $H\in C$,
then\textnormal:
there exists a normalized external polynomial $\Psi(z)\coloneqq\psi_0+\psi_1 z+\ldots+\psi_n z^n$ of degree
$n$\textnormal;
the roots $\zeta_1,\ldots,\zeta_n$ of the polynomial $H$ satisfy the sharp inequalities
    \begin{equation}
        \label{eq:Stupin:rootsOfH_abs}
        |\zeta_k|\geqslant1,
        \quad
        k=1,\ldots,n;
    \end{equation}
equality in all the inequalities simultaneously holds if and only if $\psi_0=1$ and all intermediate coefficients
$\psi_1,\ldots,\psi_{n-1}$ vanish
    \textnormal($H(z)=2h_n(1+z^n)$\textnormal)\textnormal;
    \newline
the following equalities hold:
    \begin{equation}
        \label{eq:Stupin:rootsOfH}
        h_0=2|p_n|^2\eta_0\psi_0,
        \quad
        h_n=|p_n|^2\psi_0,
        \quad
\text{where}
        \quad
        \eta_0\coloneqq(-1)^n\zeta_1\cdot\ldots\cdot\zeta_n.
    \end{equation}
\end{Theorem}

\smallskip
Note that if $\psi_0=1$, then $h_0=2|p_n|^2\eta_0$ and $h_n=|p_n|^2$ by Theorem~\ref{th:Stupin:nullsOfP_and_H}.

\begin{Theorem}
    \label{th:Stupin:nullsOfP_against_H}
If $n\in\mathbb{N}$, the polynomial $H$ is defined by formula~\eqref{eq:Stupin:polynomFromC}, $h_n>0$ and $H\in C$, and
$\zeta_1,\ldots,\zeta_n$ are the roots of $H$, then\textnormal:
    \newline
there exists a normalized external polynomial $\Psi(z)\coloneqq\psi_0+\psi_1 z+\ldots+\psi_n z^n$ of degree
$n$\textnormal;
    \newline
the following sharp inequalities hold:
    \[
        1
        \leqslant
        \frac12
        \left(
            \psi_0+\frac1{\psi_0}
        \right)
        \leqslant
        \eta_0,
        \quad
\text{where}
        \quad
        \eta_0\coloneqq(-1)^n\zeta_1\cdot\ldots\cdot\zeta_n;
    \]
equality in the first inequality holds if and only if $\psi_0=1$\textnormal;
    \newline
equality in the second inequality holds if and only if all intermediate coefficients $\psi_1,\ldots,\psi_{n-1}$
vanish\textnormal;
    \newline
equality in both inequalities simultaneously holds if and only if $\psi_0=1$ and all intermediate coefficients
$\psi_1,\ldots,\psi_{n-1}$ vanish
    \textnormal($H(z)=2h_n(1+z^n)$\textnormal)\textnormal;
    \newline
\end{Theorem}

\begin{Assertion}
If $n\in\mathbb{N}$, the polynomial $H$ of degree $n$ is defined by formula~\eqref{eq:Stupin:polynomFromC}, and $H\in
C$, then there exists an external polynomial $P(z)$ of degree $n$, and $H$ and $P$ cannot have common roots outside the
unit circle $\partial\Delta$.
\end{Assertion}

Note that the polynomials $H$ and $P$ may have common roots on the unit circle $|z|=1$.
For example, the polynomials $P(z)\coloneqq1+z$ and $H(z)\coloneqq2(1+z)$ have the same roots, and $\operatorname{Re}
H(e^{i\varphi})=|P(e^{i\varphi})|^2$ for $\varphi\in\mathbb{R}$.

\begin{Theorem}
    \label{th:Stupin:estimOnC_n}
Let $n\in\mathbb{N}$ and $H(z) \coloneqq h_0 + 2\sum\limits_{k=1}^n h_k z^k$.
If $H\in C_1$, then
    \[
        |h_k|
        \leqslant
        \cos\frac{\pi}{\lfloor\frac{n}k\rfloor + 2},
        \quad k=1,\ldots,n.
    \]
Equality is attainable for every fixed $k\in\{1,\ldots,n\}$.

\end{Theorem}

The precise equality conditions in the estimation of Theorem~\ref{th:Stupin:estimOnC_n} are considered
in~\cite{Stupin:ref16}.
Here we mention only a necessary condition for equality: there exists a point $z_0\in\partial\Delta$ such that
$\operatorname{Re}H(z_0)=0$.
This is the extremal property of polynomials in the class $C$ whose real part has zeros on $\partial\Delta$.

\smallskip
Theorem~\ref{th:Stupin:estimOnC_n} implies (see~\cite{Stupin:ref16})

\begin{Corollary}
    \label{cor:Stupin:max_h_n_C_1}
If $h\in C_1$, then $|\{h\}_n|\leqslant2$ for every positive integer $n$.
\end{Corollary}

Note that Corollary~\ref{cor:Stupin:max_h_n_C_1} gives a sharp estimation, that is, the equality $|\{h\}_n|=2$ is
attainable (see~\cite[p.~57]{Rogosinski:ref1}).
Corollary~\ref{cor:Stupin:max_h_n_C_1} is well known (see, for example,~\cite[p.~57]{Rogosinski:ref1}).
The functions for which $|\{h\}_j|=2$ are also listed in~\cite{Rogosinski:ref1}.
Note that none of them are polynomials.

\section{Further consequences}

The following assertion follows immediately from Theorem~\ref{th:Stupin:characterization_of_polynomials_from_C} and the
definition of a function of extremal type:

\begin{Corollary}
    \label{cor:Stupin:extrCoeff}
If $n\in\mathbb{N}$, $f\in\mathbb E(n)$, then there exist $p_j\in\mathbb{C}$, $j=0,\ldots,n$ such that
    \begin{equation}
        \label{eq:Stupin:f_n_p_n}
        \{f\}_{n-k}
        =
        \sum\limits_{j=0}^{n-k} p_{j+k}\overline{p_j},
        \quad
        k=0,\ldots,n.
    \end{equation}
If, in addition, $f\in B^r$, there exist $p_j\in\mathbb{R}$, $j=0,\ldots,n$ satisfying~\eqref{eq:Stupin:f_n_p_n}.
\end{Corollary}

By Theorem~\ref{th:Stupin:localExtremalIsExtremalType}, Corollary~\ref{cor:Stupin:extrCoeff} holds, in particular, for
every function $f\in\mathfrak E(n)$.

Note that formulas~\eqref{eq:Stupin:f_n_p_n}, like formulas~\eqref{eq:Stupin:characterization_of_polynomials_from_C},
are invariant under multiplication of all $p_j$ by a unimodular constant.
Thus we may assume that $p_n>0$.

\smallskip
Let $n\in\mathbb{N}$.
If $f\in\mathbb E(n)$ and $H$ is the polynomial defined by formula~\eqref{eq:Stupin:extrPolynomialH_n}, then the first
equality~\eqref{eq:Stupin:eta_0}, Corollary~\ref{cor:Stupin:extrCoeff}, and equalities~\eqref{eq:Stupin:rootsOfH}
immediately imply that
\begin{equation}
    \label{eq:Stupin:roots}
    \{f\}_n
    =
    |p_0|^2+\ldots+|p_n|^2
    =
    2\{f\}_0\eta_0
    =
    2|p_n|^2\eta_0\psi_0,
    \quad
    \{f\}_0
    =
    |p_n|^2\psi_0,
\end{equation}
where $\eta_0\coloneqq(-1)^n\zeta_1\cdot\ldots\cdot\zeta_n$ and $\zeta_1,\ldots,\zeta_n$ are the roots of $H$, and
$\psi_0\coloneqq(-1)^n z_1\cdot\ldots\cdot z_n$ and $z_1,\ldots,z_n$ are the roots of $P$.
Since for a function of extremal type $\{f\}_n>0$ and $\{f\}_0>0$, relations~\eqref{eq:Stupin:rootsOfH_abs}
and~\eqref{eq:Stupin:rootsOfH_prod_arg} also hold.
In particular, all these relations hold for every function $f\in\mathfrak E(n)$.

As mentioned in Section~\ref{sec:Stupin:review}, the Krzyz conjecture is proved for all $n=1,\ldots,6$, so for these $n$
we know the unique normalized globally extremal function $f(z)=F(-z^n,1)$, where $F$ is given by
formula~\eqref{eq:Stupin:BExtrF}.
Thus,
\[
    H(z)=\tfrac2e(1+z^n), \quad n=1,\ldots,6,
\]
and the roots $H$ are all the $n$th roots of $-1$.

\smallskip
Theorem~\ref{th:Stupin:h_0_eq_h_n_norm_Psi2} and Corollary~\ref{cor:Stupin:extrCoeff} immediately imply

\begin{Corollary}
    \label{cor:Stupin:f_n_eq_f_0_norm_Psi2}
If $n\in\mathbb{N}$, $f\in\mathbb E(n)$, and the polynomial $H$ is generated by the coefficients of the function $f$
according to formula~\eqref{eq:Stupin:extrPolynomialH_n}, then there exists a normalized external polynomial
$\Psi(z)\coloneqq\psi_0+\psi_1 z+\ldots+\psi_n z^n$ and the following relations hold\textnormal:
    \[
        \frac{\{f\}_n}{\{f\}_0}
        =
        \frac{|\psi_0|^2+\ldots+|\psi_n|^2}{\psi_0}
        \geqslant
        2,
    \]
with equality in the inequality $\{f\}_n\geqslant2\{f\}_0$ attained if and only if $\psi_0=1$,
$\psi_1=\ldots=\psi_{n-1}=0$, that is, if and only if $H(z)=2\{f\}_0(1+z^n)$.
\end{Corollary}

By Theorem~\ref{th:Stupin:localExtremalIsExtremalType}, Corollary~\ref{cor:Stupin:f_n_eq_f_0_norm_Psi2} also holds for
every function $f\in\mathfrak E(n)$.

For an arbitrary polynomial $S(z)=s_0+\ldots+s_n z^n$, define the reciprocal polynomial
$
S^*(z)\coloneqq\overline{s_0} z^n + \ldots + \overline{s_n}.
$
Clearly, $S^*(z)=z^n\overline{S\left(1/\overline z\right)}$, $z\neq0$, with the value at $z=0$ defined by continuity.
The polynomial $S$ is called self-reciprocal if the identity $S(z)=S^*(z)$ holds.

Consider the case $\psi_0=1$ in Corollary~\ref{cor:Stupin:f_n_eq_f_0_norm_Psi2}.
Since $|\psi_0|^2+\ldots+|\psi_n|^2=\|\Psi\|_{H^2}^2$, where $\|\cdot\|_{H^2}$ is the norm in the Hardy space $H^2$, we
have $\{f\}_n=\|\Psi\|_{H^2}^2\{f\}_0$ and $2\eta_0=\|\Psi\|_{H^2}^2$, with $\|\Psi\|_{H^2}^2\geqslant2$.
Note also that the condition $\psi_0=1$ is equivalent to self-reciprocity of the polynomial $\Psi$, since precisely in
this case $\psi_0=\psi_n$ and all roots of $\Psi$ lie on the unit circle.
Since $|\eta_0|\geqslant1$, Theorem~\ref{th:Stupin:unic} shows that the polynomial $\mathcal H$ is self-reciprocal if
and only if $\eta_0=1$.

Note that Corollary~\ref{cor:Stupin:f_n_eq_f_0_norm_Psi2} immediately implies inequality~\eqref{eq:Stupin:f_n_ge_2f_0},
that is, $\{f\}_n\geqslant2\{f\}_0$.

\section{On the poles of the function $h=-\ln f$, $f\in\mathbb E(n)$}

Let $f=e^{-h}$, $h\in C$.
Differentiating this equality with respect to $z$ gives $f'=-h'f$.
Substituting the Taylor expansions of the functions $f$ and $h$ into this equality gives the formula
\begin{equation}
    \label{eq:Stupin:deForBCoeff}
    \ell\{f\}_\ell
    =
    -\sum\limits_{k=1}^\ell k\{h\}_k\{f\}_{\ell-k},
    \quad
    \ell\in\mathbb N,
\end{equation}
relating the coefficients of the functions $f$ and $h$.

\smallskip
Let $n\in\mathbb N$, $f\in\mathbb E(n)$, and let the polynomial $H$ be generated by the coefficients of the function $f$
according to formula~\eqref{eq:Stupin:extrPolynomialH_n}.
That is, $H(z)=\{f\}_n+2\{f\}_{n-1}z+\ldots+2\{f\}_0 z^n$.
Consider the Laurent polynomial
\begin{equation}
    \label{eq:Stupin:zerosPolesLaurentPolynomial}
    \mathcal T(z)
    \coloneqq
    \{f\}_n+
    \sum\limits_{\ell=1}^n
    \big(\{f\}_{n-\ell}z^\ell+\overline{\{f\}_{n-\ell}}z^{-\ell}\big),
    \quad
    z\neq0.
\end{equation}
Clearly,
\begin{equation}
    \label{eq:Stupin:zerosPolesBoundaryValues}
    \mathcal T(z)=\operatorname{Re}H(z),\quad |z|=1.
\end{equation}
Thus $\mathcal T$ extends the function $\operatorname{Re}H$ from $\partial\Delta$ to $\mathbb C\setminus\{0\}$.

By the definition of a function of extremal type, the function $h=-\ln f$ has the form
\begin{equation}
    \label{eq:Stupin:zerosPolesRepresentationH}
    h(z)
    =
    \sum\limits_{j=1}^m\alpha_j\frac{1+z_j z}{1-z_j z},
    \quad
    m\leqslant n,
    ~
    \alpha_j>0,
    ~
    z_j\coloneqq e^{i\varphi_j},
    ~
    0\leqslant\varphi_1<\ldots<\varphi_m<2\pi.
\end{equation}
Define the rational function $L(z)\coloneqq\frac1z h'\big(\frac1z\big)$.
Note that
\begin{equation}
    \label{eq:Stupin:zerosPolesRationalFunction}
    L(z)
    =
    2\sum\limits_{j=1}^m\alpha_j\frac{z_jz}{(z-z_j)^2}.
\end{equation}
The following expansion also holds:
\begin{equation}
    \label{eq:Stupin:L}
    L(z)
    =
    \sum\limits_{k=1}^{\infty}k\{h\}_k z^{-k},
    \quad
    |z|>1.
\end{equation}

Here we show that each zero $z_j\in\partial\Delta$ of the function $\operatorname{Re}H$ corresponds to exactly one pole
$\overline{z_j}$ of the function $h$.
Conversely, each pole $\overline{z_j}\in\partial\Delta$ of the function $h$ corresponds to exactly one zero $z_j$ of the
function $\operatorname{Re}H$.
Moreover, every zero of the function $\varphi\mapsto\operatorname{Re}H(e^{i\varphi})$ has multiplicity two.

To this end, we establish a relation between the functions $\operatorname{Re}H(z)$, $|z|=1$, and
$\operatorname{Re}(zH'(z))$, $|z|=1$ using formula~\eqref{eq:Stupin:deForBCoeff}.
For $|z|>1$ the coefficient of $z^{-k}$ in the expansion $L$ equals $k\{h\}_k$.
Hence the coefficient of $z^\ell$, $0\leqslant\ell\leqslant n-1$, in the function $L\mathcal T$ equals the negative of
the sum on the right-hand side of this formula for the coefficient with index $n-\ell$, that is,
$-(n-\ell)\{f\}_{n-\ell}$.
This allows the function $L\mathcal T$ on the unit circle to be expressed through $\operatorname{Re}H(z)$ and
$\operatorname{Re}(zH'(z))$.

\begin{Lemma}
    \label{lem:Stupin:zerosPolesLaurentProduct}
Let $n\in\mathbb N$, $f\in\mathbb E(n)$, let the polynomial $H$ be generated by the coefficients of the function $f$
according to formula~\eqref{eq:Stupin:extrPolynomialH_n}, let the function $\mathcal T$ be given by
formula~\eqref{eq:Stupin:zerosPolesLaurentPolynomial}, and the function $L$ by
formula~\eqref{eq:Stupin:zerosPolesRationalFunction}.
Then
    \begin{equation}
        \label{eq:Stupin:zerosPolesLaurentIdentity}
        L(z)\mathcal T(z)
        =
        -n\{f\}_n-
        \sum\limits_{\ell=1}^{n-1}(n-\ell)
        \big(\{f\}_{n-\ell}z^\ell+\overline{\{f\}_{n-\ell}}z^{-\ell}\big),
        \quad
        z\neq0.
    \end{equation}
\end{Lemma}

\begin{proof}
    \textbf{1.}
We show that at every point $z_j$ the function $\mathcal T$ has a zero of order at least two.
By the definition of a function of extremal type, $H\in C^n$ and $\operatorname{Re}H(z_j)=0$.
By continuity, $\operatorname{Re}H(z)\geqslant0$ for $|z|=1$.
Thus $\mathcal T$ is nonnegative on the circle, with $\mathcal T(z_j)=0$, $j=1,\ldots,m$.
It is not identically zero, since its constant coefficient $\{f\}_n>0$.
By Assertion~\ref{st:Stupin:LaurentBoundaryZerosEven}, the multiplicity of every root $z_j$ is even and at least two.

    \textbf{2.}
The function $L$ has a second-order pole at $z_j$.
By formula~\eqref{eq:Stupin:zerosPolesRationalFunction}, in a neighborhood of $z_j$ we have
    \[
        L(z)
        =
        L\big(z_j+(z-z_j)\big)
        =
        \frac{2\alpha_jz_j^2}{(z-z_j)^2}+\frac{2\alpha_jz_j}{z-z_j}+B_j(z),
    \]
where $B_j$ is holomorphic in this neighborhood, since the points $z_1,\ldots,z_m$ are pairwise distinct.
Since $2\alpha_jz_j^2\neq0$, the point $z_j$ is a second-order pole of $L$.

    \textbf{3.}
We show that the singularities of $L\mathcal T$ at the points $z_j$ are removable.
Since at each point $z_j$ the function $L$ has a second-order pole and $\mathcal T$ a zero of order at least two, the
singularities of the product $L\mathcal T$ at these points are removable.
This product has no other possible finite poles except zero.

    \textbf{4.}
We justify formula~\eqref{eq:Stupin:zerosPolesLaurentIdentity}.
For $z\to\infty$ we have $L(z)=O(z^{-1})$, and for $z\to0$ we have $L(z)=O(z)$.
By definition, $\mathcal T(z)=O(z^n)$ for $z\to\infty$ and $\mathcal T(z)=O(z^{-n})$ for $z\to0$.
Hence $L(z)\mathcal T(z)=O(z^{n-1})$, $z\to\infty$, and $L(z)\mathcal T(z)=O(z^{-(n-1)})$, $z\to0$.

Thus the rational function $z^{n-1}L(z)\mathcal T(z)$ has no finite poles.
A rational function without finite poles is a polynomial.
The estimation as $z\to\infty$ shows that this polynomial has degree at most $2n-2$.
Hence $L\mathcal T$ is a Laurent polynomial in powers of $z$ from $-(n-1)$ to $n-1$.

Write $z_j=e^{i\varphi_j}$.
For $z=e^{i\varphi}\notin\{z_1,\ldots,z_m\}$ the identity $e^{i\varphi}-e^{i\varphi_j}=2i
e^{i\frac{\varphi+\varphi_j}2}\sin\frac{\varphi-\varphi_j}{2}$ gives
        $L(e^{i\varphi})
        =
        -\sum\limits_{j=1}^m\frac{\alpha_j}{2\sin^2\frac{\varphi-\varphi_j}2}
        \in\mathbb R$.

Write $L(z)\mathcal T(z)=\sum\limits_{\ell=-(n-1)}^{n-1}c_\ell z^\ell$.
On the unit circle, the function $L\mathcal T$ is real.
Applying Assertion~\ref{st:Stupin:T_complex} to the trigonometric polynomial $(L\mathcal T)(e^{i\varphi})$ of degree at
most $n$ with $c_{-n}=c_n=0$, formula~\eqref{eq:Stupin:T_complex_case} gives the equalities
$c_{-\ell}=\overline{c_\ell}$, $\ell=0,\ldots,n-1$.

    \textbf{5.}
We calculate the coefficients of the product $L\mathcal T$ for formula~\eqref{eq:Stupin:zerosPolesLaurentIdentity}.
By formula~\eqref{eq:Stupin:L}, we have $L(z)=\sum\limits_{k=1}^{\infty}k\{h\}_k z^{-k}$, $|z|>1$.
This series converges absolutely and locally uniformly for $|z|>1$, so it can be multiplied by the finite
sum~\eqref{eq:Stupin:zerosPolesLaurentPolynomial}.
For $\ell\geqslant0$ the power $z^\ell$ arises only from the products $k\{h\}_k z^{-k}$ and $\{f\}_{n-q}z^q$ with
$q-k=\ell$.
The conditions $1\leqslant q\leqslant n$ and $k\geqslant1$ give $1\leqslant k\leqslant n-\ell$.
For $\ell=0,\ldots,n-1$ the coefficient of $z^\ell$ in the product of~\eqref{eq:Stupin:zerosPolesLaurentPolynomial}
and~\eqref{eq:Stupin:L} equals
        $\sum\limits_{k=1}^{n-\ell}k\{h\}_k\{f\}_{n-\ell-k}
        =
        -(n-\ell)\{f\}_{n-\ell}$
by formula~\eqref{eq:Stupin:deForBCoeff} for the coefficient with index $n-\ell$.
For $\ell=0$ we obtain the real coefficient $-n\{f\}_n$.
The coefficients of negative powers are determined by conjugation.
This proves equality~\eqref{eq:Stupin:zerosPolesLaurentIdentity}.
\end{proof}

\begin{Lemma}
    \label{lem:Stupin:zerosPolesLaurentIdentity}
Let $n\in\mathbb N$, $f\in\mathbb E(n)$, let the polynomial $H$ be generated by the coefficients of the function $f$
according to formula~\eqref{eq:Stupin:extrPolynomialH_n}, let the function $\mathcal T$ be given by
formula~\eqref{eq:Stupin:zerosPolesLaurentPolynomial}, and the function $L$ by
formula~\eqref{eq:Stupin:zerosPolesRationalFunction}.
Then
    \begin{equation}
        \label{eq:Stupin:zerosPolesRadialIdentity}
        \operatorname{Re}\bigl(zH'(z)\bigr)
        =
        \bigl(n+L(z)\bigr)\operatorname{Re}H(z),
        \quad |z|=1,\quad z\notin\{z_1,\ldots,z_m\}.
    \end{equation}
\end{Lemma}

\begin{proof}
Let $|z|=1$ and $z\notin\{z_1,\ldots,z_m\}$.
Differentiating polynomial~\eqref{eq:Stupin:extrPolynomialH_n} and multiplying the result by $z$ gives
    \[
        zH'(z)=2\sum\limits_{\ell=1}^n\ell\{f\}_{n-\ell}z^\ell.
    \]
Since $|z|=1$, we have $\overline z=z^{-1}$ and
    $\overline{\{f\}_{n-\ell}z^\ell}=\overline{\{f\}_{n-\ell}}z^{-\ell}$.
Therefore, using the equality $\operatorname{Re}\big(zH'(z)\big)=\frac12\big(zH'(z)+\overline{zH'(z)}\big)$, we obtain
    \[
        \operatorname{Re}\bigl(zH'(z)\bigr)
        =
        \sum\limits_{\ell=1}^n \ell
        \big(\{f\}_{n-\ell}z^\ell+\overline{\{f\}_{n-\ell}}z^{-\ell}\big).
    \]

On the other hand, multiply~\eqref{eq:Stupin:zerosPolesLaurentPolynomial} by $n$ and
add~\eqref{eq:Stupin:zerosPolesLaurentIdentity}.
The constant terms $n\{f\}_n$ and $-n\{f\}_n$ cancel.
For $1\leqslant\ell\leqslant n-1$ the factor for each pair of opposite powers equals
    $n-(n-\ell)=\ell$.
The terms with powers $z^n$ and $z^{-n}$ occur only in $n\mathcal T(z)$ and have the factor $n$.
Hence,
    \[
        n\mathcal T(z)+L(z)\mathcal T(z)
        =
        \sum\limits_{\ell=1}^n \ell
        \left(\{f\}_{n-\ell}z^\ell+\overline{\{f\}_{n-\ell}}z^{-\ell}\right)
        =
        \operatorname{Re}\bigl(zH'(z)\bigr).
    \]
For $n=1$ the sum in~\eqref{eq:Stupin:zerosPolesLaurentIdentity} is empty, and after canceling the constant terms only
the pair of powers $z$ and $z^{-1}$ remains.

Factoring out $\mathcal T(z)$ on the left-hand side of the preceding equality and
using~\eqref{eq:Stupin:zerosPolesBoundaryValues}, we obtain
    \[
        \operatorname{Re}\bigl(zH'(z)\bigr)
        =
        n\mathcal T(z)+L(z)\mathcal T(z)
        =
        \bigl(n+L(z)\bigr)\operatorname{Re}H(z).
    \]
This proves formula~\eqref{eq:Stupin:zerosPolesRadialIdentity}.
\end{proof}

Recall Harnack's inequality for the disk (see~\cite[p.~47]{AxlerBourdonRamey:ref1}).

\begin{Theorem}[Harnack]
    \label{th:Stupin:HarnackDisk}
If the function $u$ is positive and harmonic in $\Delta$, then
    \[
        \tfrac{1-|z|}{1+|z|}u(0)
        \leqslant
        u(z)
        \leqslant
        \tfrac{1+|z|}{1-|z|}u(0),
        \quad z\in\Delta.
    \]
\end{Theorem}

\begin{Lemma}
    \label{lem:Stupin:boundaryZeroRadialDerivative}
Let $n\in\mathbb N$, $H\in C^n$, $z_0\in\partial\Delta$, and $\operatorname{Re}H(z_0)=0$.
Then
    \[
        \operatorname{Re}\bigl(z_0 H'(z_0)\bigr)
        \leqslant
        -\tfrac12H(0)
        <0.
    \]
\end{Lemma}

\begin{proof}
The function $u(z)\coloneqq\operatorname{Re}H(z)$ is positive and harmonic in $\Delta$, with $u(0)=H(0)>0$.
By Theorem~\ref{th:Stupin:HarnackDisk}, for $0<r<1$ we have $u(rz_0)\geqslant\frac{1-r}{1+r}H(0)$.
Subtracting $u(z_0)=0$ from the left-hand side of this inequality and then dividing it by $r-1<0$ gives
$\frac{u(rz_0)-u(z_0)}{r-1}\leqslant-\frac{H(0)}{1+r}$.
Since $r-1<0$ and $z_0\neq0$, we have
    \[
        \frac{u(rz_0)-u(z_0)}{r-1}
        =
        \operatorname{Re}\frac{H(rz_0)-H(z_0)}{r-1}
        =
        \operatorname{Re}\left(z_0\frac{H(rz_0)-H(z_0)}{rz_0-z_0}\right).
    \]
The polynomial $H$ is holomorphic at $z_0$, so passing to the limit as $r\to1-$ gives the required inequality.
\end{proof}

\begin{Theorem}
    \label{th:Stupin:zerosOfPAndPolesOfH}
Let $n\in\mathbb N$, $f\in\mathbb E(n)$, let the polynomial $H$ be generated by the coefficients of the function $f$
according to formula~\eqref{eq:Stupin:extrPolynomialH_n}, and let $P$ be an external polynomial for $H$.
A point $\overline{z_0}\in\partial\Delta$ is a pole of the function $h=-\ln f$ if and only if $z_0$ is a root of the
polynomial $P$.
Moreover, every root of $P$ lying on $\partial\Delta$ is simple.
\end{Theorem}

\begin{proof}
Fix $n\in\mathbb N$.
By the definition of a function of extremal type, $H\in C^n$.
Thus by continuity $\operatorname{Re}H(e^{i\varphi})\geqslant0$, $\varphi\in\mathbb R$.
By the Fejer-Riesz theorem~\ref{th:Stupin:FR_complex_case}, there exists a polynomial $P$ external for $H$.
Also, by the definition of a function of extremal type, the function $h=-\ln f$ has
representation~\eqref{eq:Stupin:zerosPolesRepresentationH}.

\textbf{Necessity.}
Let the point $\overline{z_0}$ be a pole of the function $h$.
Formula~\eqref{eq:Stupin:zerosPolesRepresentationH} defining the function $h$ implies that $z_0=z_j$ for some
$j\in\{1,\ldots,m\}$.
By the definition of a function of extremal type, the condition $\operatorname{Re}H(z_j)=0$ holds.
By Assertion~\ref{st:Stupin:zero}, it is equivalent to $P(z_j)=0$.
Hence $P(z_0)=0$.

\textbf{Sufficiency.}
Let $z_0\in\partial\Delta$ and $P(z_0)=0$.
Consider the functions $L$ and $\mathcal T$ given by formulas~\eqref{eq:Stupin:zerosPolesRationalFunction}
and~\eqref{eq:Stupin:zerosPolesLaurentPolynomial}.

    \textbf{1.}
First exclude additional zeros on $\partial\Delta$.
To this end, suppose $z_0\in\partial\Delta\setminus\{z_1,\ldots,z_m\}$ and $P(z_0)=0$.
By Assertion~\ref{st:Stupin:zero}, $P(z_0)=0$ implies $\operatorname{Re}H(z_0)=0$.
The function $L$ is finite at $z_0$ by formula~\eqref{eq:Stupin:zerosPolesRationalFunction}.
Formula~\eqref{eq:Stupin:zerosPolesRadialIdentity} gives
    $\operatorname{Re}\bigl(z_0H'(z_0)\bigr)=0$,
which contradicts Lemma~\ref{lem:Stupin:boundaryZeroRadialDerivative}.

Hence $z_0=z_j$ for some $j\in\{1,\ldots,m\}$.
In the representation of the function $h$ by the sum~\eqref{eq:Stupin:zerosPolesRepresentationH}, only the $j$th term
has a pole at $\overline{z_j}$, and its numerator there equals $2\alpha_j>0$.
Thus $\overline{z_0}$ is a pole of the function $h$.
Sufficiency is proved.

    \textbf{2.}
We establish simplicity of the boundary roots of $P$.
By formula~\eqref{eq:Stupin:zerosPolesLaurentIdentity}, the product $L\mathcal T$ is a Laurent polynomial in powers of
$z$ and hence is holomorphic at every point $z_j\neq0$.
By formulas~\eqref{eq:Stupin:zerosPolesRadialIdentity} and~\eqref{eq:Stupin:zerosPolesBoundaryValues}, we have
    \[
        \operatorname{Re}\bigl(zH'(z)\bigr)
        =
        \bigl(n+L(z)\bigr)\operatorname{Re}H(z)
        =
        n\mathcal T(z)+(L\mathcal T)(z),
        \quad |z|=1,\quad z\notin\{z_1,\ldots,z_m\}.
    \]
Since the first of these equalities is established only on the unit circle and away from the points $z_1,\ldots,z_m$, we
pass to the limit in it as $z\to z_j$, $z\neq z_j$ along this circle.
The polynomial $H$ and the Laurent polynomials $\mathcal T$ and $L\mathcal T$ are continuous at $z_j\neq0$.
Since~\eqref{eq:Stupin:zerosPolesBoundaryValues} implies $n\mathcal T(z_j)=n\operatorname{Re}H(z_j)=0$, we have
$(L\mathcal T)(z_j)=\operatorname{Re}\bigl(z_jH'(z_j)\bigr)<0$ by Lemma~\ref{lem:Stupin:boundaryZeroRadialDerivative}.
By formula~\eqref{eq:Stupin:zerosPolesRationalFunction}, the function $L$ has a second-order pole at $z_j$, since the
points $z_1,\ldots,z_m$ are pairwise distinct and $\alpha_j>0$.
Since the Laurent polynomial $L\mathcal T$ is holomorphic and nonzero at $z_j$, the zero of $\mathcal T$ there has order
two.
Since $\big.\frac{d}{d\varphi}e^{i\varphi}\big|_{\varphi=\varphi_j}=iz_j\neq0$,
formula~\eqref{eq:Stupin:zerosPolesBoundaryValues} shows that the function
$\varphi\mapsto\operatorname{Re}H(e^{i\varphi})$ has a zero of multiplicity two at $\varphi=\varphi_j$.
By the definition of an external polynomial, $\operatorname{Re}H(e^{i\varphi})=|P(e^{i\varphi})|^2$.
This equality shows that the multiplicity of this zero is twice the multiplicity of the root $z_j$ of $P$.
Hence every boundary root of $P$ is simple.
\end{proof}

\begin{Corollary}
    \label{cor:Stupin:extremalTypeZerosOfH}
Let $n\in\mathbb N$, $f\in\mathbb E(n)$, and let the polynomial $H$ be generated by the coefficients of the function $f$
according to formula~\eqref{eq:Stupin:extrPolynomialH_n}.
Then the points $z_j$, $j=1,\ldots,m$, in representation~\eqref{eq:Stupin:zerosPolesRepresentationH} of the function
$h=-\ln f$ exhaust all zeros of the function $\operatorname{Re}H$ on the unit circle $\partial\Delta$.
The number $q$ of distinct zeros equals $m\leqslant n$, and every corresponding zero of the function
$\varphi\mapsto\operatorname{Re}H(e^{i\varphi})$ has multiplicity two.
If $P$ is an external polynomial for $H$, then
    \[
        P(z)
        =
        \prod_{j=1}^m(z-z_j)R(z),
        \quad
        \deg R=n-m,
    \]
where all roots of the polynomial $R$ satisfy the strict inequality $|z|>1$.
\end{Corollary}

\begin{proof}
By Theorem~\ref{th:Stupin:zerosOfPAndPolesOfH}, the points $z_j$, $j=1,\ldots,m$, exhaust all roots of the polynomial
$P$ lying on the unit circle, and these roots are all simple.
By the definition of an external polynomial, $\operatorname{Re}H(e^{i\varphi})=|P(e^{i\varphi})|^2$.
This equality gives the statements on the number and multiplicities of the zeros.
Since $\deg P=n$, division by $\prod\limits_{j=1}^m(z-z_j)$ gives a polynomial $R$ of degree $n-m$.
It has no roots on the unit circle, and none inside it by the definition of an external polynomial.
\end{proof}

Let $\mathcal E(n)$ denote the subset of $\mathcal{B}(n)$ consisting only of functions of extremal type.
Note that $\mathcal E(n)$ is nonempty, since $F(-z^n,1)\in\mathcal E(n)$.

For $m=n$, Theorem~\ref{th:Stupin:zerosOfPAndPolesOfH} and Corollary~\ref{cor:Stupin:extremalTypeZerosOfH} give the
following corollary.

\begin{Corollary}
Let $n\in\mathbb N$, $f\in\mathcal E(n)$, let the polynomial $H$ be generated by the coefficients of the function $f$
according to formula~\eqref{eq:Stupin:extrPolynomialH_n}, and let $P$ be an external polynomial for $H$.
Then the points $z_j$, $j=1,\ldots,n$, exhaust all zeros of the function $\operatorname{Re}H$ on the unit circle, and
every corresponding zero of the function $\varphi\mapsto\operatorname{Re}H(e^{i\varphi})$ has multiplicity two.
All roots of $P$ are simple and lie on $\partial\Delta$.
A point $z_0\in\mathbb C$ is a root of the polynomial $P$ if and only if $\overline{z_0}$ is a pole of the function
$h=-\ln f$.
If $\Psi$ is a normalized external polynomial for $H$ with constant term $\psi_0$, then
    \[
        \Psi(z)=\prod_{j=1}^n(z-z_j),
        \quad
        \psi_0=(-1)^nz_1\cdots z_n=1.
    \]
\end{Corollary}

\begin{proof}
For $m=n$, the polynomial $R$ in Corollary~\ref{cor:Stupin:extremalTypeZerosOfH} is a nonzero constant.
Hence all roots of $P$ lie on the circle, and the remaining assertions on roots and poles follow from
Theorem~\ref{th:Stupin:zerosOfPAndPolesOfH}.
Since the roots of $\Psi$ are $z_1,\ldots,z_n$ and its leading coefficient is $1$, we have
$\Psi(z)=\prod\limits_{j=1}^n(z-z_j)$, whence $|\psi_0|=1$.
On the other hand, $\psi_0\geqslant1$ by~\eqref{eq:Stupin:rootsOfH_prod_arg}, so $\psi_0=(-1)^n z_1\cdots z_n=1$.
\end{proof}

\smallskip
For example, if $f(z)\coloneqq F(-z^n,1)$, then
$h(z)=\frac1n\sum\limits_{k=1}^n\frac{1+e^{i(2k-1)\pi/n}z}{1-e^{i(2k-1)\pi/n}z}
=\frac{1-z^n}{1+z^n}$.
Accordingly, $H(z)=\frac2e(1+z^n)$, and the function $\operatorname{Re}{H}$ has $n$ zeros on $\partial\Delta$ of
multiplicity $2$.

\section{On uniqueness of an extremal function}

Let $n\in\mathbb{N}$.
As mentioned above, the class $B$ and the functional $|\{f\}_n|$ are invariant under rotations in the variable planes
$z$ and $w$ ($w=f(z)$).
Hence, if $f$ is an extremal function in the Krzyz problem for the coefficient with index $n$, we may assume without
loss of generality that $\{f\}_0>0$ and $\{f\}_n>0$.
By Corollary~\ref{cor:Stupin:B_t_extremal_Th}, under this normalization $f\in\mathfrak E(n)$.
Below, uniqueness of an extremal function is considered on the larger set $\mathbb E(n)$.

\begin{Theorem}
    \label{th:Stupin:uniqunessOfExtremalFunctionInB}
Let $n\in\mathbb{N}$.
There exists a unique function $f$ in $\mathbb E(n)$ satisfying $\{f\}_n=2\{f\}_0$, and $f(z)=F(-z^n,1)$, where the
function $F$ is given by formula~\eqref{eq:Stupin:BExtrF}.
\end{Theorem}

\begin{proof}
Local extremality of the function \[F(-z^n,1)=\tfrac1e+\tfrac2e z^n+o(z^n)\] was established first in~\cite{Hummel:ref1}
and later also in~\cite{Stupin:ref8}.
Hence $F(-z^n,1)\in\mathfrak E(n)$.
By Theorem~\ref{th:Stupin:localExtremalIsExtremalType}, the function $F(-z^n,1)\in\mathbb E(n)$, with
$\{F(-z^n,1)\}_n=2\{F(-z^n,1)\}_0$.
We prove its uniqueness.
Fix $n\in\mathbb N$ and a function $f\in\mathbb E(n)$ satisfying $\{f\}_n=2\{f\}_0$.
By the definition of a function of extremal type, the polynomial $H$ generated by the coefficients of the function $f$
according to formula~\eqref{eq:Stupin:extrPolynomialH_n} belongs to the class $C^n$.
Together with $\{f\}_n=2\{f\}_0$, this means that the polynomial $H$ satisfies all the conditions of
Theorem~\ref{th:Stupin:uniqunessInC}.
Hence,
        $H(z)=\{f\}_n+2\{f\}_0 z^n$,
which implies the equalities
        $\{f\}_1=\ldots=\{f\}_{n-1}=0$.

Let $\{f\}_0=e^{-t}$, where $t>0$.
Since by assumption $\{f\}_n=2\{f\}_0$, we have
    \[
        f(z)=e^{-t}+2e^{-t}z^n+o(z^n).
    \]
Hence $h\coloneqq -\ln f=t-2z^n+o(z^n)$.
By the definition of a function of extremal type, $h\in C$.
To apply the Caratheodory-Toeplitz criterion to establish uniqueness of the function $h$ for $t>0$, calculate the minors
$M_k$ defined by formula~\eqref{eq:Stupin:C-T_minors}:
    \begin{itemize}
        \item $M_0=2t>0$;
\item $M_k=2^{k+1}t^{k+1}>0$ for $k=1,\ldots,n-1$;
\item $M_n=2^{n+1}t^{n-1}(t^2-1)\geqslant0$ and $M_n=0$ if and only if $t=1$.
    \end{itemize}
By the definition of a function of extremal type, the function $h$ satisfies the equality conditions in
Theorem~\ref{th:Stupin:scalar_prod_criteria}.
Hence this theorem implies equality~\eqref{eq:Stupin:orthogonalVects}, which gives
    \[
        0
        =
        \operatorname{Re}{(\{f\}_n\{h\}_0+\{f\}_{n-1}\{h\}_1+\ldots+\{f\}_0\{h\}_n)}
        =
        2e^{-t}t-2e^{-t},
    \]
whence $t=1$.

Thus, for the function $h(z)=1-2z^n+o(z^n)$, the minors $M_0,\ldots,M_{n-1}$ are positive, and $M_n=0$.
Therefore, by the Caratheodory-Toeplitz criterion, the extension of the polynomial $1-2z^n$ to a function in the class
$C$ is unique.
The function
    $
        \frac{1-z^n}{1+z^n}
        =
        1-2z^n+o(z^n)
    $
belongs to the class $C$, so $h(z)=\tfrac{1-z^n}{1+z^n}$.
Thus $f(z)=F(-z^n,1)$.
\end{proof}

\smallskip
By Theorem~\ref{th:Stupin:localExtremalIsExtremalType}, Theorem~\ref{th:Stupin:uniqunessOfExtremalFunctionInB} also
holds for functions $f\in\mathfrak E(n)$.

\smallskip
Let $f\in\mathbb E(n)$ and $f_k(z)\coloneqq f(e^{\frac{2\pi ki}{n}}z)$, $k=0,\ldots,n-1$.
Then $f_k\in\mathbb E(n)$, $k=0,\ldots,n-1$, that is, the class $\mathbb E(n)$ is invariant under rotations through
multiples of $\frac{2\pi}n$ in the variable plane $z$.
Below we therefore call the functions $f_k$, $k=0,\ldots,n-1$ admissible rotations of the function $f$ (in the variable
plane $z$).

\begin{Theorem}
    \label{th:Stupin:uniqunessInE}
Let $n\in\mathbb{N}$.
There exists a unique function $f$ in $\mathbb E(n)$ invariant under admissible rotations in the variable plane $z$, and
$f(z)=F(-z^n,1)$, where the function $F$ is given by formula~\eqref{eq:Stupin:BExtrF}.
\end{Theorem}

\begin{proof}
By Theorem~\ref{th:Stupin:uniqunessOfExtremalFunctionInB}, the function $F(-z^n,1)$ belongs to $\mathbb E(n)$.
Clearly, this function is invariant under admissible rotations in the variable plane $z$.

Suppose the function $f\in\mathbb E(n)$ is invariant under admissible rotations in the variable plane $z$.
For $n>1$, invariance under rotation of the argument $z$ of the function $f$ through the angle $\frac{2\pi}n$ implies
that
    $e^{\frac{2\pi ki}{n}}\{f\}_k=\{f\}_k$, $k=1,\ldots,n-1$.
Since $e^{\frac{2\pi ki}{n}}\neq1$ for all $k=1,\ldots,n-1$, we have
    $\{f\}_1=\ldots=\{f\}_{n-1}=0$.
For $n=1$ there are no intermediate coefficients.
Thus, for every $n\in\mathbb N$, the polynomial $H$ generated by the coefficients of the function $f$ according to
formula~\eqref{eq:Stupin:extrPolynomialH_n} has the form $H(z)=\{f\}_n+2\{f\}_0 z^n$.
By Corollary~\ref{cor:Stupin:f_n_eq_f_0_norm_Psi2}, we have $\{f\}_n\geqslant2\{f\}_0$.
On the other hand, by the definition of a function of extremal type, there exists a point $z_j\in\partial\Delta$ for
which
    $\operatorname{Re}H(z_j)=0$.
Hence,
    $0=\{f\}_n+2\{f\}_0\operatorname{Re}{z_j^n}\geqslant\{f\}_n-2\{f\}_0$,
whence $\{f\}_n\leqslant2\{f\}_0$.
Thus $\{f\}_n=2\{f\}_0$, and $f(z)=F(-z^n,1)$ by Theorem~\ref{th:Stupin:uniqunessOfExtremalFunctionInB}.
\end{proof}

\smallskip
The following assertion is a consequence of Theorem~\ref{th:Stupin:localExtremalIsExtremalType}:

\begin{Theorem}
    \label{th:Stupin:Krzyz1}
The set $\mathfrak E(1)$ consists of a single function $F(-z,1)$.
In particular, the Krzyz conjecture holds for $n=1$.
\end{Theorem}

\begin{proof}
Let $f\in\mathfrak E(1)$ and $h\coloneqq-\ln f$.
Then formula~\eqref{eq:Stupin:orthogonalVects}, since $h\in C$, gives
    $
        \operatorname{Re}{(\{f\}_1\{h\}_0+\{f\}_0\{h\}_1)}
        =
        0$,
whence $\{f\}_1\{h\}_0=-\{f\}_0\operatorname{Re}{\{h\}_1}$.
Formula~\eqref{eq:Stupin:f_n_ge_2f_0} gives $\{f\}_1\geqslant2\{f\}_0$, but $|\{h\}_1|\leqslant2\{h\}_0$ by
Corollary~\ref{cor:Stupin:max_h_n_C_1}, whence $\{f\}_1=2\{f\}_0$.
Theorem~\ref{th:Stupin:uniqunessOfExtremalFunctionInB} implies that $f(z)=F(-z,1)$ and that no other such functions
exist.

Since every global extremal, up to rotations in the variable planes $z$ and $w$ ($w=f(z)$), belongs to $\mathfrak E(1)$,
and this set consists of a single function $f(z)=F(-z,1)$, the maximum of $|\{f\}_1|$ on the class $B$ equals $2/e$ and
is attained by exactly one function $f(z)=F(-z,1)$.
Hence the Krzyz conjecture holds for $n=1$.
\end{proof}

\smallskip
\begin{Remark}
The conditions $f\in\mathfrak E(n)$ and $\{f\}_n=2/e$ for a global extremal $f$ are not equivalent to the Krzyz
conjecture for index $n$.
One must also show that the equality $\{f\}_n=2/e$ is possible only for the function $f(z)=F(-z^n,1)$.
Theorem~\ref{th:Stupin:uniqunessOfExtremalFunctionInB} ensures this uniqueness only under the additional condition
$\{f\}_n=2\{f\}_0$, which does not follow from the equality $\{f\}_n=2/e$.
Theorem~\ref{th:Stupin:f_n_ge_2f_0} implies only that $\{f\}_0\leqslant1/e$, with equality in this inequality possible
if and only if $f(z)=F(-z^n,1)$ by Theorem~\ref{th:Stupin:uniqunessOfExtremalFunctionInB}.
\end{Remark}

\smallskip
By Theorem~\ref{th:Stupin:uniqunessOfExtremalFunctionInB}, if $f\in\mathbb E(n)$ and $\{f\}_n=2\{f\}_0$, then
$\{f\}_k=0$ for $k=1,\ldots,n-1$.

\begin{Theorem}
    \label{th:Stupin:converseToUniqunessInB}
Let $n\in\mathbb{N}\setminus\{1\}$.
There exists a unique function $f\in\mathbb E(n)$ satisfying $\{f\}_k=0$ for $k=1,\ldots,n-1$, and $f(z)=F(-z^n,1)$.
\end{Theorem}

\begin{proof}
Fix $n>1$.
By Theorem~\ref{th:Stupin:uniqunessOfExtremalFunctionInB}, the function $F(-z^n,1)$ belongs to $\mathbb E(n)$.
Since its coefficients with indices $k=1,\ldots,n-1$ vanish, existence is proved.

Now let $f\in\mathbb E(n)$ and $\{f\}_k=0$ for $k=1,\ldots,n-1$.
Set $h\coloneqq-\ln f$.
By the definition of a function of extremal type, $h\in C$, so $h/\{h\}_0\in C_1$.
By Corollary~\ref{cor:Stupin:f_n_eq_f_0_norm_Psi2}, we have $\{f\}_n\geqslant 2\{f\}_0>0$.
By the definition of a function of extremal type, the function $h$ satisfies the equality conditions in
Theorem~\ref{th:Stupin:scalar_prod_criteria}.
Hence this theorem implies equality~\eqref{eq:Stupin:orthogonalVects}.
The condition $\{f\}_k=0$ for $k=1,\ldots,n-1$ and equality~\eqref{eq:Stupin:orthogonalVects} imply
$\{f\}_n\{h\}_0=-\{f\}_0\operatorname{Re}{\{h\}_n}$, and hence
    \[
        \{f\}_n=-\{f\}_0\frac{\operatorname{Re}{\{h\}_n}}{\{h\}_0}\leqslant2\{f\}_0,
    \]
since $|\{h\}_n/\{h\}_0|\leqslant2$ for $h/\{h\}_0 \in C_1$.
(See~\cite[p.~57]{Rogosinski:ref1} or Corollary~\ref{cor:Stupin:max_h_n_C_1}.)
Thus $\{f\}_n=2\{f\}_0$, whence Theorem~\ref{th:Stupin:uniqunessOfExtremalFunctionInB} implies $f(z)=F(-z^n,1)$.
Uniqueness is proved.
\end{proof}

By Theorem~\ref{th:Stupin:localExtremalIsExtremalType}, Theorem~\ref{th:Stupin:converseToUniqunessInB} also holds for
functions $f\in\mathfrak E(n)$.

We show that uniqueness of a global extremal in $\mathfrak E(n)$ implies the Krzyz conjecture.

\begin{Theorem}
    \label{th:Stupin:uniqueGlobalExtremalImpliesKrzyz}
Let $n\in\mathbb{N}$.
If there exists a unique globally extremal function in $\mathfrak E(n)$, then the Krzyz conjecture holds for index $n$.
\end{Theorem}

\begin{proof}
For $n=1$, the assertion follows from Theorem~\ref{th:Stupin:Krzyz1}.

Fix $n\in\mathbb{N}\setminus\{1\}$ and denote by $f$ the unique global extremal in $\mathfrak E(n)$.
We have
    \[
        f(z)
        =
        \{f\}_0+\sum\limits_{k=1}^{n-1}\{f\}_k z^k+\{f\}_n z^n+\ldots
    \]
The function
    \[
        g(z)
        \coloneqq
        f(e^{2i\pi/n}z)
        =
        \{f\}_0+\sum\limits_{k=1}^{n-1}\{f\}_k e^{2i\pi k/n}z^k+\{f\}_n z^n+\ldots
    \]
is also a global extremal and is normalized by $\{g\}_0>0$, $\{g\}_n>0$.
Hence $g\in\mathfrak E(n)$, and the uniqueness assumption gives $g=f$.
Moreover, $\{g\}_0=\{f\}_0$, $\{g\}_k=e^{2i\pi k/n}\{f\}_k$, $k=1,\ldots,n-1$, $\{g\}_n=\{f\}_n$.
Hence,
    \[
        e^{2i\pi k/n}\{f\}_k=\{f\}_k,
        \quad
        k=1,\ldots,n-1,
    \]
which is possible if and only if $\{f\}_k=0$, $k=1,\ldots,n-1$.
By Theorem~\ref{th:Stupin:converseToUniqunessInB}, we have $f(z)=F(-z^n,1)$.
Since $f$ is a global extremal,
    \[\max\limits_{g\in B}|\{g\}_n|=\{f\}_n=2/e.\]

Now let $\mathfrak f\in B$ and $|\{\mathfrak f\}_n|=2/e$.
Then $\mathfrak f$ is a global extremal.
Normalizing it by rotations in the variable planes $z$ and $w$ gives the function
    $\mathfrak g\in\mathfrak E(n)$.
By the uniqueness assumption, $\mathfrak g(z)=F(-z^n,1)$.
Returning to the function $\mathfrak f$, we conclude that it has the form $e^{i\theta}F(e^{i\varphi}z^n,1)$, where
$\theta,\varphi\in\mathbb R$.
Thus the Krzyz conjecture holds for index $n$.
\end{proof}

\smallskip
In particular, the following assertion holds:

\begin{Assertion}
If $n\in\mathbb{N}$, $f\in\mathfrak E(n)$, and $f$ is the unique global extremal in $\mathfrak E(n)$, then all its
coefficients are real.
\end{Assertion}

\begin{proof}
Fix $n\in\mathbb{N}$ such that the function $f$ is the unique global extremal normalized by $\{f\}_0>0$, $\{f\}_n>0$ in
the Krzyz problem for index $n$.
If the function $f$ has at least one nonreal coefficient, the function $f^*(z)\coloneqq\overline{f(\overline z)}$ is
also a global extremal.
Moreover, $\{f^*\}_0=\{f\}_0>0$ and $\{f^*\}_n=\{f\}_n>0$, so $f^*\in\mathfrak E(n)$.
The presence of a nonreal coefficient means $f^*\neq f$, contradicting uniqueness of the global extremal in $\mathfrak
E(n)$.
\end{proof}

\begin{Theorem}
Let $n\in\mathbb{N}$.
If for every globally extremal function $f\in\mathfrak E(n)$ all roots $\zeta_1,\ldots,\zeta_n$ of the polynomial $H$
given by formula~\eqref{eq:Stupin:extrPolynomialH_n} belong to $\partial\Delta$, then the Krzyz conjecture holds for
index $n$.
\end{Theorem}

\begin{proof}
Fix $n\in\mathbb{N}$ satisfying the hypothesis of the theorem and an arbitrary globally extremal function
    $f\in\mathfrak E(n)$.
Since all roots of the polynomial $H$ lie on the unit circle, $|\eta_0|=1$.
Equalities~\eqref{eq:Stupin:roots} and the conditions $\{f\}_n>0$, $\{f\}_0>0$ imply that
    $\eta_0=\{f\}_n/(2\{f\}_0)>0$.
Hence $\eta_0=1$ and $\{f\}_n=2\{f\}_0$.
By Theorem~\ref{th:Stupin:uniqunessOfExtremalFunctionInB}, we have $f(z)=F(-z^n,1)$.
Thus $F(-z^n,1)$ is the unique global extremal in $\mathfrak E(n)$, and
    $\max\limits_{g\in B}|\{g\}_n|=2/e$.

Now let $\mathfrak f\in B$ and $|\{\mathfrak f\}_n|=2/e$.
Then $\mathfrak f$ is a global extremal.
Normalizing it by rotations in the variable planes $z$ and $w$ gives the function
    $\mathfrak g\in\mathfrak E(n)$.
Hence $\mathfrak g(z)=F(-z^n,1)$.
Returning to the function $\mathfrak f$, we conclude that it has the form $e^{i\theta}F(e^{i\varphi}z^n,1)$, where
$\theta,\varphi\in\mathbb R$.
Thus the Krzyz conjecture holds for index $n$.
\end{proof}

We show how to reduce verification of the Krzyz conjecture to finding a function, not necessarily extremal, in the
classes $B_t$, $t\geqslant1$ (the Krzyz conjecture states that $t=1$).

\begin{Theorem}
    \label{th:Stupin:existenceInB}
Let $n\in\mathbb{N}\setminus\{1\}$.
The Krzyz conjecture for index $n$ holds if and only if the classes $B_t$, for $t\geqslant1$, contain no function
$\mathfrak f$ such that
    \[
        \{\mathfrak f\}_n=2e^{-1}
\quad\text{and}\quad
        |\{\mathfrak f\}_1|+\ldots+|\{\mathfrak f\}_{n-1}|>0.
    \]
\end{Theorem}

\begin{proof}
Fix an arbitrary index $n\in\mathbb{N}\setminus\{1\}$.

\textbf{Necessity.}
If the Krzyz conjecture for index $n$ holds, then uniqueness of the global extremal in $\mathfrak E(n)$ clearly implies
that no function with the indicated properties exists.

\textbf{Sufficiency.}
Suppose the Krzyz conjecture for index $n$ is false.
Then there exists a global extremal $f\in\mathfrak E(n)$ distinct from $F(-z^n,1)$.
Indeed, if $m_n>\frac2e$, we may take any normalized global extremal.
If $m_n=\frac2e$, falsity of the conjecture implies the existence of an equality-case function not having the form
specified in the conjecture; after normalization it belongs to $\mathfrak E(n)$ and is distinct from $F(-z^n,1)$.
Since $f$ is a global extremal, $\{f\}_n=m_n\geqslant\frac2e$.
Choose a positive constant $c\leqslant1$ so that $c\{f\}_n=\frac2e$, and set $\mathfrak f\coloneqq cf$.
Then $\mathfrak f\in B$ and $\{\mathfrak f\}_n=\frac2e$.
By inequality~\eqref{eq:Stupin:f_n_ge_2f_0}, we have $\{f\}_n\geqslant2\{f\}_0$.
Equality here is impossible by Theorem~\ref{th:Stupin:uniqunessOfExtremalFunctionInB}, so
    $
        \{\mathfrak f\}_0=c\{f\}_0<e^{-1}
    $.
Hence $\mathfrak f\in B_t$ for some $t>1$.
If $\{\mathfrak f\}_1=\ldots=\{\mathfrak f\}_{n-1}=0$, then $c>0$ would immediately imply that
    $\{f\}_1=\ldots=\{f\}_{n-1}=0$.
Then $f(z)=F(-z^n,1)$ by Theorem~\ref{th:Stupin:converseToUniqunessInB}, contradicting the choice of the extremal $f$ at
the beginning of this proof.
Thus,
    $
        |\{\mathfrak f\}_1|+\ldots+|\{\mathfrak f\}_{n-1}|>0
    $.
\end{proof}

\smallskip
The following assertion connects Theorems~\ref{th:Stupin:existenceInB} and~\ref{th:Stupin:denseExtr} with algorithmic
search in the class $\mathcal B(n)$:

\begin{Assertion}
Let $n\in\mathbb{N}$.
If $m_n>2/e$, then the Krzyz conjecture is false for index $n$ and the class $\mathcal B(n)$ contains a function $q$
such that $\{q\}_n>2/e$.
\end{Assertion}

\begin{proof}
If $m_n>2/e$, there exists an extremal function $f\in\mathfrak E(n)$ for which
    $\{f\}_n>2/e$.
Hence the Krzyz conjecture is false for index $n$.
By the definition of $\mathfrak E(n)$, we have $f\in\mathfrak B(n)$.
By Theorem~\ref{th:Stupin:denseExtr}, the set $\mathcal B(n)$ is dense in $\mathfrak B(n)$, so there exists a sequence
$q_k\in\mathcal B(n)$ converging to $f$ locally uniformly in $\Delta$.
Then $\{q_k\}_n\to\{f\}_n$, and for sufficiently large $k$ the inequality $\{q_k\}_n>2/e$ holds.
\end{proof}

\smallskip
Thus, for $n\in\mathbb N\setminus\{1\}$, the existence in the class $B$ of a function $f$ with $0<\{f\}_0<e^{-1}$ and
$\{f\}_n=2e^{-1}$ is equivalent to falsity of the Krzyz conjecture.
The Caratheodory-Toeplitz criterion easily shows that $F(-z^n,1)$ is the unique function $g\in B$ satisfying
$0<\{g\}_0<1$, $\{g\}_1=\ldots=\{g\}_{n-1}=0$, and $\{g\}_n=2e^{-1}$.
Since $f\neq F(-z^n,1)$, we have $|\{f\}_1|+\ldots+|\{f\}_{n-1}|>0$.

For $n\in\mathbb N\setminus\{1\}$, the Krzyz conjecture is also false if the class $B_1$ contains a function $f$ such
that $\{f\}_n=2\{f\}_0=2e^{-1}$ and $|\{f\}_1|+\ldots+|\{f\}_{n-1}|>0$.
Theorem~\ref{th:Stupin:uniqunessOfExtremalFunctionInB} implies that such a function $f$ does not belong to $\mathbb
E(n)$, that is, it cannot be extremal.
We state the following somewhat more general assertion:

\begin{Corollary}
If $n\in\mathbb{N}$, $f\in\mathbb E(n)$, and $f(z)\neq F(-z^n,1)$, where the function $F$ is given by
formula~\eqref{eq:Stupin:BExtrF}, then $\{f\}_n>2\{f\}_0>0$.
\end{Corollary}

\begin{proof}
By Corollary~\ref{cor:Stupin:f_n_eq_f_0_norm_Psi2}, we have $\{f\}_n\geqslant2\{f\}_0$.
Equality is impossible by Theorem~\ref{th:Stupin:uniqunessOfExtremalFunctionInB}, since $f(z)\neq F(-z^n,1)$.
Hence $\{f\}_n>2\{f\}_0$.
The inequality $\{f\}_0>0$ follows from the definition of a function of extremal type.
\end{proof}

\section{Linear invariance}

The following result (see~\cite{Ermers:ref1}) follows from the fact that the class $B$ is a linearly invariant family of
functions.

\begin{Theorem}
Let $n\in\mathbb{N}$, $f\in B$, $\{f\}_n>0$, and let $f$ be a locally extremal function in the topology of locally
uniform convergence.
Then
    \begin{equation}
        \label{eq:Stupin:Marty}
        (n+1)\{f\}_{n+1}
        =
        (n-1)\overline{\{f\}_{n-1}}.
    \end{equation}
\end{Theorem}

\begin{proof}
Let $\omega(z)\coloneqq\frac{z+\eta}{1+\overline{\eta}z}$.
Since $f\in B$, we have
    $
    f_{\eta}(z)
    \coloneqq
    f(\omega(z))
    \in
    B
$ for $\eta\in\Delta$.
As $\eta\to0$, we have $f_\eta\to f$ locally uniformly in $\Delta$.
In a neighborhood of $\eta=0$ we have
    $
    (\omega(z))^m
    =
    z^m+mz^{m-1}(\eta-\overline{\eta}z^2)+o(|\eta|)
    $.
Hence,
    \[
        f_{\eta}(z)
        =
        f(\omega(z))
        =
        \{f\}_0+
        \sum\limits_{m=1}^{\infty}
        \{f\}_m(z^m+mz^{m-1}(\eta-\overline{\eta}z^2))+o(|\eta|),
        \quad
        \eta\to0,
    \]
and
    \begin{equation*}
        \{f_{\eta}\}_n
        =
        \{f\}_n+(n+1)\{f\}_{n+1}\eta-(n-1)\{f\}_{n-1}\overline{\eta}+o(|\eta|),
        \quad
        \eta\to0.
    \end{equation*}
Since $\operatorname{Re}u\overline v=\operatorname{Re}\overline uv$ for all $u,v\in\mathbb{C}$, we have
    \begin{equation*}
        \operatorname{Re}\{f_{\eta}\}_n
        =
        \{f\}_n+\operatorname{Re}\big((n+1)\{f\}_{n+1}\eta-(n-1)\overline{\{f\}_{n-1}}\eta\big)+o(|\eta|),
        \quad
        \eta\to0.
    \end{equation*}
Since $f$ is locally extremal in the topology of locally uniform convergence,
    \begin{equation}
        \label{eq:Stupin:Marty_proof}
        \operatorname{Re}{\Big(\,\big((n+1)\{f\}_{n+1}-(n-1)\overline{\{f\}_{n-1}}\big)\,\eta\Big)}
        +
        o(|\eta|)
        \leqslant
        0,
        \quad
        \eta\to0.
    \end{equation}
Inequality~\eqref{eq:Stupin:Marty_proof} holds for every sufficiently small $\eta\in\Delta$, so we can choose $\eta$ so
small in modulus that $o(|\eta|)$ no longer affects the sign of the left-hand side of
inequality~\eqref{eq:Stupin:Marty_proof}.
However, if
    \[
        (n+1)\{f\}_{n+1}\neq(n-1)\overline{\{f\}_{n-1}},
    \]
we can choose $\eta$ so that the left-hand side of inequality~\eqref{eq:Stupin:Marty_proof} becomes positive,
contradicting~\eqref{eq:Stupin:Marty_proof}, as required.
\end{proof}

\begin{Lemma}
    \label{lem:Stupin:boundaryZeroDerivative}
Let $H$ be a polynomial, $\operatorname{Re}H(e^{i\varphi})\geqslant0$ for $\varphi\in\mathbb R$, and
    $\operatorname{Re}H(e^{i\varphi_0})=0$.
If $z_0\coloneqq e^{i\varphi_0}$, then
    $
        z_0 H'(z_0)
        =
        \overline{z_0 H'(z_0)}
$, that is, $z_0 H'(z_0)\in\mathbb{R}$, or, equivalently, $\operatorname{Im}\big(z_0 H'(z_0)\big)=0$.
\end{Lemma}

\begin{proof}
Set $T(\varphi)\coloneqq\operatorname{Re}H(e^{i\varphi})$.
By assumption, the function $T$ attains its minimum at $\varphi_0$, so $T'(\varphi_0)=0$.
On the other hand,
    \[
        T'(\varphi)
        =
        \operatorname{Re}\big(ie^{i\varphi}H'(e^{i\varphi})\big)
        =
        -\operatorname{Im}\big(e^{i\varphi}H'(e^{i\varphi})\big).
    \]
Hence $\operatorname{Im}(z_0H'(z_0))=0$, which implies the required equality.
\end{proof}

\smallskip
The following theorem establishes equality~\eqref{eq:Stupin:pointwiseMarty} independently of the choice of topology:

\begin{Theorem}
    \label{th:Stupin:pointwiseMarty}
Let $n\in\mathbb N$, $f\in\mathbb E(n)$.
Then for every $j=1,\ldots,m$ the following equality holds:
    \begin{equation}
        \label{eq:Stupin:pointwiseMarty}
        \sum\limits_{k=1}^{n+1}k\{f\}_{n+1-k}z_j^k
        =
        \sum\limits_{k=1}^{n-1}k\overline{\{f\}_{n-1-k}}z_j^{-k}.
    \end{equation}
\end{Theorem}

\begin{proof}
By the definition of a function of extremal type, $H\in C^n$, and there exist a number $m\in\{1,\ldots,n\}$ and points
    $z_j\coloneqq e^{i\varphi_j}$, $j=1,\ldots,m$, $0\leqslant\varphi_1<\ldots<\varphi_m<2\pi$,
for which equalities~\eqref{eq:Stupin:extrPolynomial_phi} hold.

Fix $j\in\{1,\ldots,m\}$.
The function $T(\varphi)\coloneqq\operatorname{Re}H(e^{i\varphi})$ is nonnegative and $T(\varphi_j)=0$, so
$H(z_j)+\overline{H(z_j)}=0$ and $z_j H'(z_j)=\overline{z_j H'(z_j)}$ by Lemma~\ref{lem:Stupin:boundaryZeroDerivative}.
These equalities give
    \[
        \{f\}_n
        +
        \sum\limits_{k=1}^n\{f\}_{n-k}z_j^k
        +
        \sum\limits_{k=1}^n\overline{\{f\}_{n-k}}z_j^{-k}
        =
        0
    \]
and
    \[
        \sum\limits_{k=1}^n k\{f\}_{n-k}z_j^k
        -
        \sum\limits_{k=1}^n k\overline{\{f\}_{n-k}}z_j^{-k}
        =
        0.
    \]

Adding these equalities and multiplying the result by $z_j$ gives equality~\eqref{eq:Stupin:pointwiseMarty}.
\end{proof}

\begin{Corollary}
If $n\in\mathbb N$ and $f\in\mathbb E(n)$, then
    $
        (n+1)\{f\}_{n+1}
        =
        (n-1)\overline{\{f\}_{n-1}}$.
\end{Corollary}

\begin{proof}
By the definition of a function of extremal type and formula~\eqref{eq:Stupin:B_t_extremal_h}, fix a representation
    \[
        h(z)
        \coloneqq
        -\ln f(z)
        =
        \sum\limits_{j=1}^m\alpha_j\frac{1+z_j z}{1-z_j z},
    \]
where $m\in\{1,\ldots,n\}$, $\alpha_j>0$, $z_j\coloneqq e^{i\varphi_j}$, $j=1,\ldots,m$,
    $0\leqslant\varphi_1<\ldots<\varphi_m<2\pi$.
Then
    \[
        \{h\}_k
        =
        2\sum\limits_{j=1}^m\alpha_j z_j^k,
        \quad
        \overline{\{h\}_k}
        =
        2\sum\limits_{j=1}^m\alpha_j z_j^{-k},
        \quad
        k\geqslant1.
    \]
Multiplying equality~\eqref{eq:Stupin:pointwiseMarty} of Theorem~\ref{th:Stupin:pointwiseMarty} by $2\alpha_j$ and
summing over $j=1,\ldots,m$ gives
    \begin{equation}
        \label{eq:Stupin:semiMarty}
        \sum\limits_{k=1}^{n+1}k\{h\}_k\{f\}_{n+1-k}
        =
        \sum\limits_{k=1}^{n-1}k\overline{\{h\}_k}\,\overline{\{f\}_{n-1-k}}.
    \end{equation}

Formula~\eqref{eq:Stupin:deForBCoeff} for $\ell=n+1$ gives
    \[
        (n+1)\{f\}_{n+1}
        =
        -\sum\limits_{k=1}^{n+1}k\{h\}_k\{f\}_{n+1-k},
    \]
and for $n>1$ and $\ell=n-1$ gives
    \[
        (n-1)\overline{\{f\}_{n-1}}
        =
        -\sum\limits_{k=1}^{n-1}k\overline{\{h\}_k}\,\overline{\{f\}_{n-1-k}}.
    \]
For $n=1$ the last equality also holds, since both sides vanish.
Hence,
    $
        (n+1)\{f\}_{n+1}
        =
(n-1)\overline{\{f\}_{n-1}}$, by equality~\eqref{eq:Stupin:semiMarty}.
\end{proof}

By Theorem~\ref{th:Stupin:localExtremalIsExtremalType}, this corollary also holds for every function $f\in\mathfrak
E(n)$.
Thus, for functions in $\mathbb E(n)$, equality~\eqref{eq:Stupin:Marty} follows directly from the properties of the
function $h$ and the polynomial $H$.

\smallskip
R.~Ermers notes in~\cite{Ermers:ref1} that the function $F(-z^n,1)$ satisfies relation~\eqref{eq:Stupin:Marty}.
It is also clear that for $n>1$ all $F(-z^n,t)$, $t>0$, trivially satisfy relation~\eqref{eq:Stupin:Marty}, since all
their coefficients whose indices are not multiples of $n$ vanish ($F(-z,t)$ satisfy equality~\eqref{eq:Stupin:Marty}
only for $t=1$).
Other functions also satisfy this relation.
For example, if $n$ is even, every even function in the class $B$ satisfies relation~\eqref{eq:Stupin:Marty}
(see~\cite{Ermers:ref1}).
Moreover, if $f\in B$ satisfies relation~\eqref{eq:Stupin:Marty}, then $cf\in B$ also satisfies this relation for
$c\in(0,1]$.

\begin{Corollary}
Let $n\in\mathbb N\setminus\{1\}$, $f\in\mathbb E(n)$.
Then the polynomial
    \[
        R_n(z)
        \coloneqq
        \sum\limits_{k=1}^{n+1}k\{f\}_{n+1-k}z^{n-1+k}
        -
        \sum\limits_{k=1}^{n-1}k\overline{\{f\}_{n-1-k}}z^{n-1-k}
    \]
has degree $2n$, and the polynomial $S_m(z)\coloneqq\prod\limits_{j=1}^m(z-z_j)$ divides $R_n$.
\end{Corollary}

\begin{proof}
The leading coefficient $r_{2n}=(n+1)\{f\}_0\neq0$, so the degree of the polynomial $R_n$ is $2n$.
Multiplying equality~\eqref{eq:Stupin:pointwiseMarty} of Theorem~\ref{th:Stupin:pointwiseMarty} by $z_j^{n-1}$ gives
$R_n(z_j)=0$ for $j=1,\ldots,m$.
This proves the assertion.
\end{proof}

\section{Inequalities between coefficients}

Let $n\in\mathbb{N}$, $f\in\mathbb E(n)$, and let the polynomial $H$ be generated by the coefficients of the function
$f$ according to formula~\eqref{eq:Stupin:extrPolynomialH_n}.
The definition of a function of extremal type implies that $\frac{H}{\{f\}_n}\in C_1^n$, with
\[
    \frac{H(z)}{\{f\}_n}
    =
    1+2\sum\limits_{k=1}^n\frac{\{f\}_{n-k}}{\{f\}_n}z^k.
\]
Applying Theorem~\ref{th:Stupin:estimOnC_n} to the polynomial $\frac{H}{\{f\}_n}$ immediately gives the following
result:

\begin{Theorem}
    \label{th:Stupin:cIneq1}
Let $n\in\mathbb{N}$.
If $f\in\mathbb E(n)$, then
    \[
        \frac{|\{f\}_k|}{\{f\}_n}
        \leqslant
        \cos\frac{\pi}{\big\lfloor\frac{n}{n-k}\big\rfloor + 2}
        <
        1,
        \quad
        k=0,\ldots,n-1.
    \]
\end{Theorem}

By Theorem~\ref{th:Stupin:localExtremalIsExtremalType}, Theorem~\ref{th:Stupin:cIneq1} also holds for every function
$f\in\mathfrak E(n)$.
Note that for $k=0$, Theorem~\ref{th:Stupin:cIneq1} gives $\tfrac{|\{f\}_0|}{\{f\}_n}\leqslant\cos\frac{\pi}3=\tfrac12$,
that is, $\{f\}_n\geqslant2\{f\}_0$.

\smallskip
Let $n\in\mathbb{N}$, and let the function $f$ be extremal.
Consider its internal variation
\[
    f_{\varepsilon,\zeta}(z)
    \coloneqq
    f(z)e^{-\varepsilon\frac{1+\zeta z^{n-k}}{1-\zeta z^{n-k}}},
    \quad
    \varepsilon>0,
    \quad
    \zeta\in\overline\Delta,
    \quad
    0\leqslant k<n/2.
\]
In the expansion of $\tfrac{1+\zeta z^{n-k}}{1-\zeta z^{n-k}}$ up to degree $z^n$, besides the constant term, only the
term $2\zeta z^{n-k}$ occurs if $2(n-k)>n$, that is, for $k<\tfrac{n}2$.
We have
$
\{f_{\varepsilon,\zeta}\}_n
=
\{f\}_n
-
(\{f\}_n+2\{f\}_k\zeta)\varepsilon
+
o(\varepsilon).
$
Since the function $f$ is extremal, there exist numbers $\delta>0$ and $0<\eta\leqslant1$ such that $|\{\widetilde
f\}_n|\leqslant\{f\}_n$ for all $\widetilde f\in V_{\delta,\eta}(f)$.
As in the proof of Theorem~\ref{th:Stupin:scalar_prod}, one can establish that $f_{\varepsilon,\zeta}\in
V_{\delta,\eta}(f)$ for all sufficiently small $\varepsilon>0$ and $\zeta\in\overline\Delta$.
Hence $\operatorname{Re}\{f_{\varepsilon,\zeta}\}_n\leqslant\{f\}_n$ for all sufficiently small $\varepsilon>0$ and
$\zeta\in\overline\Delta$.
Subtracting $\{f\}_n$ from both sides, dividing by $\varepsilon>0$, and passing to the limit as $\varepsilon\to0+$ gives
the inequality $\operatorname{Re}(\{f\}_n+2\{f\}_k\zeta)\geqslant0$ for $\zeta\in\overline\Delta$.
That is, the function $h(\zeta)\coloneqq\{f\}_n+2\{f\}_k\zeta\in C_{\{f\}_n}$.

\smallskip
Thus, for $f\in\mathfrak E(n)$ we obtain the estimations $|\{f\}_k|\leqslant\tfrac12\{f\}_n$, $0\leqslant k<\tfrac{n}2$.
Generalizing these estimations to functions of extremal type gives the following special case of
Theorem~\ref{th:Stupin:cIneq1}:

\begin{Theorem}
    \label{th:Stupin:cIneq2}
Let $n\in\mathbb{N}$.
If the function $f\in\mathbb E(n)$, then $|\{f\}_k|\leqslant\tfrac12\{f\}_n$, $0\leqslant k<\tfrac{n}2$.
\end{Theorem}

\begin{proof}
If $0\leqslant k<\tfrac{n}2$, then $1\leqslant\tfrac{n}{n-k}<2$, so $\left\lfloor\frac{n}{n-k}\right\rfloor=1$.
Theorem~\ref{th:Stupin:cIneq1} gives $\frac{|\{f\}_k|}{\{f\}_n}\leqslant\cos\frac{\pi}3=\frac12$, as required.
\end{proof}

\smallskip
Theorem~\ref{th:Stupin:cIneq2} for $k=0$ again gives $\{f\}_n\geqslant2\{f\}_0$.

\smallskip
Let $n\in\mathbb{N}$, $f\in\mathbb E(n)$.
Note that if $g$ satisfies inequality~\eqref{eq:Stupin:scalar_prod}, then for every $t>0$ the function $tg$ also
satisfies this inequality.
If $h$ satisfies equality~\eqref{eq:Stupin:orthogonalVects}, then for every $t>0$ the function $th$ also satisfies this
equality.
If $0<c<1$, then $cf\in B$.
Although $cf$ is no longer extremal (as explained in the next paragraph), it still satisfies
relations~\eqref{eq:Stupin:scalar_prod},~\eqref{eq:Stupin:f_n_ge_2f_0},
~\eqref{eq:Stupin:extrPolynomial},~\eqref{eq:Stupin:extrPolynomial_phi},~\eqref{eq:Stupin:roots},%
~\eqref{eq:Stupin:rootsOfH_abs}, \eqref{eq:Stupin:rootsOfH_prod_arg}
and relation~\eqref{eq:Stupin:Marty}, as well as the estimations of Theorems~\ref{th:Stupin:cIneq1}
and~\ref{th:Stupin:cIneq2}.

Note that equality~\eqref{eq:Stupin:orthogonalVects} is not invariant under multiplication by a constant distinct from
$1$.
Indeed, set $\mathfrak f\coloneqq cf$ and $\mathfrak h\coloneqq-\ln\mathfrak f=h-\ln c$, where $h\coloneqq-\ln f$.
Then, using equality~\eqref{eq:Stupin:orthogonalVects} for the function $f$, we obtain
\begin{multline*}
    \operatorname{Re}{(\{\mathfrak f\}_n\{\mathfrak h\}_0+\ldots+\{\mathfrak f\}_0\{\mathfrak h\}_n)}
    =\\=
    c\operatorname{Re}{(\{f\}_n\{h\}_0+\ldots+\{f\}_0\{h\}_n)}
    -c\{f\}_n\ln c
    =
    -c\{f\}_n\ln c
    >
    0.
\end{multline*}
Hence $\mathfrak f\notin\mathbb E(n)$ by Theorem~\ref{th:Stupin:scalar_prod_criteria}, and it is not extremal by
Theorem~\ref{th:Stupin:localExtremalIsExtremalType}.

\smallskip
Note that similar results for extremal functions, with a weaker estimation in Theorem~\ref{th:Stupin:cIneq1}, first
appeared implicitly in~\cite{Hummel:ref1}.
In~\cite{Ermers:ref1}, these results for extremal functions are stated explicitly, but again with a weaker first
estimation.
Here these results are generalized to functions of extremal type.

Note that all necessary conditions for local extremality discussed in this paper have been generalized to functions of
extremal type.
Every local extremal is of extremal type by Theorem~\ref{th:Stupin:localExtremalIsExtremalType}.
The converse is not proved in this paper.

\section{The set $\mathbb E(1)$}
\label{sec:Stupin:extremalTypeE1}

The class $C^1$ has a very simple structure, allowing a complete description of the set $\mathbb E(1)$.

\begin{Theorem}
    \label{th:Stupin:extremalTypeE1}
The set $\mathbb E(1)$ consists of a single function $F(-z,1)$, that is, $\mathbb E(1)=\{F(-z,1)\}$, where $F$ is
defined by formula~\eqref{eq:Stupin:BExtrF}.
\end{Theorem}

\begin{proof}
Let $f\in\mathbb E(1)$.
Then by the definition of the set $\mathbb E(1)$ there exist a point $z_1\in\partial\Delta$ and a number $t=\alpha_1>0$
such that
    \[
        f(z)
        =
        e^{-t\frac{1+z_1 z}{1-z_1 z}}
        =
        F(z_1 z,t)
        =
        e^{-t}-2te^{-t}z_1 z+o(z),
    \]
with $\{f\}_1>0$, $H\in C^1$, $\operatorname{Re}H(z_1)=0$.

Since $\{f\}_1>0$, we have $z_1=-1$.
Therefore,
    \[
        H(z)
        =
        \{f\}_1+2\{f\}_0 z
        =
        -2te^{-t}z_1+2e^{-t} z
        =
        2e^{-t}(t+z)
        \in C^1,
        \quad
        \operatorname{Re}H(-1)=0.
    \]
Hence $t=1$ and $f(z)=F(-z,1)$.
Conversely, $F(-z,1)\in\mathbb E(1)$ by Theorem~\ref{th:Stupin:uniqunessOfExtremalFunctionInB}.
\end{proof}

\smallskip
Note that for $n=1$ only rotations in the variable plane $z$ through multiples of $2\pi$ are admissible.

The set $\mathfrak E(1)$ is nonempty and, by Theorem~\ref{th:Stupin:localExtremalIsExtremalType}, is contained in
$\mathbb E(1)$.
Thus Theorem~\ref{th:Stupin:extremalTypeE1} gives another proof of the first assertion of
Theorem~\ref{th:Stupin:Krzyz1}.
Since every normalized global extremal belongs to $\mathfrak E(1)$, the maximum of $|\{f\}_1|$ on the class $B$ equals
$2/e$.
This gives another proof of the second assertion of Theorem~\ref{th:Stupin:Krzyz1}, and hence of the entire
Theorem~\ref{th:Stupin:Krzyz1}.
Moreover, $\mathfrak E(1)=\mathbb E(1)$ by Theorems~\ref{th:Stupin:Krzyz1} and~\ref{th:Stupin:extremalTypeE1}.

Note that the function $f(z)\coloneqq F(-z,1)=e^{-1}+2e^{-1}z+o(z^2)$ satisfies relation~\eqref{eq:Stupin:Marty}, that
is, $0\overline{\{f\}_0}=2\{f\}_2$, since $\{f\}_2=0$.
Conversely, for $n=1$ formula~\eqref{eq:Stupin:Marty} gives $\{f\}_2=0$ for $f\in\mathbb E(1)$.

\section{The set $\mathbb E(2)$}

Let $n\in\mathbb N$ and $f\in\mathbb E(n)$.
Then by the definition of $\mathbb E(n)$ there exists a number $m\in\{1,\ldots,n\}$ such that
$h(z)\coloneqq-\ln f(z)=\sum\limits_{j=1}^m\alpha_j\frac{1+z_j z}{1-z_j z}$,
$\alpha_j>0$, $j=1,\ldots,m$, and the points $z_1,\ldots,z_m$ are pairwise distinct.
For $m=1,\ldots,n$, let $\mathbb E_m(n)$ denote the set of functions $f\in\mathbb E(n)$ for which $h$ contains exactly
$m$ terms.
Then
\[
    \mathbb E(n)
    =
    \bigcup\limits_{m=1}^n\mathbb E_m(n),
    \quad
    \mathbb E_m(n)\cap\mathbb E_\ell(n)
    =
    \emptyset,
    \quad
    m,\ell\in\{1,\ldots,n\},
    \quad
    m\neq\ell.
\]
Indeed, every point $\overline{z_j}$ is a pole of the function $h$, and all these points are pairwise distinct.
Hence the number of terms in the representation of $h$ equals the number of poles of $h$ and is uniquely determined.

\begin{Lemma}
    \label{lem:Stupin:extremalTypeE2OneTerm}
The set $\mathbb E_1(2)$ has the following form:
    $\mathbb E_1(2)
        =
        \left\{
            F(\pm z,\tfrac{3+\sqrt5}2)
        \right\}$,
where $F$ is defined by formula~\eqref{eq:Stupin:BExtrF}.
\end{Lemma}

\begin{proof}
Let $f\in\mathbb E_1(2)$.
Then by the definition of $\mathbb E_1(2)$ there exist a point $z_1\in\partial\Delta$ and a number $t=\alpha_1>0$ such
that
    \[
        f(z)
        =
        e^{-t\frac{1+z_1 z}{1-z_1 z}}
        =
        F(z_1 z,t)
        =
        e^{-t} - 2te^{-t}z_1 z + 2t(t-1)e^{-t}z_1^2 z^2 + o(z^2),
    \]
with $\{f\}_2>0$, $H\in C^2$, $\operatorname{Re}H(z_1)=0$.

Since $\{f\}_2>0$, we have $t\neq1$, with $z_1^2=-1$ for $0<t<1$ and $z_1^2=1$ for $t>1$.
Furthermore,
    \[
        0
        =
        \operatorname{Re}H(z_1)
        =
        \operatorname{Re}\bigl(
            2t(t-1)e^{-t}z_1^2
            +2(-2t e^{-t}z_1)z_1
            +2(e^{-t})z_1^2
        \bigr)
        =
        2e^{-t}z_1^2(t^2-3t+1).
    \]
The resulting equation has two roots $t=t_1\coloneqq\tfrac{3-\sqrt5}2\approx0.382$ and
$t=t_2\coloneqq\tfrac{3+\sqrt5}2\approx2.618$.

The value $t=t_1<1$ is inadmissible, since in this case $z_1^2=-1$ and
    \[
        \operatorname{Re}H(-z_1)
        =
        \operatorname{Re}\bigl(
            2t(t-1)e^{-t}z_1^2
            +2(-2t e^{-t}z_1)(-z_1)
            +2(e^{-t})z_1^2
        \bigr)
        =
        -2e^{-t}(t^2+t+1)<0.
    \]
By continuity, $\operatorname{Re}H$ also takes negative values at some points of the disk $\Delta$, contradicting the
condition $H\in C^2$.

The case $t=t_2>1$ corresponds to $z_1=\pm1$, whence $f(z)=F(\pm z,t_2)$.
It remains to check that the functions found satisfy the condition $H\in C^2$.
Let $s\in\{-1,1\}$ and $f(z)=F(sz,t_2)$.
Then
    $
        H(z)
        =
        2e^{-t_2}\bigl(t_2(t_2-1)-2t_2sz+z^2\bigr)$.
For $z=e^{i\theta}$, where $\theta\in[0,2\pi)$, set $x\coloneqq s\cos\theta$.
Using the equality $t_2^2-3t_2+1=0$, we obtain
    \[
        \operatorname{Re}H(e^{i\theta})
        =
        2e^{-t_2}\bigl(t_2(t_2-1)-2t_2x+2x^2-1\bigr)
        =
        4e^{-t_2}(x-1)(x-t_2+1)
        \geqslant
        0.
    \]
The last inequality follows from $x\in[-1,1]$ and $t_2>2$.
Since $H(0)>0$, the minimum principle for harmonic functions gives $\operatorname{Re}H(z)>0$ for $z\in\Delta$.
Hence $H\in C^2$, so both functions found belong to $\mathbb E_1(2)$.
\end{proof}

\smallskip
Thus we have found two functions of extremal type
\begin{multline*}
    F(\pm z,t_2)
    =
    e^{-t_2}\mp2t_2e^{-t_2}z+2t_2(t_2-1)e^{-t_2}z^2
    \\
    \mp\tfrac23t_2e^{-t_2}z^3+o(z^3)
    \approx
    0.073\mp0.382z+0.618z^2\mp0.127z^3+\ldots.
\end{multline*}
For each of these functions $3\{f\}_3=\mp2t_2e^{-t_2}=\overline{\{f\}_1}$, so they satisfy
relation~\eqref{eq:Stupin:Marty}.

\begin{Lemma}
    \label{lem:Stupin:extremalTypeE2TwoTerms}
The set $\mathbb E_2(2)$ has the following form:
    \[
        \mathbb E_2(2)
        =
        \left\{F(-z^2,1),f_{-\frac1{\sqrt2}},f_{\frac1{\sqrt2}}\right\},
        \quad
        f_{\mp\frac1{\sqrt2}}(z)
        \coloneqq
        \exp\left(-2\frac{1-z^2}{1\pm\sqrt2z+z^2}\right),
    \]
where $F$ is defined by formula~\eqref{eq:Stupin:BExtrF}.
\end{Lemma}

\begin{proof}
Let $f\in\mathbb E_2(2)$.
Then by the definition of $\mathbb E_2(2)$ there exist points $z_1,z_2\in\partial\Delta$, $z_1\neq z_2$ and numbers
$\alpha_1,\alpha_2>0$, $t=\alpha_1+\alpha_2$ such that
    \begin{equation}
        \label{eq:Stupin:2f_1E_2_2}
        \begin{aligned}
        f(z)
        =&
        e^{-\alpha_1\frac{1+z_1 z}{1-z_1 z}-\alpha_2\frac{1+z_2 z}{1-z_2 z}}
        =\\=&
        e^{-t}
        -
        2e^{-t}(\alpha_1 z_1+\alpha_2 z_2) z
        +
        2e^{-t}((\alpha_1 z_1+\alpha_2 z_2)^2-(\alpha_1 z_1^2+\alpha_2 z_2^2)) z^2
        +
        o(z^2),
        \end{aligned}
    \end{equation}
with $\{f\}_2>0$, $H\in C^2$, $\operatorname{Re}H(z_1)=\operatorname{Re}H(z_2)=0$.

Since $H\in C^2$, the Fejer-Riesz theorem~\ref{th:Stupin:FR_complex_case} gives a polynomial $P$ external for $H$.
By Corollary~\ref{cor:Stupin:extremalTypeZerosOfH}, it has the form
    \[
        P(z)=p_0+p_1 z+p_2 z^2=p_2(z-z_1)(z-z_2),
        \quad
        p_2>0.
    \]
Hence $p_0=p_2 z_1 z_2$.
Formula~\eqref{eq:Stupin:f_n_p_n} for $k=2$ gives
    $
        \{f\}_0
        =
        p_2\overline{p_0}
        =
        p_2^2\overline{z_1 z_2}$.
Since $p_2>0$ and $\{f\}_0>0$, we have $z_1 z_2=1$ and $p_2^2=\{f\}_0$.
Therefore,
    \[
        z_2=\overline{z_1},
        \quad
        P(z)=p_2(1-2xz+z^2),
        \quad
        x\coloneqq\operatorname{Re}z_1=\operatorname{Re}z_2,
        \quad
        p_0=p_2,
        \quad
        p_1=-2xp_2.
    \]

Formula~\eqref{eq:Stupin:f_n_p_n} for $k=1$ gives
    \begin{equation}
        \label{eq:Stupin:1f_1E_2_2}
        \{f\}_1
        =
        p_1\overline{p_0}+p_2\overline{p_1}
        =
        -4x p_2^2
        =
        -4x\{f\}_0
    \end{equation}
Formulas~\eqref{eq:Stupin:1f_1E_2_2} and~\eqref{eq:Stupin:2f_1E_2_2} imply that $2x=\alpha_1 z_1+\alpha_2
\overline{z_1}$.
Since $z_1\neq z_2$ and $z_2=\overline{z_1}$, we have $\operatorname{Im}z_1\neq0$.
Hence comparing imaginary parts in the last formula gives $\alpha_1=\alpha_2=t/2$, and comparing real parts gives
$tx=2x$.
Thus $x=0$ or $t=2$.

Formula~\eqref{eq:Stupin:f_n_p_n} for $k=0$ gives
    \[
        \{f\}_2
        =
        |p_0|^2+|p_1|^2+|p_2|^2
        =
        (2+4x^2)p_2^2
        =
        (2+4x^2)\{f\}_0.
    \]
Since $z_1^2+\overline{z_1}^2=(z_1+\overline{z_1})^2-2|z_1|^2=4x^2-2$ (because $|z_1|=1$),
formula~\eqref{eq:Stupin:2f_1E_2_2} gives
    \[
        \{f\}_2=2t\{f\}_0\bigl(1+(t-2)x^2\bigr).
    \]
Comparing this equality with the preceding formula for $\{f\}_2$ gives
    \[
        t-1+(t^2-2t-2)x^2=0
    \]
whence, since we know that $x=0$ or $t=2$, we obtain
    \[
        x=0,\quad t=1,
\quad\text{or}\quad
        x^2=\frac12,\quad t=2.
    \]
In the first case $z_1=\pm i$, $z_2=\mp i$, so $h(z)=\frac{1-z^2}{1+z^2}$ and $f=F(-z^2,1)$.
In the second case $h(z)=2\frac{1-z^2}{1-2xz+z^2}$, so $f=f_x$, where $x=\pm\tfrac1{\sqrt2}$.

It remains to check that the functions found satisfy the condition $H\in C^2$.
Recall that $x\coloneqq\operatorname{Re}z_1=\operatorname{Re}z_2$.
For $(t,x)=(1,0)$ or $(t,x)=(2,\pm\tfrac1{\sqrt2})$, we have
    \[
        H(z)=2e^{-t}(1+2x^2-4xz+z^2),
        \quad
        \operatorname{Re}H(e^{i\theta})=4e^{-t}(\cos\theta-x)^2\geqslant0,
        \quad
        \theta\in[0,2\pi).
    \]
Since $H(0)>0$, the minimum principle for harmonic functions gives $\operatorname{Re}H(z)>0$ for $z\in\Delta$.
Hence $H\in C^2$.

Thus all three functions found belong to $\mathbb E_2(2)$.
\end{proof}

\smallskip
Thus we have found another $3$ functions of extremal type
\[
    \begin{aligned}
    F(-z^2,1)
    =&
    e^{-1}+2e^{-1}z^2+o(z^3)
    \approx
    0.368+0.736z^2+o(z^3),
    \\
    f_{\mp\frac1{\sqrt2}}(z)
    =&
    e^{-2}\pm2\sqrt2e^{-2}z+4e^{-2}z^2
    \pm\tfrac{2\sqrt2}3e^{-2}z^3+o(z^3)
    \approx \\ &\approx
    0.135\pm0.383z+0.541z^2\pm0.128z^3+o(z^3).
    \end{aligned}
\]
For the function $F(-z^2,1)$ we have $3\{f\}_3=0=\overline{\{f\}_1}$, and for the other two functions
\[
    3\{f_{\mp\frac1{\sqrt2}}\}_3
    =
    \pm2\sqrt2e^{-2}
    =
    \overline{\{f_{\mp\frac1{\sqrt2}}\}_1}.
\]
Hence all three functions satisfy relation~\eqref{eq:Stupin:Marty}.

\begin{Theorem}
    \label{th:Stupin:extremalTypeE2}
The set $\mathbb E(2)$ consists of the following five functions:
    \[
        \mathbb E(2)
        =
        \left\{
            F\left(\pm z,\tfrac{3+\sqrt5}2\right),
            F(-z^2,1),
            f_{\pm\frac1{\sqrt2}}
        \right\},
    \]
where $F$ is defined by formula~\eqref{eq:Stupin:BExtrF}.
\end{Theorem}

\begin{proof}
The assertion follows immediately from Lemmas~\ref{lem:Stupin:extremalTypeE2OneTerm}
and~\ref{lem:Stupin:extremalTypeE2TwoTerms}.
\end{proof}

\smallskip
The following assertion is a consequence of Theorem~\ref{th:Stupin:extremalTypeE2}:

\begin{Theorem}
The Krzyz conjecture holds for $n=2$.
\end{Theorem}

\begin{proof}
Let $f$ be a normalized global extremal in the Krzyz problem for the second coefficient.
Then $f$ is a local extremal and, by Theorem~\ref{th:Stupin:localExtremalIsExtremalType}, belongs to $\mathbb E(2)$.
By Theorem~\ref{th:Stupin:extremalTypeE2},
    \[
        \{f\}_2
        =
        \max\limits_{g\in\mathbb E(2)} \{g\}_2
        =
        \max
        \left\{
            \big(3+\sqrt5\big)\big(\tfrac{3+\sqrt5}2-1\big)e^{-\frac{3+\sqrt5}2},~
        \frac2e,~
        \frac4{e^2}
        \right\}
        =
        \frac2e.
    \]
Hence $f(z)=F(-z^2,1)$ and $\{f\}_2=2/e$.
If a function $g\in B$ is a global extremal, normalizing it by rotations in the variable planes $z$ and $w$ gives a
function in $\mathbb E(2)$, that is, $F(-z^2,1)$.
Thus equality is attained only by functions of the form $e^{i\theta}F(e^{i\varphi}z^2,1)$, where
$\theta,\varphi\in[0,2\pi)$.
\end{proof}

\smallskip
\smallskip
Note that for $n=2$ only rotations in the variable plane $z$ through multiples of $\pi$ are admissible.
Thus, up to replacing $z$ by $-z$, the set $\mathbb E(2)$ contains three functions.
Among all functions in $\mathbb E(2)$, only $f(z)=F(-z^2,1)$ satisfies the equality $\{f\}_2=2\{f\}_0$ (see also
Theorem~\ref{th:Stupin:uniqunessOfExtremalFunctionInB}), and only this function is a global extremal.

\section{The sets $\mathbb E(2)$ and $\mathfrak E(2)$}

It can be shown that among all functions $f$ in $\mathbb E(2)$, a local maximum of $|\{f\}_2|$ in the topology of
locally uniform convergence induced on $B_+$ is attained only by the functions
$F\left(\pm z,\tfrac{3+\sqrt5}2\right)$, $F(-z^2,1)$.
Up to rotations in the variable plane $z$, these are the functions $F\left(z,\tfrac{3+\sqrt5}2\right)$ and $F(-z^2,1)$.
It can also be shown that in the topology $\tau_{B_+}$ all functions in $\mathbb E(2)$ are locally extremal, and hence
$\mathfrak E(2)=\mathbb E(2)$.

\smallskip
We supplement Assertion~\ref{st:Stupin:segmentTopologyImpliesLocalUniform}.
To this end, note that extremality of the functions $f_{\mp\frac1{\sqrt2}}$ in the topology $\tau_{B_+}$, but not in the
topology of locally uniform convergence in $\Delta$, shows that the topology $\tau_{B_+}$ is strictly stronger than the
topology of locally uniform convergence in $\Delta$.

\smallskip
Since $\mathfrak E(1)=\mathbb E(1)$ and $\mathfrak E(2)=\mathbb E(2)$, it is natural to conjecture that $\mathfrak
E(n)=\mathbb E(n)$ for every $n\in\mathbb{N}$.

Note that for $G(z)=F(z,t)$ inequality~\eqref{eq:Stupin:f_2} becomes
\[
    |\{f\}_2|
    \leqslant
    2te^{-t}\max(1,|t-1|),
    \quad
    f\in B_t.
\]
The equality cases described after it show that every function in $\mathbb E(2)$ attains a global maximum of $|\{f\}_2|$
in the corresponding class $B_t$.
In particular, the functions $f_{\pm\frac1{\sqrt2}}$ attain a global maximum in the class $B_2$, but do not attain a
local maximum in the topology of locally uniform convergence induced on $B_+$.

\section{The set $\mathbb E_1(3)$}

\begin{Lemma}
    \label{lem:Stupin:extremalTypeE3OneTerm}
The set $\mathbb E_1(3)$ has the following form\textnormal:
    \[
        \mathbb E_1(3)
        =
        \left\{
            F(e^{i(2k+1)\frac{\pi}3}z,t_1^*)
            :
            k=0,1,2
        \right\},
    \]
where the function $F$ is given by formula~\eqref{eq:Stupin:BExtrF},
    $t_1^*\coloneqq2+\sqrt6\cos\vartheta$, $\vartheta\coloneqq\tfrac13\arccos\tfrac{5\sqrt6}{18}$.
In particular, $|\mathbb E_1(3)|=3$.
\end{Lemma}

\begin{proof}
Let $f\in\mathbb E_1(3)$.
Then there exist $t>0$ and $z_1\in\partial\Delta$ such that
    \[
        f(z)
        =
        e^{-t\frac{1+z_1 z}{1-z_1 z}}
        =
        F(z_1 z,t)
        =
        e^{-t}
        -
        2te^{-t}z_1 z
        +
        2t(t-1)e^{-t}z_1^2 z^2
        +
        \left(-2t+4t^2-\tfrac43t^3\right)e^{-t}z_1^3 z^3
        +
        o(z^3).
    \]
Since $\{f\}_3>0$ and $t>0$, we have $z_1^3=\pm1$.

We write out
        $H(z)
        =
        \left(-2t+4t^2-\tfrac43t^3\right)e^{-t}z_1^3
        +
        2(2t(t-1)e^{-t}z_1^2) z
        +
        2(-2te^{-t}z_1) z^2
        +
        2e^{-t} z^3$.
We have
        $H(z_1)
        =
        -\tfrac23e^{-t}z_1^3(2t^3-12t^2+15t-3)$.
Thus the condition $\operatorname{Re}H(z_1)=0$ is equivalent to the equation $2t^3-12t^2+15t-3=0$, which has three
positive roots
    \[
        t_1
        \coloneqq
        2+\sqrt6\cos\left(\vartheta+\tfrac{2\pi}3\right)\approx0.247,
        ~
        t_2
        \coloneqq
        2+\sqrt6\cos\left(\vartheta+\tfrac{4\pi}3\right)\approx1.395,
        ~
        t_3
        \coloneqq
        2+\sqrt6\cos\vartheta\approx4.358.
    \]

The real factor $\left(-2t+4t^2-\tfrac43t^3\right)e^{-t}$ is negative for $t=t_1,t_3$ and positive for $t=t_2$.
Thus the condition $\{f\}_3>0$ gives
    \[
z_1^3=-1~\text{for}~t=t_1~\text{and for}~t=t_3,
        \quad
z_1^3=1~\text{for}~t=t_2.
    \]

We check the condition $H\in C^3$ for each $t\in\{t_1,t_2,t_3\}$.
Write $z_1=e^{i\varphi_1}$, $z=e^{i\varphi}$, where $\varphi_1,\varphi\in[0,2\pi)$, and set
$x\coloneqq\cos(\varphi-\varphi_1)$.
Let
    \[
        G(x,t)
        \coloneqq
        e^{-t}(1-x)
        \left(
            8\big(x-\tfrac{t-1}2\big)^2+2t(t-4)
        \right).
    \]

Let $t=t_1$.
Then $z_1^3=-1$.
Substituting these values into the expression for $\operatorname{Re}H(e^{i\varphi})$ and completing the square in the
variable $x$ gives
    \[
        \operatorname{Re}H(e^{i\varphi})
        =
        G(x,t_1).
    \]
For $x=\tfrac{t_1-1}2\in[-1,1]$ the right-hand side is negative.
Thus by continuity $H\notin C^3$.

Let $t=t_2$.
Then $z_1^3=1$, and similarly
    \[
        \operatorname{Re}H(e^{i\varphi})
        =
        -G(x,t_2).
    \]
For $x=-1$ the expression in parentheses equals $(2t_2-1)^2+1>0$.
Thus the right-hand side is negative, and by continuity $H\notin C^3$.

Let $t=t_3$.
Then $z_1^3=-1$ and
    \[
        \operatorname{Re}H(e^{i\varphi})
        =
        G(x,t_3)
        \geqslant
        0,
    \]
since $x\in[-1,1]$ and $t_3>4$.
The function $\operatorname{Re}H$ is harmonic in $\Delta$, and $H(0)=\{f\}_3>0$.
By the minimum principle, $\operatorname{Re}H(z)>0$ for all $z\in\Delta$, that is, $H\in C^3$.

Thus $t=t_3=t_1^*$ and $z_1=e^{i(2k+1)\frac{\pi}3}$, $k=0,1,2$.
The calculations show that every function obtained satisfies all conditions in the definition of $\mathbb E_1(3)$.
The functions obtained are pairwise distinct, since their first coefficients are pairwise distinct.
\end{proof}

\smallskip
\smallskip
Note that for $n=3$ only rotations in the variable plane $z$ through multiples of $2\pi/3$ are admissible.
Thus, up to rotations in the variable plane $z$, we have found the function of extremal type
\label{page:Stupin:E3OneTermExpansion}
\begin{multline*}
    F(e^{i\frac\pi3}z,t_1^*)
    =
    e^{-t_1^*}-2t_1^*e^{-t_1^*}e^{i\frac\pi3}z
    +2t_1^*(t_1^*-1)e^{-t_1^*}e^{i\frac{2\pi}3}z^2
    +\\+
    \tfrac{2t_1^*}3\bigl(2(t_1^*)^2-6t_1^*+3\bigr)e^{-t_1^*}z^3
    +t_1^*(t_1^*-1)e^{-t_1^*}e^{i\frac{4\pi}3}z^4
    +o(z^4)
    \approx \\ \approx
    0.013-(0.056+0.097i)z-(0.187-0.325i)z^2+0.552z^3-(0.094+0.162i)z^4+o(z^4).
\end{multline*}
This expansion and the equality $e^{-2i(2k+1)\frac\pi3}=e^{4i(2k+1)\frac\pi3}$ imply that
\[
    4\{F(e^{i(2k+1)\frac\pi3}z,t_1^*)\}_4
    =
    2\overline{\{F(e^{i(2k+1)\frac\pi3}z,t_1^*)\}_2},
    \quad
    k=0,1,2.
\]
Hence all three functions satisfy relation~\eqref{eq:Stupin:Marty}.

\section{A system of equations for $\mathbb E_2(3)$}

Throughout what follows, unknowns interpreted as angles are chosen from the half-open interval $[0,2\pi)$ unless a
smaller set of values is specified.
This convention is necessary for assertions that the number of solutions is finite: otherwise each such solution would
generate infinitely many formally distinct solutions determining the same parameters of the same function.

\begin{Lemma}
    \label{lem:Stupin:twoCirclePointsSymmetricForm}
Let $z_1=e^{i\theta_1}$, $z_2=e^{i\theta_2}$, $0\leqslant\theta_1<\theta_2<2\pi$.
Then there exist unique $\mu\in\partial\Delta$ and $\varphi\in(0,\pi)$ such that $z_1=\mu e^{i\varphi}$, $z_2=\mu
e^{-i\varphi}$.
Moreover, $\mu=-e^{i\frac{\theta_1+\theta_2}2}$, $\varphi=\pi-\tfrac{\theta_2-\theta_1}2$, and $\mu^2=z_1 z_2$.
\end{Lemma}

\begin{proof}
Consider the system of equations $z_1=\mu e^{i\varphi}$, $z_2=\mu e^{-i\varphi}$.
Eliminating $\mu$ from the first equation gives the equivalent system
    $z_1=e^{2i\varphi}z_2$, $\mu=e^{-i\varphi}z_1$.
Its first equation implies that $\varphi=\pi k-\tfrac{\theta_2-\theta_1}2$, $k\in\mathbb Z$.
Normalizing $\varphi$ gives $\varphi=\pi-\tfrac{\theta_2-\theta_1}2$.
The inequalities $0<\theta_2-\theta_1<2\pi$ imply that $\varphi\in(0,\pi)$.
Substitution into the second equation of the system gives $\mu=-e^{i\frac{\theta_1+\theta_2}2}$.
Clearly, for $\varphi\in(0,\pi)$ this solution $(\mu, \varphi)$ is unique.

Clearly, $\mu\in\partial\Delta$, $\mu e^{i\varphi}=e^{i\theta_1}=z_1$, and
    $\mu e^{-i\varphi}=e^{i\theta_2}=z_2$.
Multiplying the last equalities gives $\mu^2=z_1 z_2$.
\end{proof}

Note that $\varphi$ determines the positions of $z_1$ and $z_2$ relative to $\mu$, while $\mu$ determines their common
rotation.

The following theorem establishes a necessary and sufficient condition for a function to belong to $\mathbb E_2(3)$.

\begin{Theorem}
    \label{th:Stupin:extremalTypeE3TwoTermsSystem}
The function $f\in\mathbb E_2(3)$ if and only if $f=e^{-h}$, where
    \begin{equation}
        \label{eq:Stupin:E3TwoTerms_h}
        h(z)
        =
        \alpha_1\frac{1+z_1 z}{1-z_1 z}+\alpha_2\frac{1+z_2 z}{1-z_2 z},
        \quad
        \alpha_1,\alpha_2>0,
        \quad
        z_{1,2}=e^{i\theta_{1,2}},
        \quad
        0\leqslant\theta_1<\theta_2<2\pi,
    \end{equation}
and there exist $\mu\in\partial\Delta$, $\varphi\in(0,\pi)$, $R\geqslant2$, $t>0$ such that $z_1=\mu e^{i\varphi}$,
$z_2=\mu e^{-i\varphi}$, and the vector $(t,R,\varphi,\gamma,\alpha_1,\alpha_2)$ is a solution of the system
    \begin{equation}
        \label{eq:Stupin:E3TwoTermsRealSystem}
        \begin{cases}
            t=\alpha_1+\alpha_2,
            \\
            -R\cos\gamma=2(t-2)\cos\varphi,
            \\
            R\sin\gamma=2(\alpha_1-\alpha_2)\sin\varphi,
            \\
            \big(\tfrac{R^2}2-1\big)\cos(2\gamma)=2(t-1)\cos(2\varphi),
            \\
            -\big(\tfrac{R^2}2-1\big)\sin(2\gamma)=2(\alpha_1-\alpha_2)\sin(2\varphi),
            \\
            \big(1-\tfrac{R^2}3\big)R\cos(3\gamma)
            =2\big(t-\tfrac23\big)\cos(3\varphi),
            \\
            -\big(1-\tfrac{R^2}3\big)R\sin(3\gamma)
            =2(\alpha_1-\alpha_2)\sin(3\varphi),
        \end{cases}
    \end{equation}
where $\gamma\coloneqq\arg\mu^3\in[0,2\pi)$.
Every solution $(t,R,\varphi,\gamma,\alpha_1,\alpha_2)\in\mathbb R^6$ of system~\eqref{eq:Stupin:E3TwoTermsRealSystem}
for which $t>0$, $R\geqslant2$, $\varphi\in(0,\pi)$, $\gamma\in[0,2\pi)$, and $\alpha_1,\alpha_2>0$ corresponds to three
pairwise distinct functions in $\mathbb E_2(3)$ related by admissible rotations in the variable plane $z$.
\end{Theorem}

\begin{proof}
\textbf{Necessity.}
    \textbf{1.}
Let $f\in\mathbb E_2(3)$.
Then the representation of the function $f=e^{-h}$ in the form~\eqref{eq:Stupin:E3TwoTerms_h} follows from the
definition of $\mathbb E_2(3)$.

By Lemma~\ref{lem:Stupin:twoCirclePointsSymmetricForm}, there exist $\mu\in\partial\Delta$ and $\varphi\in(0,\pi)$ such
that $z_1=\mu e^{i\varphi}$, $z_2=\mu e^{-i\varphi}$.
By Corollary~\ref{cor:Stupin:f_n_eq_f_0_norm_Psi2}, there exists a normalized polynomial $\Psi(z)=\psi_0+\psi_1 z+\psi_2
z^2 + \psi_3 z^3$ external for $H$, where $\psi_3=1$.
At the same time, $\Psi(z)=(z-z_1)(z-z_2)(z-w)$.
By Corollary~\ref{cor:Stupin:extremalTypeZerosOfH}, we have $|w|\geqslant1$.
Furthermore, $\psi_0=-z_1 z_2 w=-\mu^2 w$, whence $w=-\psi_0 \mu^{-2}$.
Hence,
    \begin{equation}
        \label{eq:Stupin:E3TwoTermsPsi}
        \Psi(z)
        =
        (z-z_1)(z-z_2)(z+\psi_0\mu^{-2})
        =
        \psi_0
        +(\mu^2-2\psi_0 \mu^{-1}\cos\varphi)z
        +(\psi_0 \mu^{-2}-2\mu\cos\varphi)z^2
        +z^3.
    \end{equation}
Thus the unknown $w$ is eliminated.
Comparing the coefficients of $z$ gives $\psi_1=\mu^2-2\psi_0 \mu^{-1}\cos\varphi$.
Comparing the coefficients of $z^2$ gives $\psi_2=\psi_0 \mu^{-2}-2\mu\cos\varphi$.
Formula~\eqref{eq:Stupin:rootsOfH_prod_arg} gives $\psi_0\geqslant1$.
Set $R\coloneqq \psi_0+\frac1{\psi_0}$, $\gamma\coloneqq\arg\mu^3\in[0,2\pi)$.
Note that $R\geqslant2$ by Theorem~\ref{th:Stupin:nullsOfP_against_H}.

    \textbf{2.}
Since $P=p_3\Psi$, formulas~\eqref{eq:Stupin:f_n_p_n} for $n=3$ give
    $\{f\}_0=|p_3|^2\psi_0$,
    \begin{equation}
        \label{eq:Stupin:preSys}
        \begin{split}
            \frac{\{f\}_1}{\{f\}_0}
            =&
            \frac{\psi_2\overline{\psi_0}+\psi_3\overline{\psi_1}}{\psi_0}
            =
            \mu(Re^{-i\gamma}-4\cos\varphi),\\
            \frac{\{f\}_2}{\{f\}_0}
            =&
            \frac{
                \psi_1\overline{\psi_0}
                +
                \psi_2\overline{\psi_1}
                +
                \psi_3\overline{\psi_2}
            }{\psi_0}
            =
            \mu^2\bigl(e^{-2i\gamma}+2+4\cos^2\varphi-4R\cos\varphi e^{-i\gamma}\bigr),\\
            \frac{\{f\}_3}{\{f\}_0}
            =&
            \frac{
                |\psi_0|^2
                +
                |\psi_1|^2
                +
                |\psi_2|^2
                +
                |\psi_3|^2
            }{\psi_0}
            =
            2R(1+2\cos^2\varphi)-8\cos\varphi\cos\gamma.
        \end{split}
    \end{equation}
We have also used the formulas for the coefficients $\psi_j$, $j=1,2,3$ found in the preceding step of this proof.

    \textbf{3.}
We derive the system of equations.
On the one hand,
    \[
        h(z)
        =
        -\ln f(z)
        =
        -\ln\{f\}_0
        -\tfrac{\{f\}_1}{\{f\}_0} z
        +\left(\tfrac{\{f\}_1^2}{2\{f\}_0^2}-\tfrac{\{f\}_2}{\{f\}_0}\right)z^2
        +\left(
        \tfrac{\{f\}_1\{f\}_2}{\{f\}_0^2}
        -
        \tfrac{\{f\}_1^3}{3\{f\}_0^3}
        -
        \tfrac{\{f\}_3}{\{f\}_0}
        \right)z^3
        +o(z^3).
    \]
Set $t\coloneqq\{h\}_0>0$.
Substituting formulas~\eqref{eq:Stupin:preSys} into the expressions for the coefficients of the function $h$ and
simplifying gives
    \begin{gather*}
        \{h\}_1=\mu(4\cos\varphi-Re^{-i\gamma}),\\
        \{h\}_2
        =
        \mu^2\left(
               \big(\tfrac{R^2}2-1\big)e^{-2i\gamma}+2\cos(2\varphi)
        \right),\\
        \{h\}_3
        =
        \mu^3\left(
            R\big(1-\tfrac{R^2}3\big)e^{-3i\gamma}
            +\tfrac43\cos(3\varphi)
        \right).
    \end{gather*}

On the other hand, the representation of the function $h$ in the form~\eqref{eq:Stupin:E3TwoTerms_h} implies the
equalities
    \begin{equation}
        \label{eq:Stupin:E3TwoTerms_h_k}
        \{h\}_0=\alpha_1+\alpha_2,
        \quad
        \{h\}_k=2\mu^k\bigl(t\cos(k\varphi)+i(\alpha_1-\alpha_2)\sin(k\varphi)\bigr),
        \quad
        k=1,2,3.
    \end{equation}

Equating the resulting representations of the coefficients $\{h\}_k$, $k=0,1,2,3$, and canceling the factor $\mu^k$ for
$k=1,2,3$ gives all equations of the system
    \begin{equation}
        \label{eq:Stupin:E3TwoTermsSystem}
        \begin{cases}
            t=\alpha_1+\alpha_2,
            \\
            -Re^{-i\gamma}
            =
            2\bigl((t-2)\cos\varphi+i(\alpha_1-\alpha_2)\sin\varphi\bigr),
            \\
            \big(\tfrac{R^2}2-1\big)e^{-2i\gamma}
            =
            2\bigl((t-1)\cos(2\varphi)+i(\alpha_1-\alpha_2)\sin(2\varphi)\bigr),
            \\
            \big(1-\tfrac{R^2}3\big)Re^{-3i\gamma}
            =
            2\bigl((t-\tfrac23)\cos(3\varphi)+i(\alpha_1-\alpha_2)\sin(3\varphi)\bigr).
        \end{cases}
    \end{equation}
Writing out the real and imaginary parts of the equations of system~\eqref{eq:Stupin:E3TwoTermsSystem} gives
system~\eqref{eq:Stupin:E3TwoTermsRealSystem}.

\textbf{Sufficiency.}
    \textbf{1.}
We construct the functions corresponding to a solution of the system and show that they are pairwise distinct.
Let $(t,R,\varphi,\gamma,\alpha_1,\alpha_2)$ be a solution of system~\eqref{eq:Stupin:E3TwoTermsSystem} satisfying the
conditions of the last assertion of this theorem.
Choose $\mu\in\partial\Delta$ so that $\mu^3=e^{i\gamma}$, and set
    $z_1=\mu e^{i\varphi}$, $z_2=\mu e^{-i\varphi}$.
Define $h$ by formula~\eqref{eq:Stupin:E3TwoTerms_h} and set $f=e^{-h}$.
Since $|z_1|=|z_2|=1$, the functions $h$ and $f$ are holomorphic in $\Delta$.

All three roots of the equation $\mu^3=e^{i\gamma}$ have the form
    $\mu e^{\frac{2\pi ki}3}$, $k=0,1,2$.
Replacing $\mu$ by $\mu e^{\frac{2\pi ki}3}$ rotates the points $z_1$ and $z_2$ in the variable plane $z$ through the
angle $\frac{2\pi k}3$ and replaces the function $f(z)$ by
    $f(e^{\frac{2\pi ki}3}z)$.
Thus the functions constructed are related by admissible rotations in the variable plane $z$.

We prove that the functions constructed are pairwise distinct.
If two of them coincided, then $f(e^{\frac{2\pi ki}3}z)=f(z)$ for some $k\in\{1,2\}$ and all $z\in\Delta$.
Since $e^{\frac{2\pi ki}3}\neq1$ and $e^{\frac{4\pi ki}3}\neq1$, comparing the coefficients of $z$ and $z^2$ immediately
gives $\{f\}_1=\{f\}_2=0$.
Similarly, for the function $h$ we have $\{h\}_1=\{h\}_2=0$.
Formulas~\eqref{eq:Stupin:E3TwoTerms_h_k}, the equality $\{h\}_1=0$, and the condition $t>0$ imply $\cos\varphi=0$.
Since $\varphi\in(0,\pi)$, we have $\varphi=\frac\pi2$.
Then the same formulas~\eqref{eq:Stupin:E3TwoTerms_h_k} give $\{h\}_2=-2t\mu^2\neq0$, contradicting $\{h\}_2=0$.
Hence all functions constructed are pairwise distinct.

    \textbf{2.}
Fix one of the functions constructed in the preceding step of this proof and check that this function $f$ belongs to
$\mathbb E_2(3)$.
Set $H(z)=\{f\}_3+2\{f\}_2 z+2\{f\}_1 z^2+2\{f\}_0 z^3$.

The first equation of system~\eqref{eq:Stupin:E3TwoTermsSystem} gives $\{h\}_0=t$.
Substituting the second, third, and fourth equations of this system into formulas~\eqref{eq:Stupin:E3TwoTerms_h_k}
gives, respectively,
    \[
        \begin{aligned}
            \{h\}_1&=\mu(4\cos\varphi-Re^{-i\gamma}),\\
            \{h\}_2&=\mu^2\bigl((\tfrac{R^2}2-1)e^{-2i\gamma}+2\cos(2\varphi)\bigr),\\
            \{h\}_3&=\mu^3\bigl(R(1-\tfrac{R^2}3)e^{-3i\gamma}+\tfrac43\cos(3\varphi)\bigr).
        \end{aligned}
    \]
Since $f=e^{-h}$, we have
    \[
        \{f\}_0=e^{-t},
        \quad
        \frac{\{f\}_1}{e^{-t}}=-\{h\}_1,
        \quad
        \frac{\{f\}_2}{e^{-t}}=\tfrac12\{h\}_1^2-\{h\}_2,
        \quad
        \frac{\{f\}_3}{e^{-t}}=-\tfrac16\{h\}_1^3+\{h\}_1\{h\}_2-\{h\}_3.
    \]
Substituting the coefficients $h$ found into these equalities for $\mu^3=e^{i\gamma}$ gives the expressions for the
coefficients $f$ on the right-hand sides of formulas~\eqref{eq:Stupin:preSys}.

Choose $\psi_0\geqslant1$ so that $R=\psi_0+\psi_0^{-1}$, and consider the polynomial $\Psi$ defined by
formula~\eqref{eq:Stupin:E3TwoTermsPsi}.
Set $P(z)\coloneqq\sqrt{\frac{e^{-t}}{\psi_0}}\Psi(z)=\sum\limits_{j=0}^3p_jz^j$.
Now construct the polynomial
    \[
        \widetilde H(z)
        \coloneqq
        \sum\limits_{j=0}^3|p_j|^2
        +2\sum\limits_{k=1}^3\left(\sum\limits_{j=0}^{3-k}p_{j+k}\overline{p_j}\right)z^k.
    \]
Since $p_3=\sqrt{\frac{e^{-t}}{\psi_0}}>0$, the constant coefficient $\widetilde H$ is positive.
By Theorem~\ref{th:Stupin:characterization_of_polynomials_from_C}, we have $\widetilde H\in C^3$.

We have $p_j=\sqrt{\frac{e^{-t}}{\psi_0}}\,\psi_j$, $j=0,1,2,3$.
Thus the sums $\sum\limits_{j=0}^{3-k}p_{j+k}\overline{p_j}$ for $k=0,1,2$ equal $e^{-t}$ multiplied by the
corresponding expressions~\eqref{eq:Stupin:preSys} for $\frac{\{f\}_k}{\{f\}_0}$, $k=3,2,1$.
For $k=3$ we have $p_3\overline{p_0}=e^{-t}=\{f\}_0$.
Hence all coefficients of $H$ and $\widetilde H$ coincide.
Thus $H=\widetilde H\in C^3$ and $\{f\}_3=\sum\limits_{j=0}^3|p_j|^2>0$.

Since $H=\widetilde H$, Assertion~\ref{st:Stupin:polynomialModulusIdentity} gives
    $\operatorname{Re}H(e^{i\theta})=|P(e^{i\theta})|^2$, $\theta\in\mathbb R$.
By~\eqref{eq:Stupin:E3TwoTermsPsi}, we have $P(z_j)=0$, $j=1,2$.
Hence $\operatorname{Re}H(z_j)=0$, $j=1,2$.
The points $z_1$ and $z_2$ are distinct, since $\varphi\in(0,\pi)$.
Note that the pairs $(\alpha_1,z_1)$ and $(\alpha_2,z_2)$ can be renumbered so that the arguments of the new points
$z_1$ and $z_2$ satisfy $0\leqslant\theta_1<\theta_2<2\pi$.
This does not change the functions $h$ and $f$.
Since $\alpha_1,\alpha_2>0$, we have $f\in\mathbb E_2(3)$.
This holds for all three functions constructed, since $\mu$ was chosen arbitrarily.
\end{proof}

The following lemma expresses the values of the unknowns $\alpha_1$ and $\alpha_2$ in terms of the other unknowns of
this system.

\begin{Lemma}
    \label{lem:Stupin:E3TwoTermsParameters}
Let $(t,R,\varphi,\gamma,\alpha_1,\alpha_2)\in\mathbb R^6$ and $\varphi\in(0,\pi)$.
The first and third equations of system~\eqref{eq:Stupin:E3TwoTermsRealSystem} are equivalent to the equations
    \begin{equation}
        \label{eq:Stupin:E3TwoTermsParameters}
        \alpha_1=\tfrac t2+\tfrac{R\sin\gamma}{4\sin\varphi},
        \quad
        \alpha_2=\tfrac t2-\tfrac{R\sin\gamma}{4\sin\varphi}.
    \end{equation}
\end{Lemma}

\begin{proof}
Since $\sin\varphi>0$, the third equation of system~\eqref{eq:Stupin:E3TwoTermsRealSystem} is equivalent to
$\alpha_1-\alpha_2=\tfrac{R\sin\gamma}{2\sin\varphi}$.
Together with the first equation $\alpha_1+\alpha_2=t$ of this system, it is equivalent to
equations~\eqref{eq:Stupin:E3TwoTermsParameters}, obtained by adding and subtracting these two equations and then
dividing by $2$.
\end{proof}

\section{A system of equations for $\mathbb E_2(3)$ when $\varphi=\frac\pi2$}

\begin{Assertion}
    \label{st:Stupin:extremalTypeE3TwoTermsZeroC}
The case $\varphi=\frac\pi2$ of Theorem~\ref{th:Stupin:extremalTypeE3TwoTermsSystem} corresponds to precisely the
following functions in $\mathbb E_2(3)$\textnormal:
    \[
        f_{\frac{\sqrt6}4}\big(e^{i\frac{2\pi k}3}z\big),
        \quad
        k=0,1,2,
        \quad
\text{where}
        \quad
        f_{\frac{\sqrt6}4}
        \coloneqq
        e^{-h},
        \quad
        h(z)
        \coloneqq
        \big(1+\tfrac{\sqrt6}4\big)\dfrac{1+z}{1-z} + \big(1-\tfrac{\sqrt6}4\big)\dfrac{1-z}{1+z}.
    \]
\end{Assertion}

\begin{proof}
By Theorem~\ref{th:Stupin:extremalTypeE3TwoTermsSystem}, every solution $(t,R,\frac\pi2,\gamma,\alpha_1,\alpha_2)$ of
system~\eqref{eq:Stupin:E3TwoTermsRealSystem} corresponds to three functions in $\mathbb E_2(3)$.

We solve system~\eqref{eq:Stupin:E3TwoTermsRealSystem}.
Its second equation gives $\cos\gamma=0$.
Its third equation gives $2(\alpha_1-\alpha_2)=R\sin\gamma$.
The equality $\cos\gamma=0$ implies that $\sin\gamma=\pm1$.
The fifth and sixth equations are then satisfied identically.
Using $\cos(2\gamma)=\cos(2\varphi)=-1$, the fourth equation of system~\eqref{eq:Stupin:E3TwoTermsRealSystem} gives
$\tfrac{R^2}2-1=2(t-1)$.
Using $\sin(3\gamma)=-\sin\gamma$ and $\sin(3\varphi)=-1$, its seventh equation gives
    $\big(1-\tfrac{R^2}3\big)R\sin\gamma=-2(\alpha_1-\alpha_2)$.
The resulting equation and the equation $2(\alpha_1-\alpha_2)=R\sin\gamma$ give $1-\tfrac{R^2}3=-1$.
Hence $R^2=6$.
Substituting $R^2=6$ into the equation $\tfrac{R^2}2-1=2(t-1)$ gives $t=2$.
Since $R\geqslant2$, we have $R=\sqrt6$.
Since $\varphi=\frac\pi2$, we have $z_1=-z_2=i\mu$ by Theorem~\ref{th:Stupin:extremalTypeE3TwoTermsSystem}.
Formulas~\eqref{eq:Stupin:E3TwoTermsParameters} give
    \[
        (\alpha_1,\alpha_2)
        =
        \begin{cases}
            \big(1+\tfrac{\sqrt6}4,1-\tfrac{\sqrt6}4\big), & \sin\gamma=1,\\
            \big(1-\tfrac{\sqrt6}4,1+\tfrac{\sqrt6}4\big), & \sin\gamma=-1
        \end{cases}.
    \]
By Theorem~\ref{th:Stupin:extremalTypeE3TwoTermsSystem}, to determine the function $f$ up to admissible rotations in the
variable plane $z$, it suffices for each value $\gamma$ to choose one value $\mu$ satisfying $\mu^3=e^{i\gamma}$.
Consider two cases: $\sin\gamma=1$ and $\sin\gamma=-1$.

If $\sin\gamma=1$, then $\gamma=\frac\pi2$, since $\gamma\in[0,2\pi)$.
The number $\mu=-i$ satisfies $\mu^3=e^{i\gamma}$.
Moreover, $z_1=1$.
Also, $z_2=-1$.
Together with the first pair of values $\alpha_1,\alpha_2$, this gives the function $f_{\frac{\sqrt6}4}$.

If $\sin\gamma=-1$, then $\gamma=\frac{3\pi}2$, since $\gamma\in[0,2\pi)$.
The number $\mu=i$ satisfies $\mu^3=e^{i\gamma}$.
Moreover, $z_1=-1$.
Also, $z_2=1$.
Substituting these values and the second pair $\alpha_1,\alpha_2$ into representation~\eqref{eq:Stupin:E3TwoTerms_h} of
the function $h$ and interchanging the terms again gives the function $f_{\frac{\sqrt6}4}$.

Thus both values $\gamma$ give the same function up to admissible rotations in the variable plane $z$.
Hence,
    $
        f(z)
        =
        e^{-h(e^{i\frac{2\pi k}3}z)}
        =
        f_{\frac{\sqrt6}4}\big(e^{i\frac{2\pi k}3}z\big)$,
        $k=0,1,2$.

All transformations used to solve system~\eqref{eq:Stupin:E3TwoTermsRealSystem} are equivalent, since $R>0$ and
$\sin\gamma=\pm1$, while formulas~\eqref{eq:Stupin:E3TwoTermsParameters} for $\alpha_1,\alpha_2$ are equivalent to its
first and third equations by Lemma~\ref{lem:Stupin:E3TwoTermsParameters}.
Hence both vectors found are solutions of the system, with $\alpha_1,\alpha_2>0$.
Thus, by Theorem~\ref{th:Stupin:extremalTypeE3TwoTermsSystem}, the function $f_{\frac{\sqrt6}4}$ and all its admissible
rotations in the variable plane $z$ belong to $\mathbb E_2(3)$.
\end{proof}

\smallskip
We have
\begin{multline}
    \label{eq:Stupin:fSqrt6Over4Expansion}
    f_{\frac{\sqrt6}4}(z)
    =
    e^{-2}-\sqrt6e^{-2}z-e^{-2}z^2
    +2\sqrt6e^{-2}z^3-\tfrac12e^{-2}z^4+o(z^4)
    \approx \\ \approx
    0.135-0.332z-0.135z^2+0.663z^3-0.068z^4+o(z^4).
\end{multline}
Expansion~\eqref{eq:Stupin:fSqrt6Over4Expansion} implies that
$
    4\{f_{\frac{\sqrt6}4}\}_4
    =
    -2e^{-2}
    =
    2\overline{\{f_{\frac{\sqrt6}4}\}_2}$.
Hence the function $f_{\frac{\sqrt6}4}$ and all its admissible rotations in the variable plane $z$,
$f_{\frac{\sqrt6}4}(e^{i\frac{2\pi k}3}z)$, $k=0,1,2$, satisfy relation~\eqref{eq:Stupin:Marty}.

\section{Solutions of the system for $\mathbb E_2(3)$ when $\varphi\in(0,\pi)\setminus\{\frac\pi2\}$}

\begin{Lemma}
    \label{lem:Stupin:E3TwoTermsNonzeroCPreliminary}
Let $(t,R,\varphi,\gamma,\alpha_1,\alpha_2)\in\mathbb R^6$ be a solution of
system~\eqref{eq:Stupin:E3TwoTermsRealSystem}, with $t>0$, $R\geqslant2$, $\varphi\in(0,\pi)\setminus\{\frac\pi2\}$, and
$\gamma\in[0,2\pi)$.
Then $\sin\gamma\neq0$.
\end{Lemma}

\begin{proof}
    \textbf{1.}
To exclude the case $\sin\gamma=0$, we derive a system of equations in the unknowns $t$ and $R$.
Suppose $\sin\gamma=0$.

By assumption, $\varphi\in(0,\pi)\setminus\{\frac\pi2\}$.
Hence $\sin\varphi>0$.
Also, $\cos\varphi\neq0$.
The third equation of system~\eqref{eq:Stupin:E3TwoTermsRealSystem} for $\sin\gamma=0$ implies that $\alpha_1=\alpha_2$.
The condition $\sin\gamma=0$ implies that $\cos\gamma=\pm1$.
The second equation of the system has the form $-R\cos\gamma=2(t-2)\cos\varphi$.
Hence $t\neq2$.
Solving the second equation for $\cos\varphi$ gives
    \begin{equation}
        \label{eq:Stupin:E3TwoTermsNonzeroC_cosPhi}
        \cos\varphi=-\tfrac{R\cos\gamma}{2(t-2)}.
    \end{equation}
Since $\cos\gamma=\pm1$, we have $\cos(2\gamma)=1$.
Also, $\cos(3\gamma)=\cos\gamma$.
Using the identities
    \[
        \cos(2\varphi)=2\cos^2\varphi-1,
        \quad
        \cos(3\varphi)=4\cos^3\varphi-3\cos\varphi,
    \]
we write the fourth and sixth equations of the system as
    \[
        \tfrac{R^2}2-1=2(t-1)(2\cos^2\varphi-1),
        \quad
        R\big(1-\tfrac{R^2}3\big)\cos\gamma
        =
        2\big(t-\tfrac23\big)(4\cos^3\varphi-3\cos\varphi).
    \]
Substituting expression~\eqref{eq:Stupin:E3TwoTermsNonzeroC_cosPhi} for $\cos\varphi$ into these equations and
simplifying gives the following system in the unknowns $t$ and $R$:
    \begin{equation}
        \label{eq:Stupin:E3TwoTermsR2}
        \begin{cases}
            R^2(t^2-6t+6)+2(t-2)^2(2t-3)=0,
            \\
            R^2(t^3-6t^2+9t-6)+6t(t-2)^2=0.
        \end{cases}
    \end{equation}

    \textbf{2.}
We eliminate $R^2$ from system~\eqref{eq:Stupin:E3TwoTermsR2}.
To this end, multiply the first equation of system~\eqref{eq:Stupin:E3TwoTermsR2} by $t^3-6t^2+9t-6$, the second by
$t^2-6t+6$, and subtract the second equation from the first.
We have $2(t-2)^2(2t^4-18t^3+54t^2-57t+18)=0$.
Using $t\neq2$, we obtain
    \begin{equation}
        \label{eq:Stupin:E3TwoTermsEq}
        2t^4-18t^3+54t^2-57t+18=0.
    \end{equation}

    \textbf{3.}
We narrow the set of possible values of the unknown $t$.
Since $|\cos\varphi|<1$ and $R\geqslant2$,
    $
        1
        >
        |\cos\varphi|
        =
        \tfrac{R}{2|t-2|}
        \geqslant
        \tfrac1{|t-2|}$.
Using the condition $t>0$, we obtain $t\in(0,1)\cup(3,+\infty)$.

The first equation of system~\eqref{eq:Stupin:E3TwoTermsR2} gives $R^2-4=-2t\tfrac{2t^2-9t+8}{t^2-6t+6}$.
But for $0<t<1$ we have
    \[
        2t^2-9t+8=1+(1-t)(7-2t)>0,
        \quad
        t^2-6t+6=1+(1-t)(5-t)>0.
    \]
whence $R^2<4$, which is impossible for $R\geqslant2$.
Hence $t>3$.

    \textbf{4.}
We obtain a contradiction by showing that equation~\eqref{eq:Stupin:E3TwoTermsEq} has no roots for $t>3$.
Denote the left-hand side of equation~\eqref{eq:Stupin:E3TwoTermsEq} by $Q$.
Expanding the parentheses and collecting like terms readily shows that $Q(t)=2(t-3)^4+6(t-3)^3-3(t-3)+9$.
We now prove that $Q(t)>0$ for $t>3$.
Indeed, for $3<t\leqslant4$ we have $9-3(t-3)\geqslant6$, and the other terms are positive.
For $t\geqslant4$ we have $6(t-3)^3-3(t-3)>0$, and the other terms are again positive.
Hence equation~\eqref{eq:Stupin:E3TwoTermsEq} has no roots for $t>3$.
This contradicts the fact that $t>3$ satisfies this equation.
Hence $\sin\gamma\neq0$.
\end{proof}

\begin{Lemma}
    \label{lem:Stupin:E3TwoTermsNonzeroCFormulas}
Let $(t,R,\varphi,\gamma,\alpha_1,\alpha_2)\in\mathbb R^6$, $t>0$, $R\geqslant2$,
$\varphi\in(0,\pi)\setminus\{\frac\pi2\}$, $\gamma\in[0,2\pi)$, and $\sin\gamma\neq0$.
The first five equations of system~\eqref{eq:Stupin:E3TwoTermsRealSystem} are equivalent to the condition
$3<t\leqslant4$, equations~\eqref{eq:Stupin:E3TwoTermsParameters}, and the equations
    \begin{align}
        \label{eq:Stupin:E3TwoTermsCosGamma}
        \cos\gamma
        =&
        -\tfrac{2R\cos\varphi}{R^2-2},
        \\
        \label{eq:Stupin:E3TwoTermsR2ByT}
        R^2
        =&
        2+\tfrac2{t-3},
        \\
        \label{eq:Stupin:E3TwoTermsCos2PhiByT}
        \cos^2\varphi
        =&
        \tfrac{2t^2-8t+5}{4(t-3)},
        \\
        \label{eq:Stupin:E3TwoTermsCos2GammaByT}
        \cos^2\gamma
        =&
        \tfrac{(t-2)(2t^2-8t+5)}2,
    \end{align}
Moreover, $\alpha_1\neq\alpha_2$.
\end{Lemma}

\begin{proof}
We have $\sin\varphi>0$, $\cos\varphi\neq0$, $\sin\gamma\neq0$, and $R^2-2>0$.
By Lemma~\ref{lem:Stupin:E3TwoTermsParameters}, the first and third equations of
system~\eqref{eq:Stupin:E3TwoTermsRealSystem} are equivalent to equations~\eqref{eq:Stupin:E3TwoTermsParameters}.
Moreover, the third equation of system~\eqref{eq:Stupin:E3TwoTermsRealSystem} implies
    $\alpha_1-\alpha_2=\tfrac{R\sin\gamma}{2\sin\varphi}\neq0$.

Solving the third equation of system~\eqref{eq:Stupin:E3TwoTermsRealSystem} for $\sin\gamma$ and substituting the
resulting expression into its fifth equation, after cancellation we obtain
    $\big(\tfrac{R^2}2-1\big)\cos\gamma=-R\cos\varphi$.
These transformations are reversible, since $R>0$, $\sin\varphi>0$, and $\alpha_1-\alpha_2\neq0$.
Since $R^2-2>0$, the fifth equation, with the third retained, is equivalent to
equation~\eqref{eq:Stupin:E3TwoTermsCosGamma}.

The second equation of system~\eqref{eq:Stupin:E3TwoTermsRealSystem} has the form
    $-R\cos\gamma=2(t-2)\cos\varphi$.
Substituting expression~\eqref{eq:Stupin:E3TwoTermsCosGamma} for $\cos\gamma$ into it and canceling $\cos\varphi\neq0$
gives the equivalent equation $t=2+\tfrac{R^2}{R^2-2}$.
Since $t=2+\tfrac{R^2}{R^2-2}=3+\tfrac2{R^2-2}$ and $R\geqslant2$, we have $3<t\leqslant4$.
For $t>3$, the equation $t=3+\tfrac2{R^2-2}$ is equivalent to equation~\eqref{eq:Stupin:E3TwoTermsR2ByT}.

Equation~\eqref{eq:Stupin:E3TwoTermsR2ByT} gives $t=\frac{4-3R^2}{2-R^2}$.
Substituting this expression for $t$ and expression~\eqref{eq:Stupin:E3TwoTermsCosGamma} for $\cos\gamma$ into the
fourth equation of system~\eqref{eq:Stupin:E3TwoTermsRealSystem}, using the identity $\cos(2u)=2\cos^2u-1$, transforms
it into
    $-4\cos^2\varphi+4+\tfrac4{R^2-2}-\tfrac{R^2-2}2=0$.
Since $R^2-2>0$, the substitution $t=3+\tfrac2{R^2-2}$ is admissible, and the fourth equation of
system~\eqref{eq:Stupin:E3TwoTermsRealSystem} is equivalent to the equation
    \begin{equation}
        \label{eq:Stupin:E3TwoTermsCos2Phi}
        \cos^2\varphi
        =
        1+\tfrac1{R^2-2}-\tfrac{R^2-2}8.
    \end{equation}
Substituting $R^2-2=\tfrac2{t-3}$ from equation~\eqref{eq:Stupin:E3TwoTermsR2ByT} for $t>3$ transforms
equation~\eqref{eq:Stupin:E3TwoTermsCos2Phi} into the equivalent equation~\eqref{eq:Stupin:E3TwoTermsCos2PhiByT}.

Thus the first five equations of system~\eqref{eq:Stupin:E3TwoTermsRealSystem} are equivalent to the condition
$3<t\leqslant4$, equations~\eqref{eq:Stupin:E3TwoTermsParameters},
and~\eqref{eq:Stupin:E3TwoTermsCosGamma}--\eqref{eq:Stupin:E3TwoTermsCos2PhiByT}.
Squaring equation~\eqref{eq:Stupin:E3TwoTermsCosGamma} and substituting the expressions for $R^2$ and $\cos^2\varphi$
from equations~\eqref{eq:Stupin:E3TwoTermsR2ByT},~\eqref{eq:Stupin:E3TwoTermsCos2PhiByT} gives
equation~\eqref{eq:Stupin:E3TwoTermsCos2GammaByT}.
Adding this consequence to the resulting equations does not change the solution set.
\end{proof}

\begin{Lemma}
    \label{lem:Stupin:E3TwoTermsTPolynomialUniqueRoot}
The equation
    \begin{equation}
        \label{eq:Stupin:E3TwoTermsTPolynomial}
        4t^4-44t^3+168t^2-258t+129=0.
    \end{equation}
has exactly one root on the interval $(3,4]$, and this root belongs to the interval $\big(\tfrac{23}7,\tfrac{10}3\big)$.
\end{Lemma}

\begin{proof}
Denote the left-hand side of equation~\eqref{eq:Stupin:E3TwoTermsTPolynomial} by $Q=Q(t)$.
We have $Q(0)=129$, $Q(1)=-1$, $Q(2)=-3$, $Q\big(\tfrac52\big)=\tfrac{11}4$,
    $Q\big(\tfrac{23}7\big)=\tfrac{1023}{2401}$, $Q\big(\tfrac{10}3\big)=-\tfrac{11}{81}$,
    $Q(4)=-7$, $Q\big(\tfrac92\big)=\tfrac34$.
By the Bolzano-Cauchy intermediate value theorem, the polynomial $Q$ has at least one root in each of the intervals
$(0,1)$, $\big(2,\tfrac52\big)$, $\big(\tfrac{23}7,\tfrac{10}3\big)$, and $\big(4,\tfrac92\big)$.
These roots are distinct, and the degree of $Q$ is $4$.
Hence it has no other roots.
Thus the unique root on the half-open interval $(3,4]$ belongs to the interval $\big(\tfrac{23}7,\tfrac{10}3\big)$.
\end{proof}

\begin{Lemma}
    \label{lem:Stupin:E3TwoTermsConstantsEstimates}
Let $t_2^*$ be the unique root of equation~\eqref{eq:Stupin:E3TwoTermsTPolynomial} on the half-open interval $(3,4]$,

        $R^*\coloneqq\sqrt{2+\tfrac2{t_2^*-3}}$,

        $c^*\coloneqq\sqrt{\tfrac{2(t_2^*)^2-8t_2^*+5}{4(t_2^*-3)}}$,

        $\varphi^*\coloneqq\arccos c^*$,

        $\gamma^*
        \coloneqq
        \arccos\bigl(-R^*c^*(t_2^*-3)\bigr)$.

Then
        $\tfrac{23}7<t_2^*<\tfrac{10}3$,

        $8<(R^*)^2<9$,

        $\varphi^*\in(0,\frac\pi2)$,

        $\gamma^*\in(\frac\pi2,\pi)$,

        $\sin^2\varphi^*>\tfrac7{12}$,

        $0<\tfrac{R^*\sin\gamma^*}{2\sin\varphi^*}<t_2^*$.
\end{Lemma}

\begin{proof}
We have $(R^*)^2=2+\tfrac2{t_2^*-3}$ by assumption and $\tfrac{23}7<t_2^*<\tfrac{10}3$ by
Lemma~\ref{lem:Stupin:E3TwoTermsTPolynomialUniqueRoot}.
Hence $\tfrac27<t_2^*-3<\tfrac13$.
Thus $6<\tfrac2{t_2^*-3}<7$, whence $8<(R^*)^2<9$.

The function $\tfrac{2t^2-8t+5}{4(t-3)}$ is strictly increasing for $t>3$.
Thus $\tfrac{23}7<t_2^*<\tfrac{10}3$ implies $\tfrac{15}{56}<(c^*)^2<\tfrac5{12}$.
This estimation implies that $c^*\in(0,1)$.
Hence $\varphi^*\in(0,\pi/2)$.
Since $\cos\varphi^*=c^*$, the estimation $(c^*)^2<\tfrac5{12}$ implies
    $\sin^2\varphi^*>\tfrac7{12}$.

The estimations already obtained imply
    $0<(R^*c^*(t_2^*-3))^2<9\cdot\tfrac5{12}\cdot\tfrac1{3^2}=\tfrac5{12}<1$.
Since $R^*c^*(t_2^*-3)>0$, this gives $\gamma^*\in(\pi/2,\pi)$.

The conditions $\sin^2\varphi^*>\tfrac7{12}$ and $\sin\varphi^*>0$ imply
    $\tfrac1{2\sin\varphi^*}<\sqrt{\tfrac37}$.
The conditions $(R^*)^2<9$ and $R^*>0$ imply $R^*<3$.
Since $\gamma^*\in(\frac\pi2,\pi)$, we have $0<\sin\gamma^*<1$.
Hence $0<\tfrac{R^*\sin\gamma^*}{2\sin\varphi^*}<3\sqrt{\tfrac37}<3<t_2^*$.
\end{proof}

\begin{Assertion}
    \label{st:Stupin:extremalTypeE3TwoTermsNonzeroCSystem}
Let $(t,R,\varphi,\gamma,\alpha_1,\alpha_2)\in\mathbb R^6$ be a solution of
system~\eqref{eq:Stupin:E3TwoTermsRealSystem}, with $t>0$, $R\geqslant2$, $\varphi\in(0,\pi)\setminus\{\frac\pi2\}$, and
$\gamma\in[0,2\pi)$.
Then $t=t_2^*$, $R=R^*$, and
    \begin{equation}
        \label{eq:Stupin:E3TwoTermsCosPair}
        (\cos\varphi,\cos\gamma)
        \in
        \{(c^*,\cos\gamma^*),(-c^*,-\cos\gamma^*)\}.
    \end{equation}
The values of the unknowns $\alpha_1,\alpha_2$ are determined by formulas~\eqref{eq:Stupin:E3TwoTermsParameters}.
Here $t_2^*$, $R^*$, $c^*$, and $\gamma^*$ are defined in Lemma~\ref{lem:Stupin:E3TwoTermsConstantsEstimates}.
Conversely, for $\varphi\in(0,\pi)\setminus\{\frac\pi2\}$ and $\gamma\in[0,2\pi)$, the equalities $t=t_2^*$, $R=R^*$ and
formulas~\eqref{eq:Stupin:E3TwoTermsCosPair},~\eqref{eq:Stupin:E3TwoTermsParameters} determine four solutions of
system~\eqref{eq:Stupin:E3TwoTermsRealSystem}, with $\alpha_1,\alpha_2>0$ for each of them.
\end{Assertion}

\begin{proof}
    \textbf{1.}
By Lemma~\ref{lem:Stupin:E3TwoTermsNonzeroCPreliminary}, we have $\sin\gamma\neq0$.
Denote the left-hand side of equation~\eqref{eq:Stupin:E3TwoTermsTPolynomial} by $Q(t)$.
We show that for $\sin\gamma\neq0$ system~\eqref{eq:Stupin:E3TwoTermsRealSystem} is equivalent to the system
    \begin{equation}
        \label{eq:Stupin:E3TwoTermsReducedSystem}
        \begin{cases}
            t=\alpha_1+\alpha_2,
            \\
            -R\cos\gamma=2(t-2)\cos\varphi,
            \\
            R\sin\gamma=2(\alpha_1-\alpha_2)\sin\varphi,
            \\
            \big(\tfrac{R^2}2-1\big)\cos(2\gamma)=2(t-1)\cos(2\varphi),
            \\
            -\big(\tfrac{R^2}2-1\big)\sin(2\gamma)=2(\alpha_1-\alpha_2)\sin(2\varphi),
            \\
            Q(t)=0.
        \end{cases}
    \end{equation}

Consider a vector $(t,R,\varphi,\gamma,\alpha_1,\alpha_2)$ satisfying the first five equations of
system~\eqref{eq:Stupin:E3TwoTermsRealSystem} and the constraints in the statement, with $\sin\gamma\neq0$.
Then, by Lemma~\ref{lem:Stupin:E3TwoTermsNonzeroCFormulas},
equations~\eqref{eq:Stupin:E3TwoTermsCosGamma}--\eqref{eq:Stupin:E3TwoTermsCos2GammaByT} hold.
Using the identity $\cos(3x)=(4\cos^2 x-3)\cos x$ and these formulas, we have
    \begin{equation}
        \label{eq:Stupin:E3TwoTermsCosThirdResidual}
        R\big(1-\tfrac{R^2}3\big)\cos(3\gamma)
        -2\big(t-\tfrac23\big)\cos(3\varphi)
        =
        -\tfrac{\cos\varphi}3\left(2+\tfrac2{t-3}\right)Q(t).
    \end{equation}
Using the third equation of system~\eqref{eq:Stupin:E3TwoTermsReducedSystem}, the identity $\sin(3x)=(4\cos^2x-1)\sin
x$, and the same formulas, we have
    \begin{equation}
        \label{eq:Stupin:E3TwoTermsSinThirdResidual}
        -R\big(1-\tfrac{R^2}3\big)\sin(3\gamma)
        -2(\alpha_1-\alpha_2)\sin(3\varphi)
        =
        -\tfrac{R\sin\gamma}{3(t-3)}Q(t).
    \end{equation}

Since $R\geqslant2$, $\sin\gamma\neq0$, and $t>3$, the factor $-\tfrac{R\sin\gamma}{3(t-3)}$ is nonzero.
Thus equality~\eqref{eq:Stupin:E3TwoTermsSinThirdResidual} shows that the seventh equation of
system~\eqref{eq:Stupin:E3TwoTermsRealSystem} is equivalent to $Q(t)=0$.
For $Q(t)=0$ the sixth equation is satisfied by equality~\eqref{eq:Stupin:E3TwoTermsCosThirdResidual}.
The first five equations of systems~\eqref{eq:Stupin:E3TwoTermsRealSystem} and~\eqref{eq:Stupin:E3TwoTermsReducedSystem}
coincide.
Thus these systems are equivalent for $\sin\gamma\neq0$.

    \textbf{2.}
Hence $Q(t)=0$ by formula~\eqref{eq:Stupin:E3TwoTermsReducedSystem} and $3<t\leqslant4$ by
Lemma~\ref{lem:Stupin:E3TwoTermsNonzeroCFormulas}.
Thus $t=t_2^*$ by Lemma~\ref{lem:Stupin:E3TwoTermsTPolynomialUniqueRoot}.

    \textbf{3.}
We recover the values of the remaining unknowns.
By Lemma~\ref{lem:Stupin:E3TwoTermsNonzeroCFormulas}, equation~\eqref{eq:Stupin:E3TwoTermsCosGamma} holds.
By the same lemma, equation~\eqref{eq:Stupin:E3TwoTermsR2ByT} holds.
Substituting $t=t_2^*$ into formula~\eqref{eq:Stupin:E3TwoTermsR2ByT} gives $R=R^*$.
Substituting $t=t_2^*$ into formula~\eqref{eq:Stupin:E3TwoTermsCos2PhiByT} gives
    $|\cos\varphi|=c^*$.
Since $\varphi,\varphi^*\in(0,\pi)$, this implies that $\sin\varphi=\sin\varphi^*$.
Substituting $t=t_2^*$ into formula~\eqref{eq:Stupin:E3TwoTermsCos2GammaByT} gives
    $|\cos\gamma|=|\cos\gamma^*|$.
Since $\gamma^*\in(\pi/2,\pi)$, this implies that $|\sin\gamma|=\sin\gamma^*$.

By the definition of $R^*$, we have $(R^*)^2-2=\tfrac2{t_2^*-3}$.
Thus the equality $|\cos\varphi|=c^*$ and formula~\eqref{eq:Stupin:E3TwoTermsCosGamma} imply
condition~\eqref{eq:Stupin:E3TwoTermsCosPair}.
By Lemma~\ref{lem:Stupin:E3TwoTermsParameters}, the values of the unknowns $\alpha_1,\alpha_2$ are determined by
formulas~\eqref{eq:Stupin:E3TwoTermsParameters}.

    \textbf{4.}
We prove the converse.
On the interval $(0,\pi)$ each value $\cos\varphi=\pm c^*$ corresponds to one value $\varphi$, and on the interval
$[0,2\pi)$ each value $\cos\gamma=\pm\cos\gamma^*$ corresponds to two values $\gamma$.
By the reduction formulas,
    \[
        \cos(\pi-\varphi^*)=-\cos\varphi^*,
        \quad
        \cos(2\pi-\gamma^*)=\cos\gamma^*,
        \quad
        \cos(\pi-\gamma^*)=\cos(\pi+\gamma^*)=-\cos\gamma^*.
    \]
Since $0<c^*<1$ and $\gamma^*\in(\frac\pi2,\pi)$, condition~\eqref{eq:Stupin:E3TwoTermsCosPair} on the pair
$(\cos\varphi,\cos\gamma)$ determines the following four pairs of values of the unknowns:
    \[
        (\varphi,\gamma)
        \in
        \{(\varphi^*,\gamma^*),(\varphi^*,2\pi-\gamma^*),
        (\pi-\varphi^*,\pi-\gamma^*),(\pi-\varphi^*,\pi+\gamma^*)\}.
    \]
For each of these four distinct pairs, the equalities $t=t_2^*$, $R=R^*$ and
formulas~\eqref{eq:Stupin:E3TwoTermsParameters} uniquely determine the remaining unknowns.
For each of them $\varphi\in(0,\pi)\setminus\{\frac\pi2\}$, $\sin\varphi=\sin\varphi^*$, and
$|\sin\gamma|=\sin\gamma^*$.
The equalities $t=t_2^*$ and $R=R^*$ and Lemma~\ref{lem:Stupin:E3TwoTermsConstantsEstimates} imply $t>0$ and $R>2$.
Formulas~\eqref{eq:Stupin:E3TwoTermsParameters} give
    \[
        \alpha_1+\alpha_2=t_2^*,
        \quad
        |\alpha_1-\alpha_2|
        =
        \tfrac{R^*\sin\gamma^*}{2\sin\varphi^*}.
    \]
By Lemma~\ref{lem:Stupin:E3TwoTermsConstantsEstimates}, we have
    $|\alpha_1-\alpha_2|<t_2^*$.
Hence,
    \[
        2\alpha_1
        =
        t_2^*+\alpha_1-\alpha_2
        \geqslant
        t_2^*-|\alpha_1-\alpha_2|
        >
        0,
        \quad
        2\alpha_2
        =
        t_2^*-\alpha_1+\alpha_2
        \geqslant
        t_2^*-|\alpha_1-\alpha_2|
        >
        0.
    \]

By the definitions of $R^*$, $c^*$, and $\gamma^*$, the values found satisfy the equations
    \eqref{eq:Stupin:E3TwoTermsCosGamma}--\eqref{eq:Stupin:E3TwoTermsCos2GammaByT},
and $|\sin\gamma|=\sin\gamma^*>0$.
By Lemma~\ref{lem:Stupin:E3TwoTermsNonzeroCFormulas}, the first five equations of
system~\eqref{eq:Stupin:E3TwoTermsReducedSystem} are equivalent to these equations and
equations~\eqref{eq:Stupin:E3TwoTermsParameters} for $3<t\leqslant4$.
Since $Q(t_2^*)=0$, all four vectors solve the reduced system, and, by the equivalence proved in step 1 of this proof,
also system~\eqref{eq:Stupin:E3TwoTermsRealSystem}.
\end{proof}

\section{Functions corresponding to solutions of the system for $\mathbb E_2(3)$ when
$\varphi\in(0,\pi)\setminus\{\frac\pi2\}$}

\begin{Assertion}
    \label{st:Stupin:extremalTypeE3TwoTermsNonzeroC}
In the notation of Theorem~\ref{th:Stupin:extremalTypeE3TwoTermsSystem}, the solutions of
system~\eqref{eq:Stupin:E3TwoTermsRealSystem} satisfying $t>0$, $R\geqslant2$,
$\varphi\in(0,\pi)\setminus\{\frac\pi2\}$, and $\gamma\in[0,2\pi)$ determine precisely the following functions in
$\mathbb E_2(3)$\textnormal:
        $f_2(e^{2\pi ik/3}z)$,
        $f_2^*(e^{2\pi ik/3}z)$,
        $k=0,1,2$,
where
        $f_2(z)
        \coloneqq
        \exp\bigg\{
        -\alpha_1
        \dfrac{1+e^{i(\frac{\gamma^*}3+\varphi^*)}z}
        {1-e^{i(\frac{\gamma^*}3+\varphi^*)}z}
        -\alpha_2
        \dfrac{1+e^{i(\frac{\gamma^*}3-\varphi^*)}z}
        {1-e^{i(\frac{\gamma^*}3-\varphi^*)}z}
        \bigg\}$,
$f_2^*(z)\coloneqq\overline{f_2(\overline z)}$, $\alpha_1,\alpha_2$ are determined by
formulas~\eqref{eq:Stupin:E3TwoTermsParameters}, $(t,R,\varphi,\gamma)=(t_2^*,R^*,\varphi^*,\gamma^*)$.
The numbers $t_2^*$, $R^*$, $\gamma^*$, and $\varphi^*$ are defined in
Lemma~\ref{lem:Stupin:E3TwoTermsConstantsEstimates}.
\end{Assertion}

\begin{proof}
Let $g\in\mathbb E_2(3)$ and $\varphi\in(0,\pi)\setminus\{\frac\pi2\}$.
By Theorem~\ref{th:Stupin:extremalTypeE3TwoTermsSystem}, the vector
$(t,R,\varphi,\gamma,\widetilde{\alpha}_1,\widetilde{\alpha}_2)$ is a solution of
system~\eqref{eq:Stupin:E3TwoTermsRealSystem}.
By Assertion~\ref{st:Stupin:extremalTypeE3TwoTermsNonzeroCSystem}, we have
    \[
        t=t_2^*,
        \quad
        R=R^*,
        \quad
        (\cos\varphi,\cos\gamma)
        \in
        \{(c^*,\cos\gamma^*),(-c^*,-\cos\gamma^*)\}.
    \]
For the values $\alpha_1,\alpha_2$ used in the definition of the function $f_2$, the same assertion gives $\alpha_1>0$
and $\alpha_2>0$.
By Theorem~\ref{th:Stupin:extremalTypeE3TwoTermsSystem}, the solution of the real system corresponding to the original
function $g$ determines it up to admissible rotations in the variable plane $z$.
Thus to find all possible functions it suffices to recover one representative for each of the four solutions.

By the definition of the numbers $\varphi^*$ and $\gamma^*$, the four possible pairs of values of the unknowns are
    \[
        (\varphi^*,\gamma^*),
        \quad
        (\varphi^*,2\pi-\gamma^*),
        \quad
        (\pi-\varphi^*,\pi-\gamma^*),
        \quad
        (\pi-\varphi^*,\pi+\gamma^*).
    \]
The corresponding values of $\mu$ can be chosen as follows:
    \[
        e^{i\frac{\gamma^*}3},
        \quad
        e^{-i\frac{\gamma^*}3},
        \quad
        e^{i\frac{\pi-\gamma^*}3},
        \quad
        e^{i\frac{\pi+\gamma^*}3}.
    \]
By formulas~\eqref{eq:Stupin:E3TwoTermsParameters}, the values $(\widetilde{\alpha}_1,\widetilde{\alpha}_2)$
corresponding to the four pairs are, respectively,
    \[
        (\alpha_1,\alpha_2),
        \quad
        (\alpha_2,\alpha_1),
        \quad
        (\alpha_1,\alpha_2),
        \quad
        (\alpha_2,\alpha_1).
    \]
Substituting these values and the equalities $z_1=\mu e^{i\varphi}$, $z_2=\mu e^{-i\varphi}$ into
representation~\eqref{eq:Stupin:E3TwoTerms_h} of the function $h$ gives the functions $f_2(z)$, $f_2^*(z)$,
$f_2^*(e^{4\pi i/3}z)$, and $f_2(e^{4\pi i/3}z)$, respectively.
Thus all possible functions have the form $f_2(e^{2\pi ik/3}z)$ and
    $f_2^*(e^{2\pi ik/3}z)$, $k=0,1,2$.

All four solutions of system~\eqref{eq:Stupin:E3TwoTermsRealSystem} have positive $\alpha_1,\alpha_2$ by
Assertion~\ref{st:Stupin:extremalTypeE3TwoTermsNonzeroCSystem}.
Hence, by Theorem~\ref{th:Stupin:extremalTypeE3TwoTermsSystem}, all six functions belong to $\mathbb E_2(3)$.
\end{proof}

\smallskip
Thus, up to admissible rotations in the variable plane $z$, we have found two functions of extremal type
\begin{multline}
    \label{eq:Stupin:E3TwoTermsExpansion}
    f_2(z)
    =
    e^{-t_2^*}
    +e^{-t_2^*}e^{i\frac{\gamma^*}3}
    \bigl(R^*e^{-i\gamma^*}-4c^*\bigr)z
    +\\+
    e^{-t_2^*}e^{2i\frac{\gamma^*}3}
    \bigl(e^{-2i\gamma^*}+2+4(c^*)^2
    -4R^*c^*e^{-i\gamma^*}\bigr)z^2
    +2e^{-t_2^*}\bigl(R^*(1+2(c^*)^2)-4c^*\cos\gamma^*\bigr)z^3
    +\\+
    \tfrac12e^{-t_2^*}e^{-2i\frac{\gamma^*}3}
    \bigl(e^{2i\gamma^*}+2+4(c^*)^2
    -4R^*c^*e^{i\gamma^*}\bigr)z^4
    +o(z^4)
    \approx
    \\
    \approx 0.036-(0.055+0.162i)z-(0.212-0.289i)z^2
    +\\+
    0.468z^3-(0.106+0.144i)z^4+o(z^4).
\end{multline}
The coefficients of the function $f_2^*$ are complex conjugates of those of $f_2$.
Expansion~\eqref{eq:Stupin:E3TwoTermsExpansion} shows that the coefficient of $z^4$ equals half the complex conjugate of
the coefficient of $z^2$.
This equality is preserved under complex conjugation of coefficients and the substitutions of $z$ by $e^{2\pi ik/3}z$,
$k=0,1,2$.
Hence all functions $f_2(e^{2\pi ik/3}z)$ and $f_2^*(e^{2\pi ik/3}z)$, $k=0,1,2$ satisfy
relation~\eqref{eq:Stupin:Marty}.

\section{The set $\mathbb E_2(3)$}

\begin{Assertion}
    \label{st:Stupin:extremalTypeE3TwoTerms}
The set $\mathbb E_2(3)$ consists of the functions described in Assertions~\ref{st:Stupin:extremalTypeE3TwoTermsZeroC}
and~\ref{st:Stupin:extremalTypeE3TwoTermsNonzeroC}.
In particular, $|\mathbb E_2(3)|=9$.
\end{Assertion}

\begin{proof}
Assertions~\ref{st:Stupin:extremalTypeE3TwoTermsZeroC} and~\ref{st:Stupin:extremalTypeE3TwoTermsNonzeroC} completely
describe $\mathbb E_2(3)$, since the cases $\varphi=\frac\pi2$ and $\varphi\in(0,\pi)\setminus\{\frac\pi2\}$ exhaust all
possibilities.
Expansion~\eqref{eq:Stupin:fSqrt6Over4Expansion} shows that the first coefficient of the function $f_{\frac{\sqrt6}4}$
equals $-\sqrt6e^{-2}\neq0$.
Expansion~\eqref{eq:Stupin:E3TwoTermsExpansion} shows that the first coefficient of the function $f_2$ is also nonzero,
since the imaginary part of the factor $R^*e^{-i\gamma^*}-4c^*$ equals
    $-R^*\sin\gamma^*\neq0$.
The first coefficient of the function $f_2^*$ is the complex conjugate of the first coefficient of $f_2$, so it is also
nonzero.
Thus each of these functions has three distinct admissible rotations in the variable plane $z$ corresponding to
    $k=0,1,2$.
For the function $f_{\frac{\sqrt6}4}$ the points are opposite.
For the functions $f_2$ and $f_2^*$ the points are not opposite.
The cubes of the products of the points of the functions $f_2$ and $f_2^*$ equal $e^{2i\gamma^*}$ and $e^{-2i\gamma^*}$,
respectively.
By Lemma~\ref{lem:Stupin:E3TwoTermsConstantsEstimates}, we have $\gamma^*\in(\pi/2,\pi)$.
Thus $\sin(2\gamma^*)\neq0$.
Hence the functions $f_2$ and $f_2^*$ are not related by admissible rotations in the variable plane
    $z$.
Hence the listed functions are pairwise distinct.
\end{proof}

\smallskip
Note that each of the two solutions of system~\eqref{eq:Stupin:E3TwoTermsRealSystem} found for $\varphi=\frac\pi2$
corresponds to the same class of functions related by admissible rotations in the variable plane $z$.
Conversely, every function in this class corresponds to both solutions of this system.

Similarly, each of the four solutions of system~\eqref{eq:Stupin:E3TwoTermsRealSystem} found for $\varphi\neq\frac\pi2$
determines one class of functions coinciding up to admissible rotations in the variable plane $z$.
The first and fourth solutions determine the class of the function $f_2$, and the second and third the class of $f_2^*$.
Conversely, every function in these classes corresponds to two solutions of this system.

These correspondences are explained by two distinct symmetries.
The first symmetry concerns the choice of a root of the equation $\mu^3=e^{i\gamma}$.
The transformations
\[
    \mu\longmapsto\mu e^{\frac{2\pi ki}3},
    \quad
    k=0,1,2,
\]
form a cyclic group of order three $\mathbb Z_3$.
They leave the solution of the system unchanged and replace the function $f(z)$ by
$f(e^{\frac{2\pi ki}3}z)$.
Thus one solution of the system determines one class of three functions under the action of the group $\mathbb Z_3$.

The second symmetry concerns interchanging the two terms in the representation of the function $h$ in
formula~\eqref{eq:Stupin:E3TwoTerms_h}, that is, interchanging the pairs $(\alpha_1,z_1)$ and $(\alpha_2,z_2)$.
In the notation used, it corresponds to the transformation of solutions
\[
    (t,R,\varphi,\gamma,\alpha_1,\alpha_2)
    \longmapsto
    (t,R,\pi-\varphi,\gamma',\alpha_2,\alpha_1),
\]
where $\gamma'\in[0,2\pi)$ and $\gamma'\equiv\gamma+\pi\pmod{2\pi}$.
Here $\mu$ is replaced by $-\mu$, interchanging the points $z_1$ and $z_2$.
This transformation changes the solution of the system but leaves the functions $h$ and $f=e^{-h}$ unchanged, since it
only interchanges the two terms in representation~\eqref{eq:Stupin:E3TwoTerms_h} of $h$.
Applying it twice returns the original solution.
Thus the identity transformation and this transformation form a cyclic group of order two
$\mathbb Z_2$.

The actions of the groups $\mathbb Z_3$ and $\mathbb Z_2$ commute.
Together they can be regarded as an action of the group $\mathbb Z_3\times\mathbb Z_2$ on representations of functions.

\section{The set $\mathbb E_3(3)$}

\begin{Assertion}
    \label{st:Stupin:extremalTypeE3ThreeTermsRelations}
The function $f\in\mathbb E_3(3)$ if and only if $f=e^{-h}$, where
    \begin{equation}
        \label{eq:Stupin:E3ThreeTerms_h}
        h(z)=\sum\limits_{j=1}^3\alpha_j\frac{1+z_j z}{1-z_j z},
        \quad
        \alpha_j>0,
        \quad
        z_j=e^{i\varphi_j},
        \quad
        j=1,2,3,
        \quad
        0\leqslant\varphi_1<\varphi_2<\varphi_3<2\pi,
    \end{equation}
and there exists a number $t>0$ such that for $s\coloneqq z_1+z_2+z_3$ the vector
$(t,\alpha_1,\alpha_2,\alpha_3,z_1,z_2,z_3)$ satisfies the system of equations
    \begin{equation}
        \label{eq:Stupin:E3ThreeTermsRelations}
        \begin{cases}
            t=\alpha_1+\alpha_2+\alpha_3,\\
            \sum\limits_{j=1}^3\alpha_j z_j=s,\\
            \sum\limits_{j=1}^3\alpha_j z_j^2=\frac12s^2+\overline{s},\\
            \sum\limits_{j=1}^3\alpha_j z_j^3=-1+|s|^2+\frac13s^3,\\
            z_1 z_2 z_3=-1.
        \end{cases}
    \end{equation}
\end{Assertion}

\begin{proof}
\textbf{Necessity.}
    \textbf{1.}
We find a normalized polynomial $\Psi$ external for $H$.
Let $f\in\mathbb E_3(3)$.
By the definition of $\mathbb E_3(3)$, the function $h\coloneqq-\ln f$ has
representation~\eqref{eq:Stupin:E3ThreeTerms_h}.
Set $t\coloneqq\{h\}_0>0$.
Then $\{f\}_0=e^{-t}$.
By Corollary~\ref{cor:Stupin:f_n_eq_f_0_norm_Psi2}, there exists a normalized polynomial $\Psi$ external for $H$.
By Corollary~\ref{cor:Stupin:extremalTypeZerosOfH}, it has the form
    $\Psi(z)=\prod\limits_{j=1}^3(z-z_j)$.
For the external polynomial $P=p_3\Psi$, formula~\eqref{eq:Stupin:f_n_p_n} for $k=3$ gives
$e^{-t}=p_3\overline{p_0}=|p_3|^2\overline{\Psi(0)}$.
Thus $\Psi(0)>0$.
Since $|\Psi(0)|=|z_1 z_2 z_3|=1$, we have $\Psi(0)=1$ and $z_1 z_2 z_3=-1$.
Also, $|p_3|^2=e^{-t}$.
The condition $|z_j|=1$ implies that
    \[
        z_1 z_2+z_1 z_3+z_2 z_3
        =z_1 z_2 z_3(\overline{z_1}+\overline{z_2}+\overline{z_3})
        =-\overline{s}.
    \]
Hence,
    \begin{equation}
        \label{eq:Stupin:E3ThreeTermsPsiProof}
        \Psi(z)=1-\overline{s}z-sz^2+z^3.
    \end{equation}

    \textbf{2.}
We find the coefficients of the function $f$.
By formula~\eqref{eq:Stupin:E3ThreeTermsPsiProof}, the coefficients of the polynomial $P=p_3\Psi$ satisfy the following
relations:
    \[
        p_0=p_3,
        \quad
        p_1=-p_3\overline{s},
        \quad
        p_2=-p_3s.
    \]
Substituting them into formulas~\eqref{eq:Stupin:f_n_p_n} and using $|p_3|^2=e^{-t}$ gives
    \begin{equation}
        \label{eq:Stupin:E3ThreeTermsFCoefficients}
        \begin{aligned}
            \{f\}_0&=p_3\overline{p_0}=e^{-t},\\
            \{f\}_1&=p_2\overline{p_0}+p_3\overline{p_1}=-2se^{-t},\\
            \{f\}_2&=p_1\overline{p_0}+p_2\overline{p_1}+p_3\overline{p_2}
            =(s^2-2\overline{s})e^{-t},\\
            \{f\}_3&=|p_0|^2+|p_1|^2+|p_2|^2+|p_3|^2=2(1+|s|^2)e^{-t}.
        \end{aligned}
    \end{equation}

    \textbf{3.}
We derive the system of equations.
Using formula~\eqref{eq:Stupin:deForBCoeff} for $\ell=1,2,3$, express the coefficients of the function $h$ and
substitute expressions~\eqref{eq:Stupin:E3ThreeTermsFCoefficients} for the coefficients of $f$:
    \begin{equation}
        \label{eq:Stupin:E3ThreeTermsHCoefficients}
        \begin{aligned}
            \{h\}_0&=-\ln\{f\}_0=t,\\
            \{h\}_1&=-\tfrac{\{f\}_1}{\{f\}_0}=2s,\\
            \{h\}_2&=\tfrac{\{f\}_1^2}{2\{f\}_0^2}-\tfrac{\{f\}_2}{\{f\}_0}
            =
            2s^2-(s^2-2\overline{s})=s^2+2\overline{s},\\
            \{h\}_3&=\tfrac{\{f\}_1\{f\}_2}{\{f\}_0^2}
            -\tfrac{\{f\}_1^3}{3\{f\}_0^3}-\tfrac{\{f\}_3}{\{f\}_0}
            =
            -2s(s^2-2\overline{s})+\tfrac83s^3-2(1+|s|^2)
            =
            -2+2|s|^2+\tfrac23s^3.
        \end{aligned}
    \end{equation}
On the other hand, representation~\eqref{eq:Stupin:E3ThreeTerms_h} implies that
    \[
        \{h\}_0=\alpha_1+\alpha_2+\alpha_3,
        \quad
        \{h\}_k=2\sum\limits_{j=1}^3\alpha_jz_j^k,
        \quad
        k=1,2,3.
    \]
Equating the resulting expressions for $\{h\}_k$, $k=0,1,2,3$ gives the first four equations of
system~\eqref{eq:Stupin:E3ThreeTermsRelations}.
The fifth equation was derived in the first step of the proof when finding the polynomial $\Psi$.

\textbf{Sufficiency.}
    \textbf{1.}
We construct the function corresponding to a solution of the system and find its coefficients.
Let $(t,\alpha_1,\alpha_2,\alpha_3,z_1,z_2,z_3)$ be a solution of system~\eqref{eq:Stupin:E3ThreeTermsRelations} with
the constraints in the statement.
Define $h$ by formula~\eqref{eq:Stupin:E3ThreeTerms_h} and set $f=e^{-h}$.
Since $|z_j|=1$, $j=1,2,3$, the functions $h$ and $f$ are holomorphic in $\Delta$.
Since $\alpha_j>0$, we have $h\in C$.
Hence $f\in B$.
The first equation of the system gives $\{h\}_0=t$.
Substituting the second, third, and fourth equations of the system into the formulas
$\{h\}_k=2\sum\limits_{j=1}^3\alpha_jz_j^k$, $k=1,2,3$ gives the values of the coefficients $h$ on the right-hand sides
of equalities~\eqref{eq:Stupin:E3ThreeTermsHCoefficients}.
Since $f=e^{-h}$, formula~\eqref{eq:Stupin:deForBCoeff} gives
    \[
        \begin{aligned}
            \{f\}_0&=e^{-t},\\
            \{f\}_1&=-e^{-t}\{h\}_1,\\
            \{f\}_2&=e^{-t}\big(\tfrac12\{h\}_1^2-\{h\}_2\big),\\
            \{f\}_3&=e^{-t}\big(-\tfrac16\{h\}_1^3+\{h\}_1\{h\}_2-\{h\}_3\big).
        \end{aligned}
    \]
Substituting the coefficients $h$ found into these equalities gives the expressions for the coefficients $f$ on the
right-hand sides of formulas~\eqref{eq:Stupin:E3ThreeTermsFCoefficients}.

    \textbf{2.}
We check that the constructed function $f$ belongs to $\mathbb E_3(3)$.
Set $H(z)=\{f\}_3+2\{f\}_2z+2\{f\}_1z^2+2\{f\}_0z^3$.
The fifth equation of system~\eqref{eq:Stupin:E3ThreeTermsRelations} and the conditions $|z_j|=1$ give
    $\prod\limits_{j=1}^3(z-z_j)=1-\overline{s}z-sz^2+z^3$.
Now denote this polynomial by $\Psi$ and set $P(z)\coloneqq\sqrt{e^{-t}}\Psi(z)=\sum\limits_{j=0}^3 p_j z^j$.
Construct the polynomial
        $\widetilde H(z)
        \coloneqq
        \sum\limits_{j=0}^3|p_j|^2
        +2\sum\limits_{k=1}^3\left(\sum\limits_{j=0}^{3-k}p_{j+k}\overline{p_j}\right)z^k$.
Since $p_3=\sqrt{e^{-t}}>0$, the constant coefficient $\widetilde H$ is positive, so $\widetilde H\in C^3$ by
Theorem~\ref{th:Stupin:characterization_of_polynomials_from_C}.

We have $p_0=p_3=\sqrt{e^{-t}}$, $p_1=-\sqrt{e^{-t}}\overline{s}$, $p_2=-\sqrt{e^{-t}}s$.
Thus the sums calculated in~\eqref{eq:Stupin:E3ThreeTermsFCoefficients} give
    \[
        \begin{aligned}
            \sum\limits_{j=0}^3|p_j|^2&=2e^{-t}(1+|s|^2)=\{f\}_3,\\
            \sum\limits_{j=0}^2p_{j+1}\overline{p_j}&=e^{-t}(s^2-2\overline{s})=\{f\}_2,\\
            \sum\limits_{j=0}^1p_{j+2}\overline{p_j}&=-2e^{-t}s=\{f\}_1,\\
            p_3\overline{p_0}&=e^{-t}=\{f\}_0.
        \end{aligned}
    \]
Hence all coefficients of $H$ and $\widetilde H$ coincide.
Thus $H=\widetilde H$ and $H\in C^3$, with $\{f\}_3=\sum\limits_{j=0}^3|p_j|^2>0$.

By Assertion~\ref{st:Stupin:polynomialModulusIdentity}, we have
    $\operatorname{Re}H(e^{i\theta})=|P(e^{i\theta})|^2$, $\theta\in\mathbb R$.
By the definition of the polynomials $P$ and $\Psi$, we have $P(z_j)=0$, $j=1,2,3$.
Hence $\operatorname{Re}H(z_j)=0$, $j=1,2,3$.
The points $z_1,z_2,z_3$ are distinct, and the parameters $\alpha_1,\alpha_2,\alpha_3$ in
representation~\eqref{eq:Stupin:E3ThreeTerms_h} of the function $h$ are positive.
Hence $f\in\mathbb E_3(3)$.
\end{proof}

\begin{Corollary}
    \label{cor:Stupin:extremalTypeE3ThreeTermsRelations}
Let $f\in\mathbb E_3(3)$, and let $s$ be defined in Assertion~\ref{st:Stupin:extremalTypeE3ThreeTermsRelations}.
The normalized polynomial external for $H$ has the form
    \begin{equation}
        \label{eq:Stupin:E3ThreeTermsPsi}
        \Psi(z)=1-\overline{s}z-sz^2+z^3.
    \end{equation}
\end{Corollary}

\begin{proof}
Let $f\in\mathbb E_3(3)$.
Then, by Assertion~\ref{st:Stupin:extremalTypeE3ThreeTermsRelations}, the vector
$(t,\alpha_1,\alpha_2,\alpha_3,z_1,z_2,z_3)$ determining the function $f$ is a solution of
system~\eqref{eq:Stupin:E3ThreeTermsRelations}.
Formula~\eqref{eq:Stupin:E3ThreeTermsPsi} was obtained in the proof of
Assertion~\ref{st:Stupin:extremalTypeE3ThreeTermsRelations}.
(See~\eqref{eq:Stupin:E3ThreeTermsPsiProof}.)
\end{proof}

\begin{Corollary}
    \label{cor:Stupin:E3ThreeTermsT}
Let $f\in\mathbb E_3(3)$, and let $s$ and $t$ be defined in Assertion~\ref{st:Stupin:extremalTypeE3ThreeTermsRelations}.
Then $s$ and $t$ satisfy the system of equations
    \begin{equation}
        \label{eq:Stupin:E3ThreeTermsT}
        \begin{cases}
            t=1+|s|^2+\tfrac16s^3,\\
            t=2-\tfrac12\tfrac{s^3}{|s|^2}, & s\neq0.
        \end{cases}
    \end{equation}
\end{Corollary}

\begin{proof}
Let $f\in\mathbb E_3(3)$.
Then, by Assertion~\ref{st:Stupin:extremalTypeE3ThreeTermsRelations}, the vector
$(t,\alpha_1,\alpha_2,\alpha_3,z_1,z_2,z_3)$ determining the function $f$ is a solution of
system~\eqref{eq:Stupin:E3ThreeTermsRelations}.
By Corollary~\ref{cor:Stupin:extremalTypeE3ThreeTermsRelations}, the polynomial $\Psi$ has the
form~\eqref{eq:Stupin:E3ThreeTermsPsi}.

Dividing $\Psi(z_j)=0$ by $z_j\neq0$ gives $z_j^{-1}=-z_j^2+sz_j+\overline{s}$.
Multiplying this equality by $\alpha_j$ and summing over $j=1,2,3$ gives
    \begin{equation}
        \label{eq:Stupin:E3ThreeTermsInverseSum}
        \sum\limits_{j=1}^3\alpha_j z_j^{-1}
        =
        -\sum\limits_{j=1}^3\alpha_j z_j^2
        +s\sum\limits_{j=1}^3\alpha_j z_j
        +\overline{s}\sum\limits_{j=1}^3\alpha_j.
    \end{equation}
Since $|z_j|=1$ and $\alpha_j>0$, using the second equation of system~\eqref{eq:Stupin:E3ThreeTermsRelations}, we
transform the left-hand side of~\eqref{eq:Stupin:E3ThreeTermsInverseSum} as follows:
        $\sum\limits_{j=1}^3\alpha_j z_j^{-1}
        =
        \sum\limits_{j=1}^3\alpha_j\overline{z_j}
        =
        \overline{\sum\limits_{j=1}^3\alpha_j z_j}
        =
        \overline{s}$.
Substituting the third, second, and first equations of system~\eqref{eq:Stupin:E3ThreeTermsRelations} into the
right-hand side of equation~\eqref{eq:Stupin:E3ThreeTermsInverseSum} transforms
equation~\eqref{eq:Stupin:E3ThreeTermsInverseSum} into $\overline{s}=-\frac12s^2-\overline{s}+s^2+\overline{s}t$.
Thus $\overline{s}(t-2)=-\frac{s^2}2$.
For $s\neq0$, this implies $t=2-\frac{s^2}{2\overline{s}}=2-\frac{s^3}{2|s|^2}$.

The equation $\Psi(z_j)=0$ gives $z_j^3=sz_j^2+\overline{s}z_j-1$.
Multiplying by $\alpha_j$ and summing over $j=1,2,3$ gives
        $-1+|s|^2+\frac13s^3
        =
        s(\frac12s^2+\overline{s})+|s|^2-t$.
Hence $t=1+|s|^2+\frac16s^3$.
\end{proof}

\begin{Assertion}
    \label{st:Stupin:extremalTypeE3ThreeTermsZeroS}
Let $f=e^{-h}\in\mathbb E_3(3)$ and $s=0$.
Then $\alpha_1=\alpha_2=\alpha_3=\frac13$, and the points $z_1$, $z_2$, $z_3$ in
representation~\eqref{eq:Stupin:E3ThreeTerms_h} of the function $h$ are roots of the equation $z^3=-1$.
That is, $f(z)=F(-z^3,1)$.
\end{Assertion}

\begin{proof}
Let $f=e^{-h}\in\mathbb E_3(3)$ and $s=0$.
Then $\Psi(z)=1+z^3$ by Corollary~\ref{cor:Stupin:extremalTypeE3ThreeTermsRelations}.
It follows that the points $z_j$ are roots of the equation $z^3=-1$.

The second and third equations of system~\eqref{eq:Stupin:E3ThreeTermsRelations} imply $\sum\limits_{j=1}^3\alpha_j
z_j=0$ and $\sum\limits_{j=1}^3\alpha_j z_j^2=0$, respectively.
Since $z_1=e^{i\frac\pi3}$, $z_2=-1$, $z_3=e^{i\frac{5\pi}3}$, we have $\alpha_1=\alpha_2=\alpha_3$.
The first equation of system~\eqref{eq:Stupin:E3ThreeTermsT} gives $t=1$.
The first equation of system~\eqref{eq:Stupin:E3ThreeTermsRelations} gives $\alpha_1=\alpha_2=\alpha_3=\frac13$.
Hence $f(z)=F(-z^3,1)$.
\end{proof}

\begin{Lemma}
    \label{lem:Stupin:extremalTypeE3ThreeTermsCandidates}
Let $f=e^{-h}\in\mathbb E_3(3)$ and $s\neq0$.
Then, after a rotation in the variable plane $z$ through a multiple of $\frac{2\pi}3$, the points $z_j$ in
representation~\eqref{eq:Stupin:E3ThreeTerms_h} of the function $h$ have the form
    \begin{equation}
        \label{eq:Stupin:E3ThreeTermsPoints}
        z_1=e^{i\varphi_c},
        \quad
        z_2=-1,
        \quad
        z_3=e^{-i\varphi_c},
        \quad
        \varphi_c=\arccos c\in(0,\pi),
    \end{equation}
where $c$ is a root of the equation
    \begin{equation}
        \label{eq:Stupin:E3ThreeTermsCubic}
        2c^3+3c^2-3c-1=0,
    \end{equation}
lying in $(0,1)$ or in $(-1/2,0)$.
Moreover, $s=2c-1$.
\end{Lemma}

\begin{proof}
Since $t>0$, any equation of system~\eqref{eq:Stupin:E3ThreeTermsT} implies that $s^3$ is real.
Thus, after a rotation in the variable plane $z$ through a multiple of $\frac{2\pi}3$, we may assume that $s$ is real.
In this case, formula~\eqref{eq:Stupin:E3ThreeTermsPsi} shows that the polynomial $\Psi$ has real coefficients.
Hence its set of roots is invariant under complex conjugation.
Since the points $z_1,z_2,z_3$ are pairwise distinct,
$\{z_1,z_2,z_3,\}\subset\{e^{i\varphi},e^{i(2\pi-\varphi)},-1,1\}$, where $\varphi\in(0,\pi)$.
This and the fifth equation of system~\eqref{eq:Stupin:E3ThreeTermsRelations} imply that one root equals $-1$.
Number them as follows: $z_1=e^{i\varphi}$, $z_2=-1$, $z_3=e^{-i\varphi}$.

Set $c\coloneqq\cos\varphi$.
Then $s=z_1+z_2+z_3=2c-1$.
Equating the right-hand sides of the equations of system~\eqref{eq:Stupin:E3ThreeTermsT} gives $s^3+6s^2+3s-6=0$.
Substituting $s=2c-1$ gives $2c^3+3c^2-3c-1=0$.
Hence $c$ satisfies equation~\eqref{eq:Stupin:E3ThreeTermsCubic}.
By Cardano's formula, this equation has three real roots, only two of which belong to the segment
    $[-1,1]$.
Since the values of the left-hand side of equation~\eqref{eq:Stupin:E3ThreeTermsCubic} at $-1/2$, $0$, and $1$ equal
$1$, $-1$, and $1$, respectively, these roots lie in the intervals $(-1/2,0)$ and $(0,1)$.
Since $0<\varphi<\pi$, we have $\varphi=\varphi_c\coloneqq\arccos c$.
\end{proof}

\begin{Lemma}
    \label{lem:Stupin:E3ThreeTermsParameters}
Let $f=e^{-h}\in\mathbb E_3(3)$, $s\neq0$, and let the points $z_j$ in representation~\eqref{eq:Stupin:E3ThreeTerms_h}
of the function $h$ have the form~\eqref{eq:Stupin:E3ThreeTermsPoints}, where $c$ satisfies the conditions of
Lemma~\ref{lem:Stupin:extremalTypeE3ThreeTermsCandidates}.
Then
    \[
        \alpha_1=\alpha_3=\tfrac14\tfrac{2c+3}{c+1}>0,
        \quad
        \alpha_2=1+\tfrac{c}2\tfrac{2c+3}{c+1}-2c>0.
    \]
Conversely, for every such $c$ the vector $(t,\alpha_1,\alpha_2,\alpha_3,z_1,z_2,z_3)$, where $t=\frac52-c$, the numbers
$\alpha_j$ are given above, and $z_j$ by formula~\eqref{eq:Stupin:E3ThreeTermsPoints}, is a solution of
system~\eqref{eq:Stupin:E3ThreeTermsRelations} satisfying the constraints of
Assertion~\ref{st:Stupin:extremalTypeE3ThreeTermsRelations}.
\end{Lemma}

\begin{proof}
    \textbf{1.}
We find the parameters $\alpha_j$, $j=1,2,3$ in representation~\eqref{eq:Stupin:E3ThreeTerms_h} of the function $h$.
Set $\varphi\coloneqq\varphi_c=\arccos c$.
By Lemma~\ref{lem:Stupin:extremalTypeE3ThreeTermsCandidates}, we have $s=2c-1$ and $0<\varphi<\pi$.
The second equation of system~\eqref{eq:Stupin:E3ThreeTermsRelations} gives
    \begin{equation}
        \label{eq:Stupin:E3ThreeTermsFirstMoment}
        \alpha_1 e^{i\varphi}-\alpha_2+\alpha_3 e^{-i\varphi}=2c-1.
    \end{equation}

Take the imaginary part of equation~\eqref{eq:Stupin:E3ThreeTermsFirstMoment}.
Since $\sin\varphi>0$ and $\operatorname{Im}(\alpha_1 e^{i\varphi}+\alpha_3 e^{-i\varphi})=0$, $\alpha_1=\alpha_3$.

Denote the common value of $\alpha_1$ and $\alpha_3$ by $\alpha$.
Taking the real part of equation~\eqref{eq:Stupin:E3ThreeTermsFirstMoment} gives $2c\alpha-\alpha_2=2c-1$.
Substituting the first equation of system~\eqref{eq:Stupin:E3ThreeTermsRelations} and $s=2c-1$ into the second equation
of system~\eqref{eq:Stupin:E3ThreeTermsT} gives $2\alpha+\alpha_2=2-\frac{2c-1}2$.

Solving these two equations for the unknowns $\alpha$ and $\alpha_2$ gives
    \[
        \alpha_1=\alpha_3=\tfrac14\tfrac{2c+3}{c+1}>0,
        \quad
        \alpha_2=1+\tfrac{c}2\tfrac{2c+3}{c+1}-2c>0.
    \]
Since $-1/2<c<1$, we have $c+1>0$ and $2c+3>0$, whence $\alpha_1=\alpha_3>0$.
Furthermore, $\alpha_2=\frac{-2c^2+c+2}{2(c+1)}$ and $-2c^2+c+2>0$ on $[-1/2,1]$, so $\alpha_2>0$.

    \textbf{2.}
We show equivalence of the transformations used to solve system~\eqref{eq:Stupin:E3ThreeTermsRelations}.
For $0<\varphi<\pi$ and $c+1>0$, the linear equations above are equivalent to the resulting formulas for $\alpha_j$, and
the first equation of system~\eqref{eq:Stupin:E3ThreeTermsRelations} gives $t=\frac52-c>0$.
For the values $\alpha_j$ found, the third equation is satisfied identically, since
    \[
        \sum\limits_{j=1}^3\alpha_jz_j^2
        =2\alpha(2c^2-1)+\alpha_2=2c^2-\tfrac12=\tfrac12s^2+s.
    \]
The fourth equation is equivalent to equation~\eqref{eq:Stupin:E3ThreeTermsCubic}, since
    \[
        \sum\limits_{j=1}^3\alpha_jz_j^3-\big(-1+s^2+\tfrac13s^3\big)
        =2\alpha(4c^3-3c)-\alpha_2+1-s^2-\tfrac13s^3
        =\tfrac23(2c^3+3c^2-3c-1).
    \]
The fifth equation is satisfied by the choice of $z_j$ in~\eqref{eq:Stupin:E3ThreeTermsPoints}.
Thus, for these $z_j$, the original system is equivalent to the formulas found for $t$, $\alpha_j$ and
equation~\eqref{eq:Stupin:E3ThreeTermsCubic}.
Positivity of $\alpha_j$ has already been established, and the arguments of $z_j$ equal $\varphi$, $\pi$, $2\pi-\varphi$
and are in the required order.
\end{proof}

Denote the roots of equation~\eqref{eq:Stupin:E3ThreeTermsCubic} lying on $[-1,1]$ by $c_1^*$ and $c_2^*$.
Let $c_1^*\in(0,1)$ and $c_2^*\in(-1/2,0)$.
Calculations give $c_1^*\approx0.873$, $c_2^*\approx-0.272$.

\begin{Assertion}
    \label{st:Stupin:extremalTypeE3ThreeTerms}
The set $\mathbb E_3(3)$ has the following form\textnormal:
    \[
        \mathbb E_3(3)
        =
        \big\{F(-z^3,1)\big\}
        \cup
        \big\{
            f_{c_j^*}(\zeta^kz)
            :
            j=1,2,\quad
            k=0,1,2
        \big\},
    \]
where $f_c(z)
    \coloneqq
    e^{
            -\alpha_c\frac{1+e^{i\varphi_c}z}{1-e^{i\varphi_c}z}
            -\beta_c\frac{1-z}{1+z}
            -\alpha_c\frac{1+e^{-i\varphi_c}z}{1-e^{-i\varphi_c}z}
}$, $c\in\{c_1^*,c_2^*\}$, $c_1^*$ and $c_2^*$ are the roots of equation~\eqref{eq:Stupin:E3ThreeTermsCubic},
    $\zeta\coloneqq e^{2\pi i/3}$, $\varphi_c\coloneqq\arccos c$, $\alpha_c\coloneqq\tfrac14\tfrac{2c+3}{c+1}$,
    $\beta_c\coloneqq1+2c(\alpha_c-1)$.
In particular, $|\mathbb E_3(3)|=7$.
\end{Assertion}

\begin{proof}
By Assertion~\ref{st:Stupin:extremalTypeE3ThreeTermsZeroS} and
Lemmas~\ref{lem:Stupin:extremalTypeE3ThreeTermsCandidates},~\ref{lem:Stupin:E3ThreeTermsParameters}, every function in
$\mathbb E_3(3)$, after a rotation in the variable plane $z$ through a multiple of $\frac{2\pi}3$, coincides with one of
the functions $F(-z^3,1)$, $f_{c_1^*}(z)$, or $f_{c_2^*}(z)$.
Conversely, the function $f(z)=F(-z^3,1)$ belongs to $\mathbb E(3)$ by Theorem~\ref{th:Stupin:converseToUniqunessInB}.
For it, $h=-\ln f$ is given by formula~\eqref{eq:Stupin:E3ThreeTerms_h} for $\alpha_1=\alpha_2=\alpha_3=\frac13$ and
$z_j=e^{i(2j-1)\pi/3}$, $j=1,2,3$, so $f\in\mathbb E_3(3)$.
By Lemma~\ref{lem:Stupin:E3ThreeTermsParameters}, for every function $f_{c_j^*}$, $j=1,2$, the corresponding vector
$(t,\alpha_1,\alpha_2,\alpha_3,z_1,z_2,z_3)$ is a solution of system~\eqref{eq:Stupin:E3ThreeTermsRelations} satisfying
the constraints of Assertion~\ref{st:Stupin:extremalTypeE3ThreeTermsRelations}.
Thus, by Assertion~\ref{st:Stupin:extremalTypeE3ThreeTermsRelations}, these functions belong to $\mathbb E_3(3)$.
Rotations in the variable plane $z$ through multiples of $\frac{2\pi}3$ preserve membership in $\mathbb E_3(3)$.

Since $2c_j^*-1\neq0$, the first coefficients of the functions $f_{c_j^*}$, $j=1,2$ are nonzero.
Thus each of these functions gives three distinct functions when $z$ is replaced by $\zeta^k z$, $k=0,1,2$, and none
coincides with $F(-z^3,1)$, whose first coefficient vanishes.
Also, $f_{c_1^*}(0)\neq f_{c_2^*}(0)$, since $c_1^*\neq c_2^*$.
Hence the listed functions are pairwise distinct.
\end{proof}

\smallskip
Thus, up to rotations in the variable plane $z$, we have found three functions of extremal type.
\label{page:Stupin:E3ThreeTermsExpansions}
For the first of them $F(-z^3,1)=e^{-1}+2e^{-1}z^3+o(z^4)$.
For the other two functions, set $s_j\coloneqq2c_j^*-1$ and $t_j\coloneqq\tfrac52-c_j^*$, $j=1,2$.
Then
\begin{multline*}
    f_{c_j^*}(z)
    =
    e^{-t_j}-2s_je^{-t_j}z+(s_j^2-2s_j)e^{-t_j}z^2+
    \\
    +2(1+s_j^2)e^{-t_j}z^3+\tfrac12(s_j^2-2s_j)e^{-t_j}z^4
    +o(z^4),
    \quad
    j=1,2.
\end{multline*}
For the function $F(-z^3,1)$ both sides of relation~\eqref{eq:Stupin:Marty} vanish, and for the other two functions
    $4\{f_{c_j^*}\}_4
    =
    2(s_j^2-2s_j)e^{-t_j}
    =
    2\overline{\{f_{c_j^*}\}_2}$,
    $j=1,2$.
These equalities are preserved when $z$ is replaced by $\zeta^k z$, $k=0,1,2$.
Hence all seven functions satisfy relation~\eqref{eq:Stupin:Marty}.

\section{The set $\mathbb E(3)$}

\begin{Theorem}
    \label{th:Stupin:extremalTypeE3}
The set $\mathbb E(3)$ has the form
    $\mathbb E(3)=\mathbb E_1(3)\mathbin{\dot\cup}\mathbb E_2(3)\mathbin{\dot\cup}\mathbb E_3(3)$
and consists of $19$ functions.
\end{Theorem}

\begin{proof}
The required equality follows directly from the definition of the sets $\mathbb E_1(3)$, $\mathbb E_2(3)$, and $\mathbb
E_3(3)$.
The number of functions follows from
Assertions~\ref{lem:Stupin:extremalTypeE3OneTerm},~\ref{st:Stupin:extremalTypeE3TwoTerms}
and~\ref{st:Stupin:extremalTypeE3ThreeTerms}, which describe these sets.
\end{proof}

Note that comparing the description obtained with the sharp estimations on the classes $B_t$ and $B_t^r$
from~\cite{Prokhorov:ref1, Stupin:ref15} shows that the functions in the sets $\mathbb E_1(3)$, $\mathbb E_2(3)$, and
the function $F(-z^3,1)$ attain a global maximum of $|\{f\}_3|$ in the classes $B_t$ containing these functions.

At the same time, the functions $f_{c_j^*}(\zeta^k z)$, $j=1,2$, $k=0,1,2$ do not attain a global maximum in the classes
$B_{t_j}$, where $t_j=\tfrac52-c_j^*$.

For example, for $t_j=\tfrac52-c_j^*$ and $d_j=\sqrt{t_j-\tfrac12}$ the functions
    $g_j(z)\coloneqq e^{
        -\frac{t_j+d_j}{2}\frac{1+z}{1-z}
        -\frac{t_j-d_j}{2}\frac{1-z}{1+z}
}$, $j=1,2$, belong to $B_{t_j}$, respectively, with
$\{g_1\}_3\approx0.627>\{f_{c_1^*}\}_3\approx0.612$ and $\{g_2\}_3\approx0.571>\{f_{c_2^*}\}_3\approx0.423$.

Thus, unlike the set $\mathbb E(2)$, not every function in $\mathbb E(3)$ attains a global maximum for a fixed value of
$f(0)$.

\begin{Theorem}
The Krzyz conjecture holds for $n=3$.
\end{Theorem}

\begin{proof}
Let $f$ be a normalized global extremal in the Krzyz problem for the coefficient with index $n=3$.
Then by Theorem~\ref{th:Stupin:localExtremalIsExtremalType}, the function $f\in\mathbb E(3)$.
By Theorem~\ref{th:Stupin:extremalTypeE3},
        $\mathbb E(3)
        =
        \mathbb E_1(3)\mathbin{\dot\cup}\mathbb E_2(3)\mathbin{\dot\cup}\mathbb E_3(3)$.
Use the expansions given on p.~\pageref{page:Stupin:E3OneTermExpansion}, in
formulas~\eqref{eq:Stupin:fSqrt6Over4Expansion}, \eqref{eq:Stupin:E3TwoTermsExpansion}, and on
p.~\pageref{page:Stupin:E3ThreeTermsExpansions}.
For functions in $\mathbb E_1(3)$ we have
    \[
        \{F(e^{i(2k+1)\frac{\pi}3}z,t_1^*)\}_3
        =
        \tfrac23t_1^*\bigl(2(t_1^*)^2-6t_1^*+3\bigr)e^{-t_1^*}
        \approx
        0.552,
        \quad
        k=0,1,2.
    \]
For functions in $\mathbb E_2(3)$ we have
    \begin{gather*}
        \{f_{\frac{\sqrt6}4}(e^{i\frac{2\pi k}3}z)\}_3
        =
        2\sqrt6e^{-2}
        \approx
        0.663,
        \quad
        k=0,1,2,
        \\
        \{f_2(e^{2\pi ik/3}z)\}_3
        =
        \{f_2^*(e^{2\pi ik/3}z)\}_3
        =
        2R^*e^{-t_2^*}\bigl(1+(4t_2^*-10)(c^*)^2\bigr)
        \approx
        0.468,
        \quad
        k=0,1,2.
    \end{gather*}
For functions in $\mathbb E_3(3)$ we have
    \begin{gather*}
        \{F(-z^3,1)\}_3
        =
        \frac2e
        \approx
        0.736,
        \\
        \{f_{c_1^*}(\zeta^kz)\}_3
        =
        2\bigl(1+(2c_1^*-1)^2\bigr)e^{-(5/2-c_1^*)}
        \approx
        0.612,
        \quad
        k=0,1,2,
        \\
        \{f_{c_2^*}(\zeta^kz)\}_3
        =
        2\bigl(1+(2c_2^*-1)^2\bigr)e^{-(5/2-c_2^*)}
        \approx
        0.423,
        \quad
        k=0,1,2.
    \end{gather*}

The largest of all the listed values equals $\frac2e$ and is attained only by the function $F(-z^3,1)$.
Thus $\{f\}_3=\frac2e$ and $f(z)=F(-z^3,1)$.
Now let $g\in B$ and $|\{g\}_3|=\frac2e$.
Then $g$ is a global extremal.
Rotations in the variable planes $z$ and $w$ can make its coefficients with indices $0$ and $3$ positive.
The resulting function is a local extremal and, by Theorem~\ref{th:Stupin:localExtremalIsExtremalType}, belongs to
$\mathbb E(3)$.
Hence it coincides with $F(-z^3,1)$.
Returning to the function $g$, we conclude that there exist numbers $\theta,\varphi\in[0,2\pi)$ for which
$g(z)=e^{i\theta}F(e^{i\varphi}z^3,1)$.
This proves the assertion.
\end{proof}

\section*{Supplement}

In Section~\ref{sec:Stupin:extremalTypeE1} and thereafter, studying the Krzyz problem for $n=1,2,3$ was reduced, in
particular, to solving algebraic equations with real coefficients and comparing the coefficients $\{f\}_n$ of the
functions $f$ found.
In the cases $n=1,2,3$ considered above, elimination of unknowns for each $m\in\{1,\ldots,n\}$ gives finitely many
solutions of the corresponding system.
This approach can also be used for other indices $n$.
It does not require an explicit description of $\mathfrak E(n)$ or $\mathbb E(n)$, as given above for $\mathbb E(n)$.
It suffices that $\mathfrak E(n)$ be contained in the set of functions determined by solutions of
systems~\eqref{eq:Stupin:conclusionSystem} for $m=1,\ldots,n$, and that all these functions belong to the class $B$.

Fix $n\in\mathbb N$ and $m\in\{1,\ldots,n\}$.
Let $\alpha_j>0$, $z_j\coloneqq e^{i\varphi_j}$, $j=1,\ldots,m$,
$0\leqslant\varphi_1<\ldots<\varphi_m<2\pi$, $f=e^{-h}$, $h(z)=\sum\limits_{j=1}^m\alpha_j\frac{1+z_jz}{1-z_jz}$.
Then $\{f\}_0=e^{-t}$ for some $t>0$.
The coefficients $\{f\}_\ell$, $\ell=1,\ldots,n$ are calculated successively by formula~\eqref{eq:Stupin:deForBCoeff}.
Define the polynomial $H$ by formula~\eqref{eq:Stupin:extrPolynomialH_n}.

Set $\mathfrak E_m(n)\coloneqq\mathfrak E(n)\cap\mathbb E_m(n)$.
If $f\in\mathfrak E_m(n)$ or $f\in\mathbb E_m(n)$, then $\operatorname{Re}H(e^{i\varphi})\geqslant0$,
$\varphi\in[0,2\pi)$, and $\operatorname{Re}H(z_j)=0$, $j=1,\ldots,m$.
Lemma~\ref{lem:Stupin:boundaryZeroDerivative} essentially gives the necessary condition $\operatorname{Im}\bigl(z_j
H'(z_j)\bigr)=0$, $j=1,\ldots,m$ for the function $f$ to belong to the class $\mathfrak E_m(n)$ and the class $\mathbb
E_m(n)$.
Replacing the condition $\operatorname{Re}H\geqslant0$, $|z|=1$ by this necessary condition slightly weakens the
conditions used earlier.
Including the equation for $t$ gives the following system with unknowns $t,\alpha_j,z_j$, $j=1,\ldots,m$:
\begin{equation}
    \label{eq:Stupin:conclusionSystem}
    \begin{cases}
        t=\sum\limits_{j=1}^m\alpha_j,\\
        \alpha_j>0,\quad |z_j|=1, &j=1,\ldots,m,\\
        z_j\neq z_k, &1\leqslant j<k\leqslant m,\\
        \{f\}_n>0,\\
        \operatorname{Re}H(z_j)=0,&j=1,\ldots,m,\\
        \operatorname{Im}\bigl(z_jH'(z_j)\bigr)=0,&j=1,\ldots,m.
    \end{cases}
\end{equation}
Thus every function $f\in\mathfrak E_m(n)$ and every function $f\in\mathbb E_m(n)$ corresponds to a solution of
system~\eqref{eq:Stupin:conclusionSystem}.
Every solution of this system determines a function in $B$.
Hence the set of functions determined by solutions of systems~\eqref{eq:Stupin:conclusionSystem} for $m=1,\ldots,n$
contains $\mathfrak E(n)$ and $\mathbb E(n)$, since $\mathfrak E(n)=\bigcup\limits_{m=1}^n\mathfrak E_m(n)$ and $\mathbb
E(n)=\bigcup\limits_{m=1}^n\mathbb E_m(n)$.

Make the substitution $z_j=x_j+iy_j$ and set $\widehat H\coloneqq e^t H=H/\{f\}_0$.
By formula~\eqref{eq:Stupin:deForBCoeff}, the ratios $\{f\}_\ell/\{f\}_0$, $\ell=0,\ldots,n$ are polynomials in
$\alpha_j$ and $z_j$, $j=1,\ldots,m$.
Thus, after this substitution, the real and imaginary parts of the coefficients $\widehat H$ become polynomials with
real coefficients in $\alpha_j,x_j,y_j$.
Using $e^t>0$ and $H(0)=\{f\}_n$, we obtain an equivalent system of algebraic equations and inequalities:
\begin{equation}
    \label{eq:Stupin:conclusionAlgebraicSystem}
    \begin{cases}
        t=\sum\limits_{j=1}^m\alpha_j,\\
        \alpha_j>0, \quad x_j,y_j\in\mathbb R,&j=1,\ldots,m,\\
        x_j^2+y_j^2=1,&j=1,\ldots,m,\\
        (x_j-x_k)^2+(y_j-y_k)^2>0,&1\leqslant j<k\leqslant m,\\
        \operatorname{Im}\widehat H(0)=0,\quad \operatorname{Re}\widehat H(0)>0,\\
        \operatorname{Re}\widehat H(x_j+iy_j)=0,&j=1,\ldots,m,\\
        \operatorname{Im}\bigl((x_j+iy_j)\widehat H'(x_j+iy_j)\bigr)=0,&j=1,\ldots,m.
    \end{cases}
\end{equation}
The numbering of the triples $(\alpha_j,x_j,y_j)$ does not affect the function $f$.
For every solution of system~\eqref{eq:Stupin:conclusionAlgebraicSystem}, set $x_j+iy_j=e^{i\varphi_j}$, where
$\varphi_j\in[0,2\pi)$, $j=1,\ldots,m$.
These triples can be renumbered so that
$0\leqslant\varphi_1<\ldots<\varphi_m<2\pi$.

The unknown $t$ is uniquely determined by the first equation and can be eliminated by substitution.
Consider the subsystem of algebraic equations of system~\eqref{eq:Stupin:conclusionAlgebraicSystem}.
If it can be reduced to one nontrivial algebraic equation in one unknown, and each of its roots corresponds to only
finitely many vectors of values of the remaining unknowns, then this subsystem has finitely many complex solutions.
We call this equation the determining equation.
Recall that above we obtained determining equations by eliminating unknowns.
Then among the solutions of this subsystem we select those satisfying all inequalities of
system~\eqref{eq:Stupin:conclusionAlgebraicSystem}.
The remaining solutions of the subsystem are discarded.

Note that the degree of the determining equation obtained determines the number of its complex roots counted with
multiplicity.
The number of real roots and their isolating intervals can be determined exactly.
This allows numerical calculations with control of completeness and error.
Comparing the values of $\{f\}_n$ for functions determined by all solutions of
systems~\eqref{eq:Stupin:conclusionAlgebraicSystem} for $m=1,\ldots,n$ solves the Krzyz problem for index $n$; proving
the Krzyz conjecture for index $n$ also requires analysis of all equality cases.
Finiteness of the solution set of the subsystem of algebraic equations and the possibility of eliminating unknowns for
arbitrary $m,n$ require a separate investigation.

\section*{Conclusion}

To prove the Krzyz conjecture, it is not enough to show that the inequality $|\{f\}_n|\leqslant\frac2e$ holds for all
$f\in B$ and all $n\in\mathbb N$; uniqueness of the global extremal must also be proved for every $n$.
Uniqueness of the global extremal may follow from a proof that $|\{f\}_n|\leqslant\frac2e$, but does not follow directly
from this fact alone.

On the other hand, uniqueness of the global extremal for some index $n$ implies (see
Theorem~\ref{th:Stupin:uniqueGlobalExtremalImpliesKrzyz}) that $|\{f\}_n|\leqslant\frac2e$.

For $n=1$ the Krzyz conjecture holds by Theorem~\ref{th:Stupin:Krzyz1}.
Let $n\in\mathbb N\setminus\{1\}$.
One promising approach to proving the conjecture for index $n$ is to prove uniqueness of the global extremal $f$ in
$\mathfrak E(n)$.
Invariance of the function $f$ under admissible rotations in the variable plane $z$ then implies
$\{f\}_1=\ldots=\{f\}_{n-1}=0$, after which Theorem~\ref{th:Stupin:converseToUniqunessInB} gives $f(z)=F(-z^n,1)$,
implying both the required estimation and the equality cases.

Another possibility is to prove that for every global extremal all roots of the polynomial $H$ lie on $\partial\Delta$.
Then $\{f\}_n=2\{f\}_0$, and Theorem~\ref{th:Stupin:uniqunessOfExtremalFunctionInB} gives $f(z)=F(-z^n,1)$.
However, membership of $f$ in $\mathbb E(n)$ does not imply that all roots of $H$ lie on $\partial\Delta$.

Systems~\eqref{eq:Stupin:conclusionSystem} for all $m=1,\ldots,n$ allow all functions in $\mathfrak E(n)$ to be
obtained.
Thus solving them completely for $m=1,\ldots,n$ and comparing all values $\{f\}_n$ can prove or disprove the conjecture
for fixed $n$.

Finally, by Theorem~\ref{th:Stupin:existenceInB}, a complete proof for $n>1$ is equivalent to the absence in the classes
$B_t$, $t\geqslant1$ of a function with coefficient $\{f\}_n=2/e$ and at least one nonzero coefficient among
$\{f\}_1,\ldots,\{f\}_{n-1}$.

\bigskip
\noindent\textbf{Information about the author}

\noindent Denis L. Stupin,
Candidate of Physics and Mathematics,
Associate Professor of the Fundamental Mathematics and Digital Technologies Department,
Tver State University, Tver, Russian Federation.

\noindent E-mail: \href{mailto:dstupin@mail.ru}{dstupin@mail.ru}

\noindent ORCID: \url{https://orcid.org/0000-0002-9183-9543}

\noindent There is no conflict of interests.
\end{document}